\documentclass[12pt,a4paper]{article}
\usepackage{multirow}
\usepackage{graphicx}
\usepackage{latexsym,amsmath,amsfonts,amscd,amssymb}
\usepackage{amsmath}
\usepackage{bm}
\usepackage{algorithm}
\usepackage{algorithmic}

\usepackage{epsfig}
\usepackage{epstopdf}
\usepackage{subfigure}
\usepackage{float}
\usepackage{changebar}
\usepackage{indentfirst}
\usepackage{verbatim}
\usepackage{color}
\usepackage{colordvi}
\allowdisplaybreaks[4]
\usepackage{appendix}
\usepackage{cite}


\newcommand{\bx}{\mathbf{x}}

\newtheorem{remark}{Remark}[section]

\begin{document}

\baselineskip=2pc
\begin{center}

{\Large \bf High order finite difference Hermite WENO fixed-point fast sweeping method for static Hamilton-Jacobi equations}
\end{center}
\centerline{Yupeng Ren\footnote{School of Mathematical Sciences, Xiamen University,
Xiamen, Fujian 361005, P.R. China. E-mail: ypren@stu.xmu.edu.cn. This work was carried out when Y. Ren was visiting Department of Mathematics, The Ohio State University.},
Yulong Xing\footnote{Department of Mathematics, The Ohio State University, Columbus, OH 43210, USA.
 E-mail: xing.205@osu.edu. The work of Y. Xing is partially supported by the NSF grant DMS-1753581.},
 and Jianxian Qiu\footnote{School of Mathematical Sciences and Fujian Provincial
Key Laboratory of Mathematical Modeling and High-Performance
Scientific Computing, Xiamen University,
Xiamen, Fujian 361005, P.R. China. E-mail: jxqiu@xmu.edu.cn. The work of J. Qiu is partially supported by NSFC grant 12071392.
}}

\baselineskip=1.8pc
\vspace{1cm}
\centerline{\bf Abstract}
\bigskip
In this paper, we combine the nonlinear HWENO reconstruction in \cite{newhwenozq} 
and the fixed-point iteration with Gauss-Seidel fast sweeping strategy, to solve the static Hamilton-Jacobi equations in a novel HWENO framework recently developed in \cite{mehweno1}. The proposed HWENO frameworks enjoys several advantages. First, compared with the traditional HWENO framework, the proposed methods do not need to introduce additional auxiliary equations to update the derivatives of the unknown function $\phi$. They are now computed from the current value of $\phi$ and the previous spatial derivatives of $\phi$. This approach saves the computational storage and CPU time, which greatly improves the computational efficiency of the traditional HWENO scheme. In addition, compared with the traditional WENO method, reconstruction stencil of the HWENO methods becomes more compact, their boundary treatment is simpler, and the numerical errors are smaller on the same mesh.
Second, the fixed-point fast sweeping method is used to update the numerical approximation. It is an explicit method and does not involve the inverse operation of nonlinear Hamiltonian,
therefore any Hamilton-Jacobi equations with complex Hamiltonian can be solved easily. It also resolves some known issues, including that the iterative number is very sensitive to the parameter $\varepsilon$ used in the nonlinear weights, as observed in previous studies.
Finally, in order to further reduce the computational cost, a hybrid strategy is also presented. Extensive numerical experiments are performed on two-dimensional problems, which demonstrate the good performance of the proposed fixed-point fast sweeping HWENO methods.


\section{Introduction}
\setcounter{equation}{0}
\setcounter{figure}{0}
\setcounter{table}{0}

The static Hamilton-Jacobi (HJ) equations often appear in many different application fields, for instance in optimal control, computer vision, differential game, geometric optics, image processing and so on \cite{Huang1,Yxia}. The general static HJ equations have the form
\begin{equation}\label{Y1.0}
\begin{cases}
H(\nabla\phi,\mathbf{x})=0,    & \mathbf{x}\in\Omega \setminus \Gamma,\\
\phi(\mathbf{x})=g(\mathbf{x}),                         & \mathbf{x}\in\Gamma\subset\Omega,
\end{cases}
\end{equation}
where $\Omega$ is the computational domain in $\mathbb{R}^{d}$,
$\phi(\bx)$ is the unknown function in $\Omega$, the Hamiltonian $H$ is a nonlinear Lipschitz continuous function depending on
$\nabla\phi$ and $\bx$, and the boundary condition is given by $\phi(\mathbf{x})=g(\mathbf{x})$ on the subset $\Gamma\subset\Omega$.
Eikonal equation is a prototype example of the static HJ equations, taking the form
\begin{equation}\label{eq:eikonal}
\begin{cases}
|\nabla\phi|=f(\mathbf{x}),    & \mathbf{x}\in\Omega \setminus \Gamma,\\
\phi(\mathbf{x})=g(\mathbf{x}),                         & \mathbf{x}\in\Gamma\subset\Omega ,
\end{cases}
\end{equation}
where $f(\mathbf{x})>0$. It can be derived from Maxwell's electromagnetic equations and provides a link between physical optics and geometric optics.

In general, the global $C^{1}$ solution does not exist for the time-dependent nonlinear HJ equations, even if the initial condition is sufficiently smooth. Singularities in the form of discontinuities would appear in the derivatives of the unknown function, hence it is necessary to define a ``weak solution'' for the HJ equations. The viscosity solutions of the HJ equations were first introduced by Crandall and Lions in \cite{viscosity}.

There are mainly two types of numerical methods to solve the static HJ equations. The first one is to solve the following equation
\[
  \phi_{t}+H(\nabla \phi)=f(\textbf{x}),\quad \mathbf{x}\in \mathbb{R}^{d},
\]
with pseudo-time iteration. The equation is evolved in time \cite{shutimehj} until the numerical solution converges to a steady state. However, such method requires a very large number of iterations to obtain the convergence of the solution in the entire domain, with the main reason being the finite speed of propagation, and the restrictive CFL time step requirement for stability. The second popular algorithm is to treat the problem as a stationary boundary value problem, so that the fast marching method (FMM) \cite{Tsitsiklis,Sethian,Helmsen} or the fast sweeping method (FSM) \cite{KaoLF,qianzhangzhao,Tsai,Zhaofsm1} can be applied. FSM can be constructed to be high order accurate, and becomes a class of popular and effective methods for solving static HJ equations nowadays.
The FSM was first proposed in \cite{Boue} by Bou\'e and Dupuis when solving a deterministic control problem with quadratic running cost using Markov chain approximation. Later, Zhao \cite{Zhaofsm1} studied the FSM for the Eikonal equation, and demonstrate the efficiency and effectiveness of the method.
Since then, many high order FSM have been proposed to solving static HJ equations. In the framework of finite difference methods, Zhang et al. \cite{zhang2006} proposed the third order weighted essentially non-oscillatory (WENO) FSM scheme \cite{WENOJP}, and its fifth order extension was studied in \cite{xiong}. High order accurate boundary treatments that are consistent with high order FSM, including Richardson extrapolation and Lax-Wendroff type procedure, have been developed for the inflow boundary conditions in \cite{huangnum,xiong}. In the setting of finite element methods, some researchers have combined high order FSM with the finite element discontinuous Galerkin (DG) method to solve Eikonal equation in  \cite{dgli,dgwu,dgzhang}, and their numerical performance were shown to be effective and robust.

Inspired by these successful methods, Zhang et al. \cite{fixhjzhang} proposed to embed FSM techniques into a time-marching scheme to solve static HJ equations, by introducing the fixed-point iterative schemes. These methods take the explicit form and do not have to involve any inverse operation of the nonlinear Hamiltonian. Hence, they are very flexible and can be used in solving very general static HJ equations with complex Hamiltonian. The fixed-point fast sweeping method has also been extended to solve the steady state hyperbolic conservation laws in \cite{fixhyperbolicwu,fixhyperbolicchen,fixhyperbolicli}.

Recently, high order Hermite WENO (HWENO) methods \cite{qiushuhwenohj2005,lqfdhweno1} have gained many attention in solving hyperbolic partial differential equations. Both the  WENO and HWENO methods can achieve the high order accuracy and preserve the essentially non-oscillatory property,
but the HWENO scheme uses the Hermite interpolation in reconstructing polynomials, that involves both the unknown variable $\phi$ and its first order spatial derivative (or first moment). Therefore, the HWENO reconstruction stencil becomes more compact and their boundary treatment is much simpler, although more storage and some additional work are needed to evaluate the spatial derivatives.
The HWENO scheme was first proposed in the construction of a suitable limiter for the DG method \cite{hwenolim1,hwenolim2} due to its compact stencil, and was later used to solve the time-dependent HJ equation in \cite{qiushuhwenohj2005}. Numerical results demonstrate that the HWENO scheme has smaller errors than the traditional WENO method with the same mesh and order. The HWENO scheme was also extended to solve the hyperbolic conservation law in the finite difference framework \cite{lqfdhweno1}.
Since then, many HWENO schemes have been developed to solve hyperbolic conservation laws \cite{taohwenosta,7thhweno,zhaohyhweno} and time-dependent HJ equations on structured \cite{zhengfdhwenohj}  and triangular meshes \cite{hwenounmesh,zhengfdhwenohjtr}. In addition, it was observed in \cite{HWENO_AP} that the finite volume HWENO scheme enjoys the asymptotic preserving property, when applied to the steady-state discrete ordinates ($S_N$) transport equation.
Recently, a new HWENO scheme (denoted by HWENO-ZZQ) was proposed by Zhu et al. in \cite{newhwenozq} for solving the time-dependent HJ equations. The main difference between the standard HWENO and HWENO-ZZQ methods lies in the stencil used in the reconstruction procedure. The HWENO-ZZQ reconstruction involves one big stencil and two small stencils, and their linear weights can be any positivity number summing up to 1. In additional to the simplicity, this method is shown in \cite{newhwenozq} to yield small errors with the same high order accuracy in the smooth areas, and maintain sharp transitions and non-oscillatory property near discontinuity.

In this paper, built upon the high order HWENO-ZZQ schemes and different time discretization schemes, we design a class of fixed-point methods with Gauss-Seidel fast sweeping strategy for efficiently solving the static HJ equations. It is worth mentioning that we are no longer using the traditional HWENO framework, namely, using one equation to update $\phi$ and several additional auxiliary equations to update its derivatives. Instead, we use the new approach presented in our recent work \cite{mehweno1}, to use only one equation to update $\phi$. The derivatives of $\phi$ will be obtained from applying the HWENO-ZZQ reconstruction on the updated values of $\phi$ and the previous values of the derivatives. This approach saves the computational storage and CPU costs, which improves the computational efficiency of the traditional HWENO scheme.
We will present four different sweeping methods to discrete the temporal derivatives, which include the forward Euler (FE), Runge-Kutta (RK) time-marching methods, as well as these two methods combined with Gauss-Seidel fast sweeping technique.
Extensive numerical examples suggest that the FE time discretization with Gauss-Seidel fast sweeping strategy is the most effective method for solving the static HJ equations. In addition, a hybrid strategy which combines both linear and HWENO-ZZQ reconstruction is also proposed and numerically validated to provide additional savings in computational time.

It was numerically observed in \cite{xiong,mehweno1} that, in some test cases, the convergence of the WENO/HWENO FSM scheme is very sensitive to the parameter $\varepsilon$ used in the nonlinear weights of HWENO reconstruction procedure. In order to address this issue, they adjust the parameter $\varepsilon$ according to the mesh size, so that the scheme can converge quickly and provide the expected high order. However, the manual adjustment of $\varepsilon$ will impact the application of the proposed scheme as the best choice of $\varepsilon$ cannot be known a priori and may be problem and mesh dependent.
An interesting observation is that the proposed HWENO fixed-point fast sweeping methods resolve this issue, and all of the numerical examples work well with a fixed parameter $\varepsilon$.
In addition, previous work \cite{me,mehweno1} suggested that the HWENO FSM does not converge to machine epsilon when the Godunov flux is used for some numerical examples, and requires additional work to address them. The numerical examples demonstrate that the HWENO fixed-point fast sweeping methods proposed in this paper work well for these tests.

The rest of the paper is organized as follows. In Section \ref{sec:method}, we present HWENO fixed-point fast sweeping methods and the hybrid strategy. The numerical tests are presented in Section \ref{sec:numerical} to demonstrate the effectiveness and efficiency of our schemes.  Conclusion remarks are provided in Section \ref{sec:conclusion}.

\section{High order fixed-point fast sweeping methods}\label{sec:method}
In this section, we present high order HWENO fixed-point fast sweeping methods to solve the static HJ equations. For ease of presentation, we consider the following two-dimensional static HJ equation
\begin{equation}\label{Y1.1}
\begin{cases}
H(\phi_{x},\phi_{y})=f(x,y),    & (x,y)\in\Omega \setminus \Gamma,\\
\phi(x,y)=g(x,y),                         & (x,y)\in\Gamma\subset\Omega.
\end{cases}
\end{equation}
Suppose the computational domain $\Omega$ is discretized into uniform rectangular meshes $\Omega_{h}=\{(x_i, y_j), 1\leq i\leq
N_x, 1\leq j\leq N_y\}$, with $(x_i, y_j)$ being a grid point in $\Omega_{h}$.
$\Delta{x}$ and $\Delta{y}$ stand for the grid sizes in the $x$ and $y$ directions, respectively.
We denote the numerical solution at the grid point $(x_i, y_j)$ by $\phi_{i,j}$. In addition, we introduce the notations $u=\phi_{x}(x, y)$ and $v=\phi_{y}(x, y)$ as the first order partial derivatives of $\phi$ with respect to the variables $x$ and $y$, respectively.

\subsection{Numerical Hamiltonian}\label{sec:approach0}

The numerical approximation of the Hamiltonian in \eqref{Y1.1} is given by a monotone numerical Hamiltonian $\widehat{H}$:
\begin{equation}\label{hamiltonian}
\left.H(\phi_{x},\phi_{y})\right|_{ij} \approx \widehat{H}(\phi_{x}^{-},\phi_{x}^{+},\phi_{y}^{-},\phi_{y}^{+})_{ij} .
\end{equation}
Such numerical Hamiltonian takes inputs $\phi^{\pm}_x$ and $\phi^{\pm}_y$ at the corresponding
grid point, which will be reconstructed from its neighboring point values using the high order reconstruction procedure as detailed in the next subsection.

Two types of numerical Hamiltonian are often considered for Hamilton-Jacobi equations.
For general static HJ equations, we adopt the Lax-Friedrichs (LF) numerical Hamiltonian\cite{lfflux}
\begin{equation}\label{lfflux}
  \widehat{H}^{LF}_{i,j}=H\left(\frac{u_{i,j}^{-}+u_{i,j}^{+}}{2},\frac{v_{i,j}^{-}+v_{i,j}^{+}}{2}\right)-
\frac{1}{2}\alpha(u_{i,j}^{+}-u_{i,j}^{-})-
\frac{1}{2}\beta(v_{i,j}^{+}-v_{i,j}^{-}),
\end{equation}
where
\begin{equation}\label{alphabeta}
\alpha=\max_{u,v}|H_{1}(u,v)|, \qquad \beta=\max_{u,v}|H_{2}(u,v)|,
\end{equation}
and $H_{\ell}(u, v)\, (\ell=1,2)$ denotes the partial derivative of $H$ with respect to the $\ell$-th argument.

The other commonly used numerical Hamiltonian is the Godunov numerical Hamiltonian, often employed in the approximation of the
Eikonal equation. As discussed in \cite{Oshershu1991}, the Godunov numerical Hamiltonian for general Hamiltonian $H(u,v)$ takes the form
\begin{equation} \label{eq:generalGodunov}
  \widehat{H}^{G}(u^{-},u^{+},v^{-},v^{+})=ext_{u\in I(u^{-},u^{+})}ext_{v\in I(v^{-},v^{+})}H(u,v),
\end{equation}
where $I(a,b)=[\min(a,b),\max(a,b)]$ and
the function \emph{ext} is defined by
\begin{equation*}
  ext_{u\in I(a,b)}= \begin{cases}
  \min\limits_{a\leq u\leq b}, \quad if~a\leq b,\\
    \max\limits_{b\leq u\leq a}, \quad if~a\geq b.
  \end{cases}
\end{equation*}
For the two-dimensional Eikonal equation
\begin{equation}\label{eikonal}
  \begin{cases}
 \displaystyle \sqrt{\phi_{x}^2+\phi_{y}^{2}}=f(x,y),&(x,y)\in\Omega,\\
 \displaystyle \phi(x,y)=g(x,y),&(x,y)\in\Gamma\subset\Omega,
  \end{cases}
\end{equation}
the Godunov numerical Hamiltonian \eqref{eq:generalGodunov} reduces to
\begin{equation}\label{godunovflux}
\widehat{H}^{G}(u^{-},u^{+},v^{-},v^{+})=\sqrt{\max\{(u^{-})^{+},(u^{+})^{-}\}^{2}+\max\{(v^{-})^{+},(v^{+})^{-}\}^{2}},
\end{equation}
where $x^{+}=\max(x,0)$ and $x^{-}=-\min(x,0)$.

\subsection{The review of HWENO-ZZQ reconstruction}\label{sec:approach1}

In this subsection, we will review the finite difference HWENO-ZZQ reconstruction recently developed in \cite{newhwenozq}.
To save space, we only illustrate the reconstruction of $(\phi_{x})^{\pm}_{i,j}$ along $x$-direction here. The approximation of $(\phi_{y})^{\pm}_{i,j}$ along $y$-direction can be obtained similarly, and we refer to \cite{newhwenozq} for the details. In this paper, we focus on the fifth order finite difference HWENO-ZZQ reconstruction.
\begin{itemize}
  \item Reconstruction of $(\phi_{x})_{i,j}^{-}$ from the upwind information: \\
Take a big stencil $S_{0}=\{x_{i-2}, x_{i-1}, x_{i}, x_{i+1}\}$ and two small stencils $S_{1}=\{x_{i-2},  x_{i-1}, x_{i}\}$, $S_{2}=\{x_{i-1}, x_{i}, ,x_{i+1}\}$, we compose a Hermite quintic polynomial $p^{-}_{1}(x)$, and two quadratic polynomials $p^{-}_{2}(x)$, $p^{-}_{3}(x)$ satisfying
\begin{equation*}
\begin{split}
 p^{-}_{1}(x_{k})&=\phi_{k,j},\quad k=i-2,\cdots,i+1,\,\,\mathrm{ and }\,\,~(p^{-}_{1})'|_{x_{k}}=u_{k,j},\quad k=i-1,i+1;\\
p^{-}_{2}(x_{k})&=\phi_{k,j},\quad k=i-2,i-1,i;\,\,\quad \mathrm{ and }\,\,~p^{-}_{3}(x_{k})=\phi_{k,j},\quad k=i-1,i,i+1;
\end{split}
\end{equation*}
The values of their first-order derivative at $x=x_i$ can be evaluated as
\begin{subequations}
\begin{align}
&(\phi_{x})_{i,j}^{-,1}=(p^{-}_{1})'|_{x_{i}}=\frac{\phi_{i-2,j}+18\phi_{i-1,j}-9\phi_{i,j}-10\phi_{i+1,j}+
9\Delta x u_{i-1,j}+3\Delta x u_{i+1,j}}{-18\Delta x};	\label{hybridlinearf}\\
&(\phi_{x})_{i,j}^{-,2}=(p^{-}_{2})'|_{x_{i}}=\frac{\phi_{i-2,j}-4\phi_{i-1,j}+3\phi_{i,j}}{2\Delta x};\\	
 &(\phi_{x})_{i,j}^{-,3}=(p^{-}_{3})'|_{x_{i}}=\frac{-\phi_{i-1,j}+\phi_{i+1,j}}{2\Delta x}.
\end{align}
\end{subequations}
  \item Reconstruction of $(\phi_{x})_{i,j}^{+}$ from the downwind information: \\
Take a big stencil $\widetilde{S}_{0}=\{x_{i-1}, x_{i}, x_{i+1},x_{i+2}\}$ and two small stencils $\widetilde{S}_{1}=\{x_{i-1}, x_{i},x_{i+1}\}$, $\widetilde{S}_{2}=\{x_{i}, x_{i+1}, ,x_{i+2}\}$, we compose a Hermite quintic polynomial $p^{+}_{1}(x)$, and two quadratic polynomials $p^{+}_{2}(x)$, $p^{+}_{3}(x)$ such that
\begin{equation*}
\begin{split}
 p^{+}_{1}(x_{k})&=\phi_{k,j},\quad k=i-1,\cdots,i+2, \,\,\mathrm{and}\,\,~(p^{+}_{1})'|_{x_{k}}=u_{k,j},\quad k=i-1,i+1;\\
p^{+}_{2}(x_{k})&=\phi_{k,j},\quad k=i-1,i,i+1;\,\,\quad \mathrm{ and }\,\,~p^{+}_{3}(x_{k})=\phi_{k,j},\quad k=i,i+1,i+2;
\end{split}
\end{equation*}
The values of their first-order derivative at $x=x_i$ can be evaluated as
\begin{subequations}
\begin{align}
&(\phi_{x})_{i,j}^{+,1}=(p^{+}_{1})'|_{x_{i}}=\frac{10\phi_{i-1,j}+9\phi_{i,j}-18\phi_{i+1,j}-\phi_{i+2,j}+3\Delta xu_{i-1,j}+9\Delta xu_{i+1,j}}{-18\Delta x};	\label{hybridlinearz}\\
&(\phi_{x})_{i,j}^{+,2}=(p^{+}_{2})'|_{x_{i}}=\frac{\phi_{i+1,j}-\phi_{i-1,j}}{2\Delta x};	\\
 &(\phi_{x})_{i,j}^{+,3}=(p^{+}_{3})'|_{x_{i}}=\frac{-3\phi_{i,j}+4\phi_{i+1,j}-\phi_{i+2,j}}{2\Delta x}.
\end{align}
\end{subequations}
\end{itemize}

In the nonlinear HWENO reconstructions, $(\phi_{x})_{i,j}^{\pm}$ are computed as a convex combination of these three corresponding values \cite{Levy1,newhwenozq}
\begin{equation}\label{weno}
(\phi_{x})_{i,j}^{\pm}=\omega_{1}^{\pm}\left(\frac{1}{\gamma_{1}}(\phi_{x})_{i,j}^{\pm,1}-
\frac{\gamma_{2}}{\gamma_{1}}(\phi_{x})_{i,j}^{\pm,2}-\frac{\gamma_{3}}{\gamma_{1}}(\phi_{x})_{i,j}^{\pm,3}\right)
+\omega_{2}^{\pm}(\phi_{x})_{i,j}^{\pm,2}+\omega_{3}^{\pm}(\phi_{x})_{i,j}^{\pm,3},
\end{equation}
where the parameters $\omega_n~(n=1,2,3)$ are called the nonlinear weights. The parameters $\gamma_n $ could be any positive constants satisfying $\gamma_1+\gamma_2+\gamma_3=1$. The nonlinear weights $\omega_{n}$ can be computed from
\begin{equation}\label{epsilon}
\omega_{n}^{\pm}=\frac{\overline{\omega}_{n}^{\pm}}{\sum_{l=1}^{3}\overline{\omega}_{l}^{\pm}},\qquad
\overline{\omega}_{n}=\gamma_{n}\left(1+\frac{\tau^{\pm}}{\varepsilon+\beta_{n}^{\pm}}\right),\qquad n=1,2,3,
\end{equation}
where $\varepsilon$ is a small positive number to avoid the denominator from becoming $0$, and
$$\tau^{\pm}=\left(\frac{|\beta_{1}^{\pm}-\beta_{2}^{\pm}|+|\beta_{1}^{\pm}-\beta_{3}^{\pm}|}{2}\right)^{2},$$
with $\beta_{n}^{\pm}$ being the so-called smoothness indicators
$$\beta^{\pm}_{n}=
\sum_{\alpha=2}^{r}\int^{x_{i+\frac{1}{2}}}_{x_{i-\frac{1}{2}}}\Delta x^{2\alpha-3}\left(\frac{d^{\alpha}p^{\pm}_{n}(x)}{dx^{\alpha}}\right)^2dx,
\quad n=1,2,3,
$$
measuring the smoothness of the derivative functions of $p^{\pm}_{n}(x)$ near the target point $x_{i}$. The parameter $r$ is the degree of polynomial $p_{n}^{\pm}(x)$, and here we set $r=5$ for $n=1$, and $r=2$ for $n=2,3$, respectively.

\subsection{HWENO fixed-point fast sweeping methods}\label{fourschemes}

\indent After the spatial HWENO discretization, the model \eqref{Y1.1} reduces to a large nonlinear system of algebraic equations, and its size is determined by the number of spatial grid points. There are different approaches available to solve this nonlinear system. Here, we propose to use the fixed-point fast sweeping methods, and discuss four different sweeping methods below.

{\bf Forward Euler Jacobi method:} We could view the static HJ equation as the steady state version of the time-dependent HJ equation. For the time-dependent problem, the simple FE method could be used for the temporal discretization, and the resulting explicit time marching scheme can be written as follows
\begin{equation}\label{preiterative}
\phi_{i,j}^{n+1}=\phi_{i,j}^{n}+\Delta t \left[f_{i,j}-\widehat{H}(\phi_{x}^{-},\phi_{x}^{+},\phi_{y}^{-},\phi_{y}^{+})^n_{ij}\right],
\end{equation}
where the LF numerical Hamiltonian \eqref{lfflux} or Godunov numerical Hamiltonian \eqref{godunovflux} is used, and the time step size is given by
$$\Delta t=\gamma\left(\frac{1}{\frac{\alpha}{\Delta x}+\frac{\beta}{\Delta y}}\right),$$
with $\gamma$ being the CFL number, and $\alpha$, $\beta$ defined in \eqref{alphabeta}. Obviously, $\alpha=\beta=1$ for the Eikonal equation \eqref{eikonal}.
Here $\phi_{x}^{\pm},\phi_{y}^{\pm}$ are obtained through the high order HWENO reconstruction procedure based on the values of $\phi$, $u$, $v$. We introduce the following operator
$$
L\Big(\{\phi^{n}\}_{i-2,j}^{i+2,j},u_{i\pm1,j}^{n};\{\phi^{n}\}_{i,j-2}^{i,j+2},v_{i,j\pm1}^{n}\Big) := f_{i,j}-\widehat{H}(\phi_{x}^{-},\phi_{x}^{+},\phi_{y}^{-},\phi_{y}^{+})^n_{ij}
$$
to represent this dependence.
From the perspective of iterative schemes, the scheme \eqref{preiterative} can be considered as a forward Euler Jacobi (FE-Jacobi) type fixed-point iterative scheme
\begin{equation}\label{fejacobi}
\phi_{i,j}^{new}=\phi_{i,j}^{old}+\frac{\gamma}{\frac{\alpha}{\Delta x}+\frac{\beta}{\Delta y}}L\Big(\{\phi^{old}\}_{i-2,j}^{i+2,j},u_{i\pm1,j}^{old};\{\phi^{old}\}_{i,j-2}^{i,j+2},v_{i,j\pm1}^{old}\Big),
\end{equation}
where $\phi_{i,j}^{new}$ denotes the updated numerical approximations of $\phi$ at the grid point
$(x_i, y_j )$, $\phi_{i,j}^{old}$ denotes the previous value of $\phi$, and $u_{i,j}^{old}$, $v_{i,j}^{old}$ denote the previous values of $u$ and $v$ at the same grid point.

The values of $\phi$, $u$ and $v$ are used in the HWENO-ZZQ reconstruction.
In the traditional HWENO framework, one needs to take the spatial derivatives of \eqref{Y1.1}, to obtain two new equations involving $u$ and $v$, which will be used to update $u^{new}$ and $v^{new}$. In our recent work \cite{mehweno1}, we presented a novel HWENO method for static HJ equations to simplify this procedure, which will be used in this paper. Compared with traditional HWENO framework, the new method does not involve any additional auxiliary equations, and is more computationally efficient.
The main idea is to reuse the step of reconstructing $\phi_{x}^{\pm}$ and $\phi_{y}^{\pm}$ from $\phi$, $u$ and $v$, which was already available during the HWENO-ZZQ reconstruction procedure. More specifically, we use the updated $\phi^{new}$ (computed by \eqref{fejacobi} with suitable numerical Hamiltonian) and the previous spatial derivatives $u^{old}$, $v^{old}$ to reconstruct $\phi_{x}^{\pm}$ and $\phi_{y}^{\pm}$ by HWENO-ZZQ reconstruction, and then define $u^{new}$ and $v^{new}$ as
\begin{equation}\label{uvupdate}
  u_{i,j}^{new}=\begin{cases}
  (\phi_{x})_{i,j}^{-}, &\text{if}~ (\phi_{x})_{i,j}^{\pm}>0,\\
  (\phi_{x})_{i,j}^{+}, &\text{if}~ (\phi_{x})_{i,j}^{\pm}<0,\\
  u_{i,j}^{old},& \text{otherwise},
  \end{cases}
  \qquad
       v_{i,j}^{new}=\begin{cases}
  (\phi_{y})_{i,j}^{-}, &\text{if}~ (\phi_{y})_{i,j}^{\pm}>0,\\
  (\phi_{y})_{i,j}^{+}, &\text{if}~ (\phi_{y})_{i,j}^{\pm}<0,\\
  v_{i,j}^{old},& \text{otherwise}.
  \end{cases}
\end{equation}
Here $u_{i,j}^{new}$ and $v_{i,j}^{new}$ denote the updated numerical approximations of $u$ and $v$ at the grid point
$(x_i, y_j )$, respectively.

This finishes the description of the FE-Jacobi iterative scheme. The pseudo code of this method is presented in Algorithm \ref{fejacobipseudo}.
\begin{algorithm}
\caption{FE-Jacobi scheme}
\label{fejacobipseudo}
\begin{algorithmic}
\FOR {$i=1$ \TO $N_{x}$, $j=1$ \TO $N_{y}$}
\STATE $\phi_{i,j}^{new}=\phi_{i,j}^{old}+\gamma\left(\frac{1}{\frac{\alpha}{\Delta_{x}}+\frac{\beta}{\Delta_{y}}}\right)
L\Big(\{\phi^{old}\}_{i-2,j}^{i+2,j},u_{i\pm1,j}^{old};\{\phi^{old}\}_{i,j-2}^{i,j+2},v_{i,j\pm1}^{old}\Big)$,
\ENDFOR
\FOR {$i=1$ \TO $N_{x}$, $j=1$ \TO $N_{y}$}
\STATE Update $u_{i,j}^{new}$ and $v_{i,j}^{new}$ by \eqref{uvupdate} from $\phi_{i,j}^{new}$, $u^{old}_{i,j}$ and $v^{old}_{i,j}$.
\ENDFOR
\end{algorithmic}
\end{algorithm}

{\bf Forward Euler fast sweeping method:}
While the simple FE-Jacobi method converges for most of the problems, it suffers from linear stability problems when coupled with high-order spatial discretization, and requires a small CFL number $\gamma$, hence a lot of iteration steps, to converge. In order to accelerate its convergence under the simple explicit framework of fixed-point iteration \eqref{fejacobi},
the Gauss-Seidel (GS) sweeping strategy could be applied to scheme \eqref{fejacobi}. According to the GS philosophy, the newest available numerical values of $\phi$ are always used in the interpolation stencils as long as they are available, and the GS iterative scheme takes the form
\begin{equation}\label{fegauss}
\phi_{i,j}^{new}=\phi_{i,j}^{old}+\frac{\gamma}{\frac{\alpha}{\Delta x}+\frac{\beta}{\Delta y}}L\Big(\{\phi^{*}\}_{i-2,j}^{i+2,j},u_{i\pm1,j}^{*};\{\phi^{*}\}_{i,j-2}^{i,j+2},v_{i,j\pm1}^{*}\Big),
\end{equation}
where $\phi^{*}_{i,j}$ represents the most up-to-date point values of $\phi$ at the point $(x_i,y_j)$. Here, we further propose to combine the GS iteration with the fast sweeping idea, and proceed the sweeping in the following four alternating directions repeatedly
$$(1)\quad i=1:N_{x}, ~j=1:N_{y}; \qquad (2)\quad i=N_{x}:1,~ j=1:N_{y};$$
$$(3)\quad i=N_{x}:1,~ j=N_{y}:1; \qquad (4)\quad i=1:N_{x},~ j=N_{y}:1,$$
which leads to the FE type fixed-point fast sweeping method (FE-FSM). This method allows a larger CFL number and can reduce the number of iterations significantly. As observed in the FSM applied to solve the static HJ equation, this strategy utilized the directions of characteristic lines, and leads to an acceleration of the convergence speed significantly, as observed in our numerical experiment section.
The pseudo codes for FE-FSM are presented in Algorithm \ref{fegspseudo}.
\begin{algorithm}
\caption{FE-FSM}
\label{fegspseudo}
\begin{algorithmic}
\FOR {$i=1$ \TO $N_{x}$, $j=1$ \TO $N_{y}$}
\STATE $\phi_{i,j}^{new}=\phi_{i,j}^{old}+\gamma\left(\frac{1}{\frac{\alpha}{\Delta_{x}}+\frac{\beta}{\Delta_{y}}}\right)
L\Big(\{\phi^{*}\}_{i-2,j}^{i+2,j},u_{i\pm1,j}^{*};\{\phi^{*}\}_{i,j-2}^{i,j+2},v_{i,j\pm1}^{*}\Big)$,
\STATE update $u_{i,j}^{new}$ and $v_{i,j}^{new}$ by \eqref{uvupdate} from $\phi_{i,j}^{new}$, $u^{old}_{i,j}$ and $v^{old}_{i,j}$.
\ENDFOR
\STATE Repeat the above process in the other three sweeping directions.
\end{algorithmic}
\end{algorithm}

{\bf Runge-Kutta Jacobi method:}  We only considered the FE time discretization scheme so far. In the temporal discretization of differential equations, the high order RK methods are frequently used. Combining them with the Jacobi iteration, we have the third order RK Jacobi (RK-Jacobi) type fixed-point method of the form
\begin{subequations}\label{rkjacobi}
\begin{align}
&\phi_{i,j}^{(1)}=\phi_{i,j}^{old}+\Delta t L\Big(\{\phi^{old}\}_{i-2,j}^{i+2,j},u_{i\pm1,j}^{old};\{\phi^{old}\}_{i,j-2}^{i,j+2},v_{i,j\pm1}^{old}\Big),\\
\displaystyle&\phi_{i,j}^{(2)}=\frac{3}{4}\phi_{i,j}^{old}+\frac{1}{4}\phi_{i,j}^{(1)}+\frac{1}{4}\Delta t L\Big(\{\phi^{(1)}\}_{i-2,j}^{i+2,j},u_{i\pm1,j}^{(1)};\{\phi^{(1)}\}_{i,j-2}^{i,j+2},v_{i,j\pm1}^{(1)}\Big),\\
\displaystyle&\phi_{i,j}^{new}=\frac{1}{3}\phi_{i,j}^{old}+\frac{2}{3}\phi_{i,j}^{(2)}+\frac{2}{3}\Delta t L\Big(\{\phi^{(2)}\}_{i-2,j}^{i+2,j},u_{i\pm1,j}^{(2)};\{\phi^{(2)}\}_{i,j-2}^{i,j+2},v_{i,j\pm1}^{(2)}\Big).
\end{align}
 \end{subequations}

{\bf Runge-Kutta fast sweeping method:} Similar to FE-FSM, by combining the RK fixed-point method and Gauss-Seidel fast sweeping strategy, we can obtain the third order RK type fixed-point fast sweeping method (RK-FSM), which takes the form
 \begin{subequations}\label{rkgauss}
  \begin{align}
&\phi_{i,j}^{(1)}=\phi_{i,j}^{old}+\Delta t L\Big(\{\phi^{*}\}_{i-2,j}^{i+2,j},u_{i\pm1,j}^{*};\{\phi^{*}\}_{i,j-2}^{i,j+2},v_{i,j\pm1}^{*}\Big),\\
\displaystyle&\phi_{i,j}^{(2)}=\phi_{i,j}^{(1)}+\frac{1}{4}\Delta t L\Big(\{\phi^{*}\}_{i-2,j}^{i+2,j},u_{i\pm1,j}^{*};\{\phi^{*}\}_{i,j-2}^{i,j+2},v_{i,j\pm1}^{*}\Big),\\
\displaystyle&\phi_{i,j}^{new}=\phi_{i,j}^{(2)}+\frac{2}{3}\Delta t L\Big(\{\phi^{*}\}_{i-2,j}^{i+2,j},u_{i\pm1,j}^{*};\{\phi^{*}\}_{i,j-2}^{i,j+2},v_{i,j\pm1}^{*}\Big).
\end{align}
\end{subequations}
Note that this method contains three sub-iterations, and again, each iteration proceeds the sweeping in four alternating directions repeatedly. In the implementation, the sweeping directions of these three sub-iterations are the same in each iteration step. The pseudo codes of RK-Jacobi and RK-FSM are given in Algorithm \ref{rkjacobipseudo} and Algorithm \ref{rkgspseudo}, respectively.
\begin{algorithm}
\caption{RK-Jacobi scheme}
\label{rkjacobipseudo}
\begin{algorithmic}
\FOR {$i=1$ \TO $N_{x}$, $j=1$ \TO $N_{y}$}
\STATE $\phi_{i,j}^{(1)}=\phi_{i,j}^{old}+\gamma\left(\frac{1}{\frac{\alpha}{\Delta_{x}}+\frac{\beta}{\Delta_{y}}}\right) L\Big(\{\phi^{old}\}_{i-2,j}^{i+2,j},u_{i\pm1,j}^{old};\{\phi^{old}\}_{i,j-2}^{i,j+2},v_{i,j\pm1}^{old}\Big)$,
\ENDFOR
\FOR {$i=1$ \TO $N_{x}$, $j=1$ \TO $N_{y}$}
\STATE update $u_{i,j}^{(1)}$ and $v_{i,j}^{(1)}$ from $\phi_{i,j}^{(1)}$, $u^{old}_{i,j}$, $v^{old}_{i,j}$ by \eqref{uvupdate}.
\ENDFOR
\FOR {$i=1$ \TO $N_{x}$, $j=1$ \TO $N_{y}$}
\STATE $\phi_{i,j}^{(2)}=\frac{3}{4}\phi_{i,j}^{old}+\frac{1}{4}\phi_{i,j}^{(1)}+\frac{1}{4}\gamma\left(\frac{1}{\frac{\alpha}{\Delta_{x}}+\frac{\beta}{\Delta_{y}}}\right) L\Big(\{\phi^{(1)}\}_{i-2,j}^{i+2,j},u_{i\pm1,j}^{(1)};\{\phi^{(1)}\}_{i,j-2}^{i,j+2},v_{i,j\pm1}^{(1)}\Big)$,
\ENDFOR
\FOR {$i=1$ \TO $N_{x}$, $j=1$ \TO $N_{y}$}
\STATE update $u_{i,j}^{(2)}$ and $v_{i,j}^{(2)}$ from $\phi_{i,j}^{(2)}$, $u^{(1)}_{i,j}$, $v^{(1)}_{i,j}$ by \eqref{uvupdate}.
\ENDFOR
\FOR {$i=1$ \TO $N_{x}$, $j=1$ \TO $N_{y}$}
\STATE $\phi_{i,j}^{new}=\frac{1}{3}\phi_{i,j}^{old}+\frac{2}{3}\phi_{i,j}^{(2)}+\frac{2}{3}\gamma\left(\frac{1}{\frac{\alpha}{\Delta_{x}}+\frac{\beta}{\Delta_{y}}}\right) L\Big(\{\phi^{(2)}\}_{i-2,j}^{i+2,j},u_{i\pm1,j}^{(2)};\{\phi^{(2)}\}_{i,j-2}^{i,j+2},v_{i,j\pm1}^{(2)}\Big),$
\ENDFOR
\FOR {$i=1$ \TO $N_{x}$, $j=1$ \TO $N_{y}$}
\STATE update $u_{i,j}^{new}$ and $v_{i,j}^{new}$ from $\phi_{i,j}^{new}$, $u^{(2)}_{i,j}$, $v^{(2)}_{i,j}$ by \eqref{uvupdate}.
\ENDFOR
\end{algorithmic}
\end{algorithm}

\begin{algorithm}
\caption{RK-FSM} 
\label{rkgspseudo}
\begin{algorithmic}
\FOR {$i=1$ \TO $N_{x}$, $j=1$ \TO $N_{y}$}
\STATE $\phi_{i,j}^{(1)}=\phi_{i,j}^{old}+\gamma\left(\frac{1}{\frac{\alpha}{\Delta_{x}}+\frac{\beta}{\Delta_{y}}}\right) L\Big(\{\phi^{*}\}_{i-2,j}^{i+2,j},u_{i\pm1,j}^{*};\{\phi^{*}\}_{i,j-2}^{i,j+2},v_{i,j\pm1}^{*}\Big)$,
\STATE update $u_{i,j}^{(1)}$ and $v_{i,j}^{(1)}$ from $\phi_{i,j}^{(1)}$, $u^{old}_{i,j}$, $v^{old}_{i,j}$ by \eqref{uvupdate}.
\ENDFOR
\FOR {$i=1$ \TO $N_{x}$, $j=1$ \TO $N_{y}$}
\STATE $\phi_{i,j}^{(2)}=\phi_{i,j}^{(1)}+\frac{1}{4}\gamma\left(\frac{1}{\frac{\alpha}{\Delta_{x}}+\frac{\beta}{\Delta_{y}}}\right) L\Big(\{\phi^{*}\}_{i-2,j}^{i+2,j},u_{i\pm1,j}^{*};\{\phi^{*}\}_{i,j-2}^{i,j+2},v_{i,j\pm1}^{*}\Big)$,
\STATE update $u_{i,j}^{(2)}$ and $v_{i,j}^{(2)}$ from $\phi_{i,j}^{(2)}$, $u^{(1)}_{i,j}$, $v^{(1)}_{i,j}$ by \eqref{uvupdate}.
\ENDFOR
\FOR {$i=1$ \TO $N_{x}$, $j=1$ \TO $N_{y}$}
\STATE $\phi_{i,j}^{new}=\phi_{i,j}^{(2)}+\frac{2}{3}\gamma\left(\frac{1}{\frac{\alpha}{\Delta_{x}}+\frac{\beta}{\Delta_{y}}}\right) L\Big(\{\phi^{*}\}_{i-2,j}^{i+2,j},u_{i\pm1,j}^{*};\{\phi^{*}\}_{i,j-2}^{i,j+2},v_{i,j\pm1}^{*}\Big),$
\STATE update $u_{i,j}^{new}$ and $v_{i,j}^{new}$ from $\phi_{i,j}^{new}$, $u^{(2)}_{i,j}$, $v^{(2)}_{i,j}$ by \eqref{uvupdate}.
\ENDFOR
\STATE Repeat the above process in the other three sweeping directions.
\end{algorithmic}
\end{algorithm}

\begin{remark}
There are several different forms of the RK-FSM in addition to \eqref{rkgauss}, since we could mix the usage of $\phi^{old}$, $\phi^{(1)}$ and $\phi^{(2)}$ for this steady state problem. 
Numerically, we also tested other RK type fixed-point fast sweeping methods, including the following RK3 fixed-point fast sweeping scheme
 \begin{subequations}\label{RK-FSM-T}
  \begin{align}
&\phi_{i,j}^{(1)}=\phi_{i,j}^{old}+\Delta t L\Big(\{\phi^{*}\}_{i-2,j}^{i+2,j},u_{i\pm1,j}^{*};\{\phi^{*}\}_{i,j-2}^{i,j+2},v_{i,j\pm1}^{*}\Big),\\
\displaystyle&\phi_{i,j}^{(2)}=\frac{3}{4}\phi_{i,j}^{old}+\frac{1}{4}\phi_{i,j}^{(1)}+\frac{1}{4}\Delta t L\Big(\{\phi^{*}\}_{i-2,j}^{i+2,j},u_{i\pm1,j}^{*};\{\phi^{*}\}_{i,j-2}^{i,j+2},v_{i,j\pm1}^{*}\Big),\\
\displaystyle&\phi_{i,j}^{new}=\frac{1}{3}\phi_{i,j}^{old}+\frac{2}{3}\phi_{i,j}^{(2)}+\frac{2}{3}\Delta t L\Big(\{\phi^{*}\}_{i-2,j}^{i+2,j},u_{i\pm1,j}^{*};\{\phi^{*}\}_{i,j-2}^{i,j+2},v_{i,j\pm1}^{*}\Big).
\end{align}
\end{subequations}
The numerical experiments suggest that the scheme in \eqref{rkgauss} yields the fastest convergence among all these RK type fixed-point fast sweeping schemes.
\end{remark}

\subsection{The flowchart of four HWENO fixed-point sweeping schemes}\label{sec:flowchart}

We have proposed four HWENO fixed-point sweeping approaches to solve the static HJ equations. Here we will summarize the detailed procedure of these approaches, and provide a flowchart for them. We start by labeling the computational nodal points $\{(x_i, y_j)\}$ into several categories as in \cite{me}:

\noindent\emph{Category I}: For points on the boundary $\Gamma$, the exact values are assigned for these points.

\noindent\emph{Category II}: For ghost points (exterior of the boundary), we use the high order extrapolation to compute their numerical
solution $\phi_{i,j}$.

\noindent\emph{Category III}: For points near the $\Gamma$ (whose distances to $\Gamma$ are less than or equal to $2h$), the numerical boundary treatment from \cite{xiong,huangnum} could be used (i.e., Richardson extrapolation for a single point or a set of isolated points, and Lax-Wendroff type procedure for continuous $\Gamma$).
Since our main focus is on HWENO fixed-point sweeping method itself in this paper, the exact solutions are used on these points.

\noindent\emph{Category IV}: All the remaining points, which will be updated by fixed-point sweeping methods.

Note that only \emph{Category IV} points need to be updated by following numerical method. We now summarize our flowchart of these methods as follows:

\noindent \textbf{Step 1}. \emph{Initialization}:
The numerical solution from the first order fast sweeping method \cite{Zhaofsm1} is taken as the initial guess of $\phi$. The forward or backward difference of the resulting $\phi$ is used as the initial guess of $u$ and $v$.

\noindent\textbf{Step 2}. \emph{Update $\phi^{new}$}:
We can choose one of these four schemes \eqref{fejacobi},\eqref{fegauss},\eqref{rkjacobi} and \eqref{rkgauss} to update $\phi^{new}_{i,j}$ at the grid point $(x_{i},y_{j})$. The pseudo codes of these schemes have been given in Section \ref{fourschemes}.
The values at ghost points will be updated by high order extrapolations in the iterative methods.

\noindent\textbf{Step 3}. \emph{Convergence}: In general, the iteration will stop if, for two consecutive iteration steps, the error satisfies
$$\delta=||\phi^{new}-\phi^{old}||_{L_{1}}<10^{-14}.$$
\subsection{Comments and remarks}\label{sec:comments}
We need to emphasize that the local solver of fixed-point sweeping methods derived from different time marching schemes with either Jacobi or Gauss-Seidel fast sweeping strategy, are given in the explicit form and do not involve solving nonlinear equations. In principle, the approaches can be applied to any general static HJ equations with complicated Hamiltonian. If the Godunov Hamiltonian is used, the numerical method is different from the fast sweeping methods discussed in \cite{KaoLF,zhang2006,Zhaofsm1,me,mehweno1}, where the methods are implicit and need to solve nonlinear equations. However, it can be shown that the FE-FSM with LF Hamiltonian are equivalent to the fast sweeping method with LF Hamiltonian in \cite{mehweno1,me}, if the same high order reconstruction is used.

\begin{remark}
In Remark 2.1 of the recent work \cite{mehweno1} to study HWENO FSM, we observed that such method does not converge to machine epsilon when the Godunov flux is used for some examples. To fix this issue, it was proposed there to update the solution by
$$\phi^{new} = \omega \phi^{new} + (1-\omega) \phi^{old}, \qquad 0<\omega<1,$$
which may increase the number of iteration steps. The numerical examples suggest that there is no need to introduce any other parameters for HWENO fixed-point sweeping methods proposed in this paper.
\end{remark}

\begin{remark}
It was numerically observed in \cite{xiong,mehweno1} that, in some test cases (e.g. Example 6 and 7 in Section \ref{sec:numerical}), the convergence of the HWENO FSM scheme is very sensitive to the parameter $\varepsilon$ used in the nonlinear weights of HWENO reconstruction procedure. In order to address this issue, the parameter $\varepsilon$ was adjusted according to the mesh size, so that the scheme can converge quickly and provide the expected high order. However, the manual adjustment of $\varepsilon$ will impact the application of the proposed scheme as the best choice of $\varepsilon$ cannot be known a priori and may be problem and mesh dependent.
An interesting observation is that the proposed HWENO fixed-point sweeping method resolves this issue, and all of the numerical examples work well with a fixed parameter $\varepsilon$.
\end{remark}

\subsection{Hybrid strategy}
In the HWENO reconstruction procedure, compared with simple linear reconstruction, the evaluation of the smoothness indicators occupies most of the additional computational cost. Here we explore a hybrid strategy to combine the linear and HWENO method, which is similar to the hybrid fast sweeping WENO method studied in \cite{me,mehweno1}. In this paper, we directly follow the approach in \cite{mehweno1}, and apply the fifth order linear reconstruction if the numerical solution is monotonic on the big stencil $S_0$ or $\widetilde{S_{0}}$, and apply the HWENO-ZZQ reconstruction in other cases.

In order to describe the fixed-point fast sweeping method with hybrid strategy, we separate the the points $\{(x_i, y_j)\}$ in \emph{Category IV} (defined in Section \ref{sec:flowchart}) into the following two subcategories, which will be handled slightly different in the hybrid method. \\
\noindent\emph{Category IV.1}: For points whose distances to \emph{Category III} are less than or equal to $2h$ (excluding those in \emph{Category I}). \\
\noindent\emph{Category IV.2}: All the remaining points in Category IV.

The specific flowchart of fixed-point sweeping method with hybrid strategy is similar to that in Section \ref{sec:flowchart}, except the \textbf{Step 2}. During the \textbf{Step 2}, for the points in \emph{Category IV.1}:  the nonlinear HWENO-ZZQ reconstruction (\ref{weno}) is applied to evaluate $(\phi_{x})_{i,j}^{\pm}$. For the points in \emph{Category IV.2}, the hybrid strategy is applied, by using either the linear or nonlinear reconstruction based on the following criteria:
\begin{equation}\label{signf}
(\phi_{x})_{i,j}^{\pm}=
\begin{cases}
\eqref{hybridlinearf} ~\mathrm{or}~ \eqref{hybridlinearz}, &\mathrm{if~} \{u_{i,j}\}\mathrm{~have~the~same~sign~on~}S_{0}~or~ \widetilde{S_{0}},\\
\eqref{weno}, &      \mathrm{otherwise}.
\end{cases}
\end{equation}
The similar procedure for $(\phi_{y})_{i,j}^{\pm}$ is used. The rest of the algorithm is the same as that in Section \ref{sec:flowchart}.

\section{Numerical results}\label{sec:numerical}

\setcounter{equation}{0}\setcounter{figure}{0}\setcounter{table}{0}
In this section, we will present extensive numerical examples by testing the proposed fifth order finite difference HWENO fixed-point sweeping methods on the Eikonal equation and general static HJ equations in two dimensions. We will compare the numerical results of four fixed-point sweeping methods with/without the hybrid strategy, and list their errors, convergence rates and the numbers of iterations. In all the numerical examples, $\varepsilon$ in \eqref{epsilon} is taken as $10^{-6}$, and linear weights are taken as $\gamma_{1}=0.98$ and $\gamma_{2}=\gamma_{3}=0.01$. The total number of grid points is assumed to be $N_x=N_y=N$. We use ``iter'' to denote the number of iterations in all the tables. Note that one iteration means that all of the point values are updated once. Therefore, one GS fast sweeping iteration involves the sweeping in four alternating directions and would count as four iterations.  All the computations are implemented by using MATLAB 2020a on ThinkPad computer with 1.80 GHz Intel Core i7 processor and 16GB RAM.

\noindent \textbf{Example 1}. We solve the Eikonal equation with
$$f(x,y)=\frac{\pi}{2}\sqrt{\sin^{2}\Big(\pi+\frac{\pi}{2}x\Big)+\sin^{2}\Big(\pi+\frac{\pi}{2}y\Big)},$$
on the computational domain $[-1,1]^2$, with the inflow boundary $\Gamma=(0,0)$. The exact solution is given by
\begin{equation*}
\phi(x,y)=\cos\Big(\pi+\frac{\pi}{2}x\Big)+\cos\Big(\pi+\frac{\pi}{2}y\Big).
\end{equation*}

The Godunov numerical Hamiltonian \eqref{godunovflux} is used in this example. Table \ref{tab1} lists the numerical results of four fixed-point sweeping methods, including the numerical errors, convergence order, number of iterations and CPU time. The numerical results without hybrid strategy are reported on the left side of Table \ref{tab1}, and those with the hybrid strategy are reported on the right side.

We can observe that the FE-Jacobi method requires the smaller CFL number of value $0.1$, but still fails to converge on the refined mesh. The reason is the that FE time discretization coupled with high-order linear spatial discretization suffers from linear stability problems. However, when the third RK time discretization is used, the RK-Jacobi scheme can take a larger CFL number of value $1$. Also, the number of iterations of RK-Jacobi method is smaller than FE-Jacobi method on the same mesh. The fast sweeping technique (in the FE-FSM and RK-FSM) can improve the convergence of the Jacobi scheme, as observed in Table \ref{tab1}. In addition, when the fast sweeping technique is used, the FE-FSM can now use a larger CFL number than FE-Jacobi scheme. On the same refined mesh, we observe that the RK-FSM only takes about $50\%$ CPU time of the RK-Jacobi scheme. Furthermore, the FE-FSM costs even less CPU time than RK-FSM. The numerical results after the hybrid strategy is applied can be seen on the right side of Table \ref{tab1}, which suggests that the hybrid strategy can save $70\%-80\%$ of the CPU time on the refined mesh.
In summary, the FE-FSM performs the best out of these four methods, and this observation is consistent with that in \cite{fixhjzhang,fixhyperbolicwu}. The hybrid strategy can further reduce the computational cost.

\begin{figure}[H]
\centering
\subfigure[FE-FSM-Convergence.]{
\includegraphics[width=5.1cm]{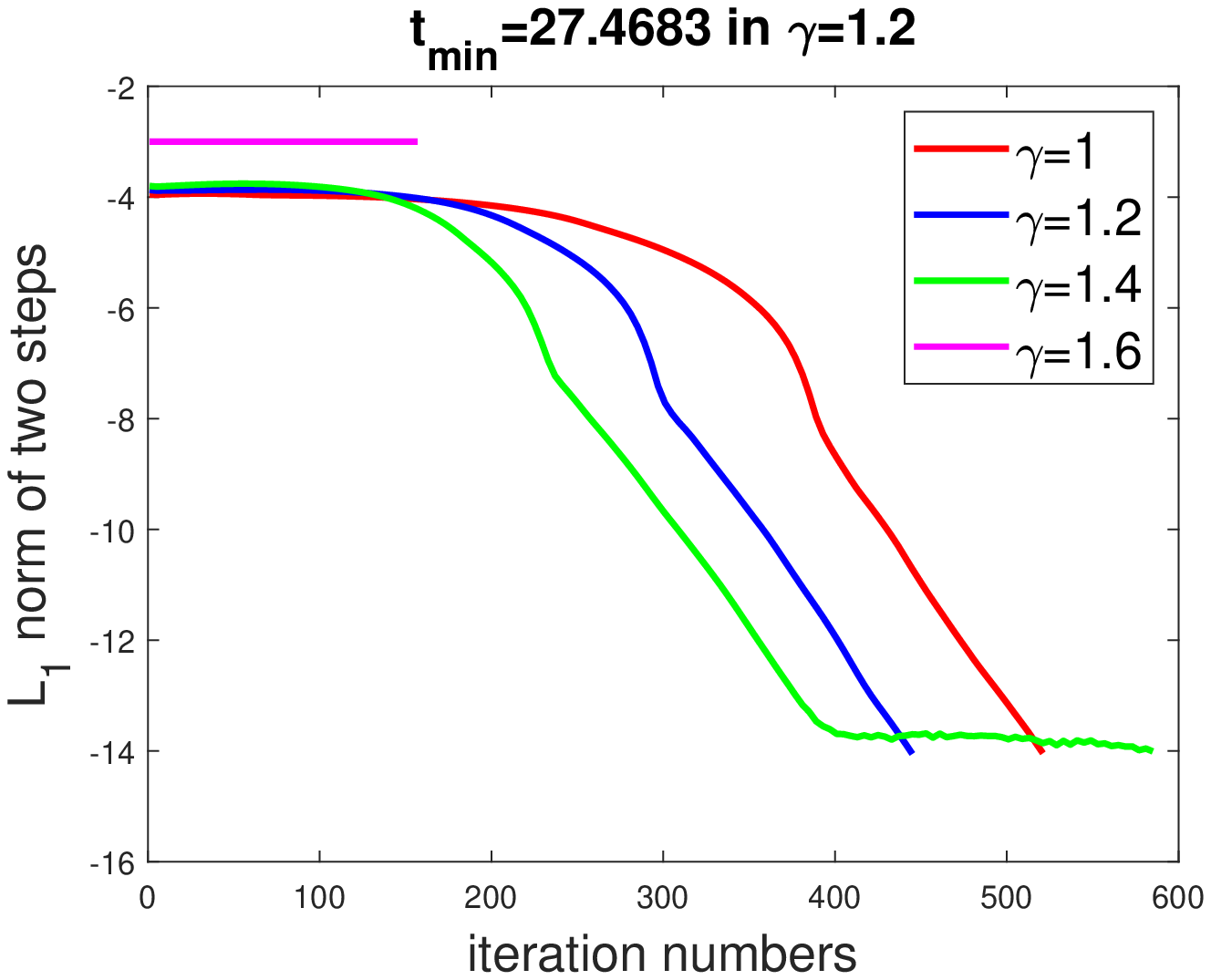}
}
\subfigure[FE-FSM-Residual.]{
\includegraphics[width=5.1cm]{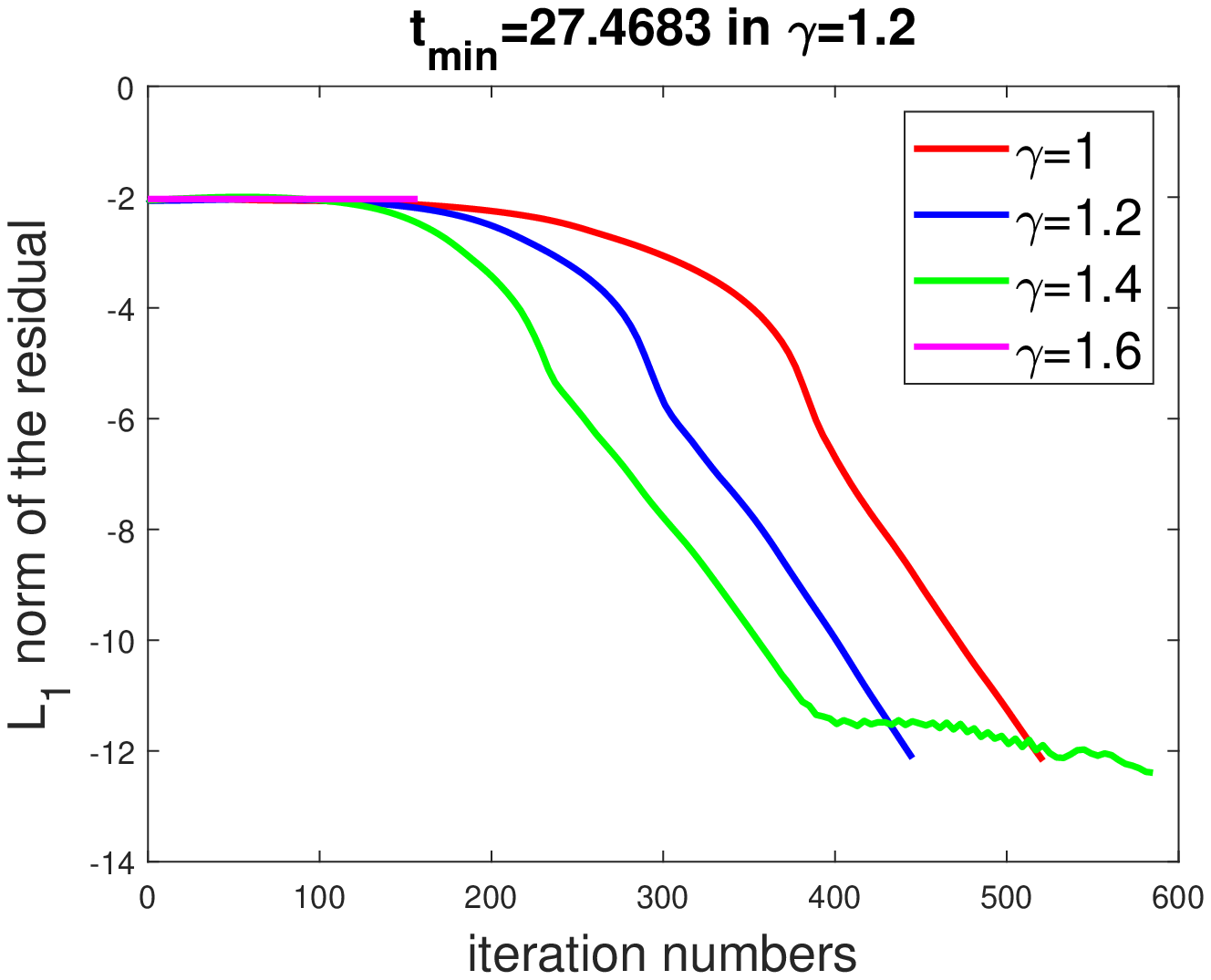}
}
\subfigure[FE-FSM-Numerical error.]{
\includegraphics[width=5.1cm]{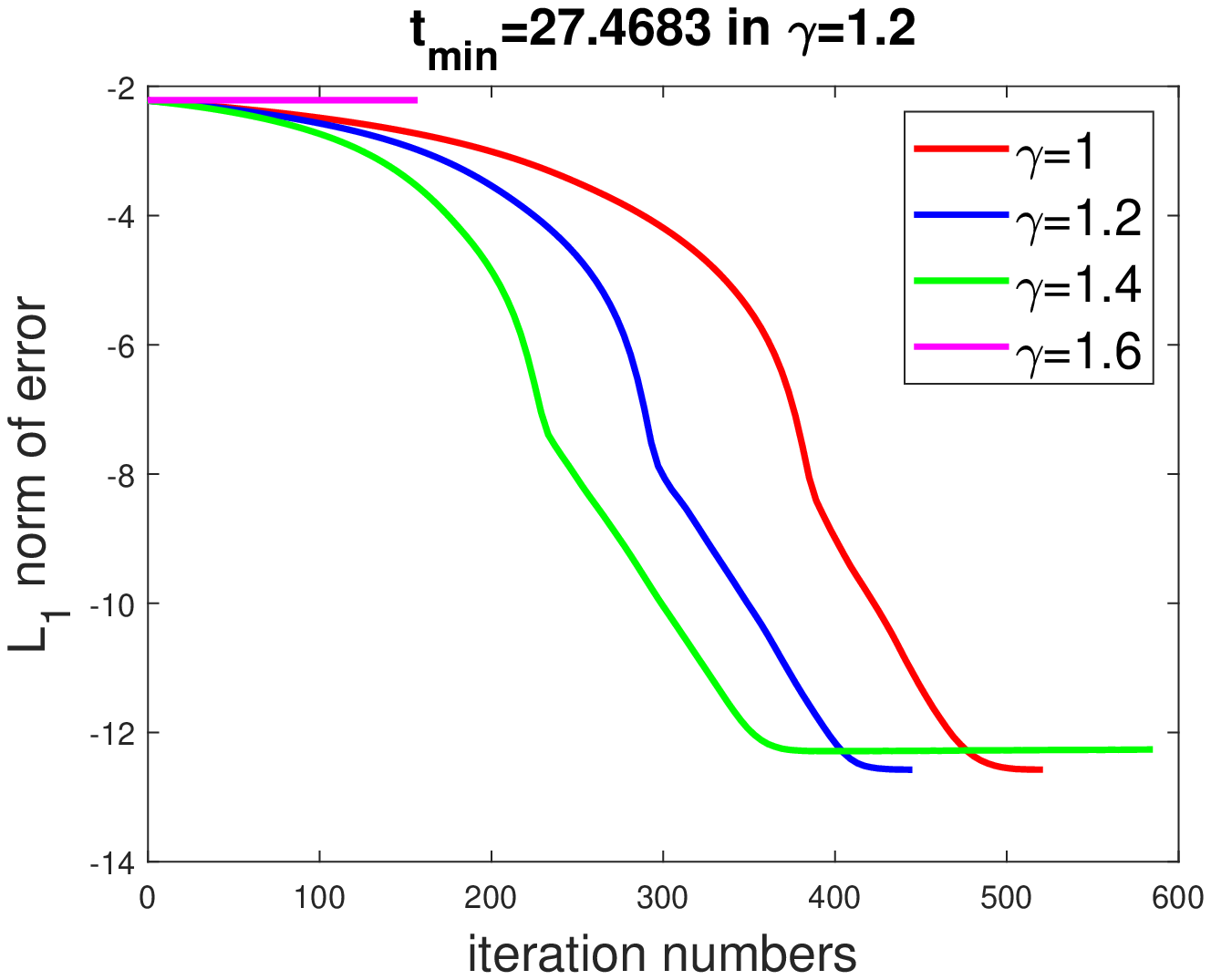}
}
\subfigure[RK-FSM-Convergence.]{
\includegraphics[width=5.1cm]{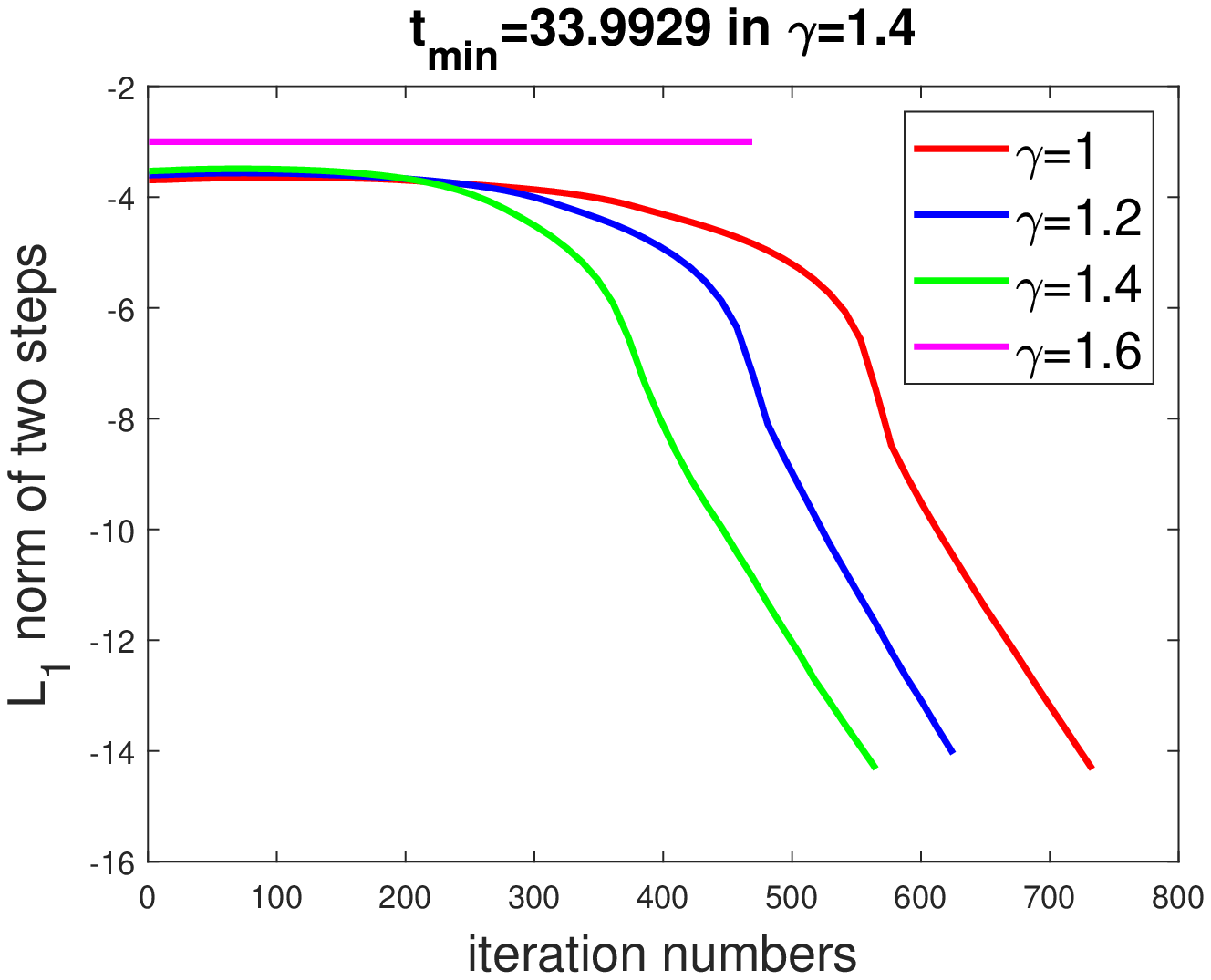}
}
\subfigure[RK-FSM-Residual.]{
\includegraphics[width=5.1cm]{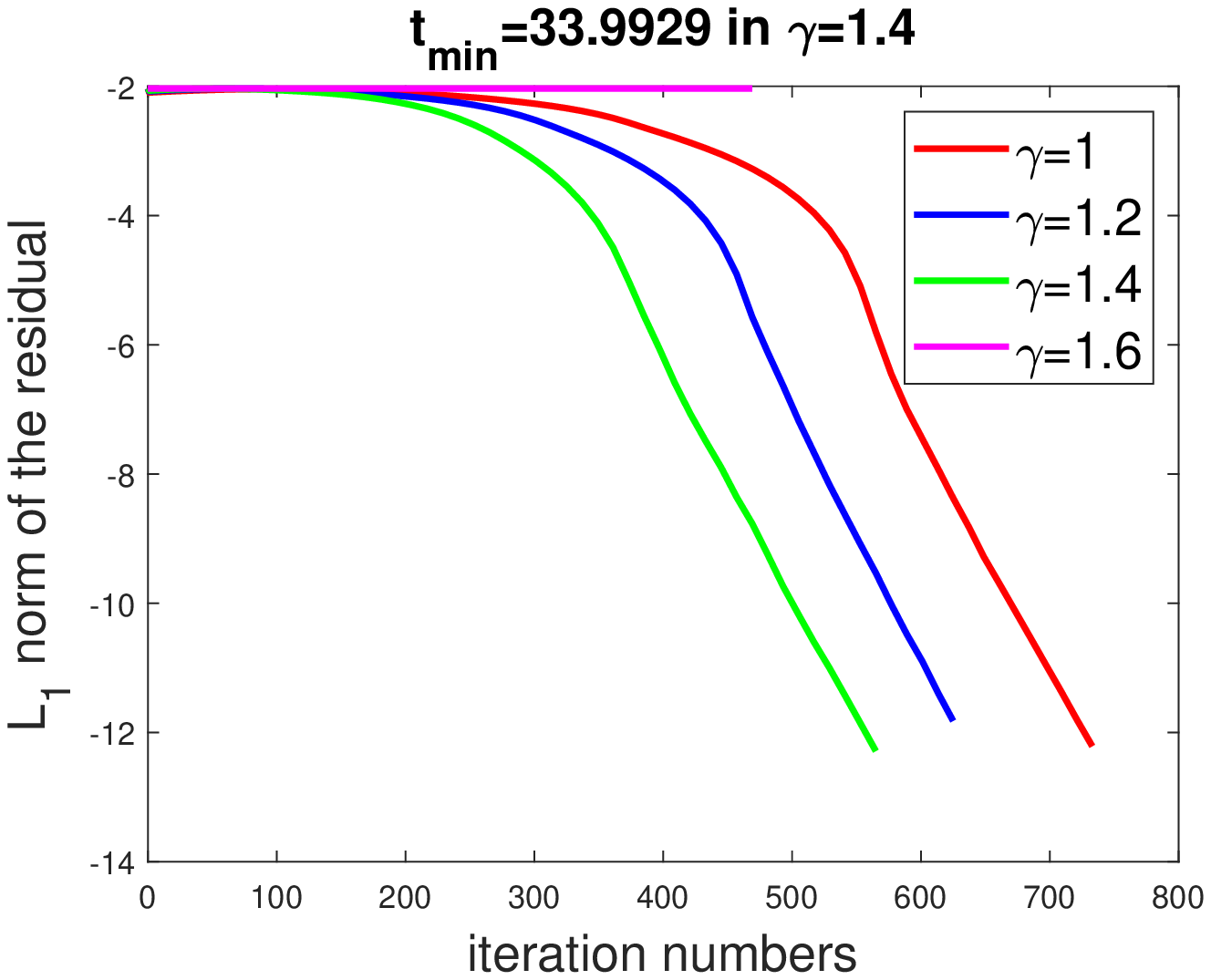}
}
\subfigure[RK-FSM-Numerical error.]{
\includegraphics[width=5.1cm]{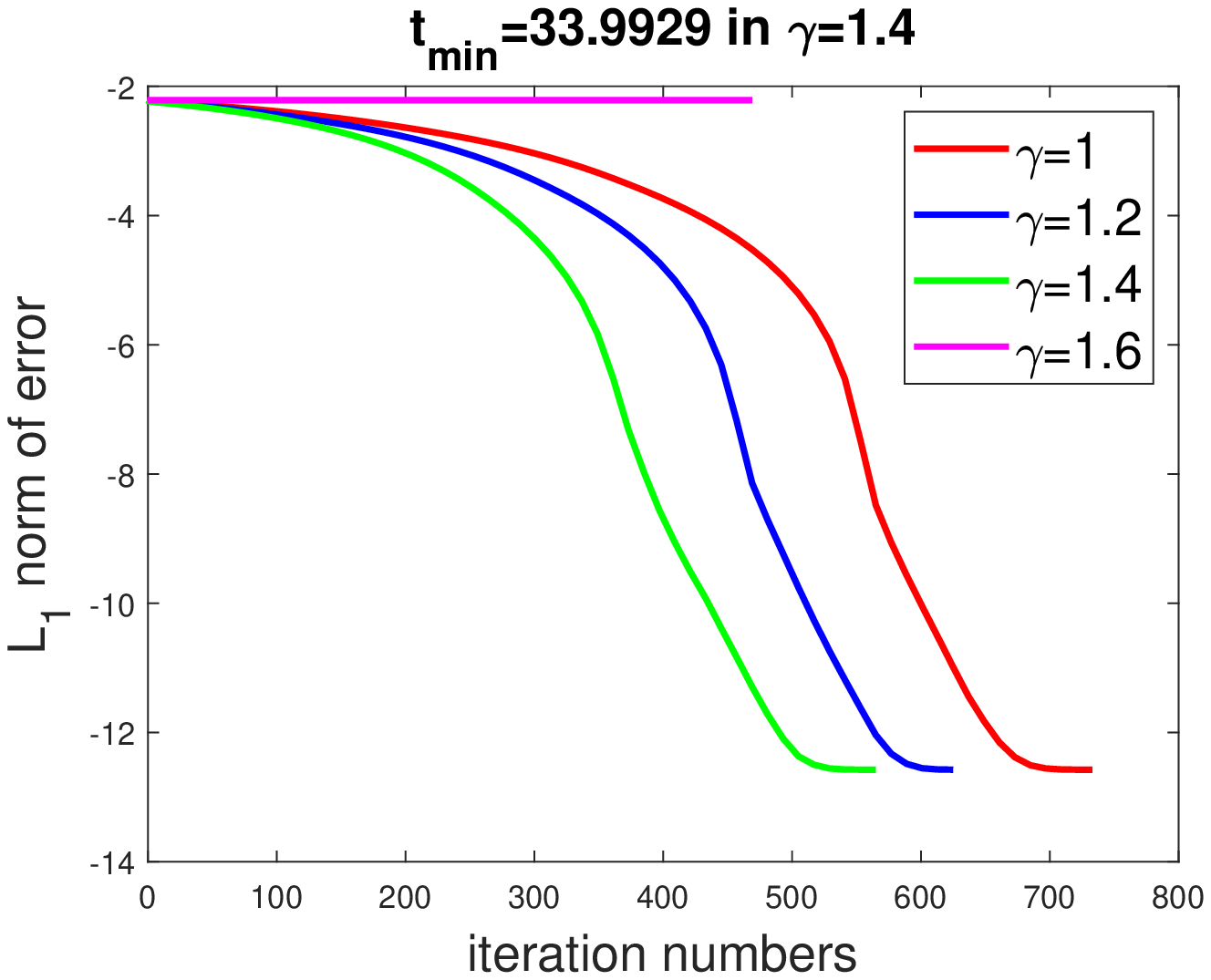}
}
\subfigure[FE-FSM-h-Convergence.]{
\includegraphics[width=5.1cm]{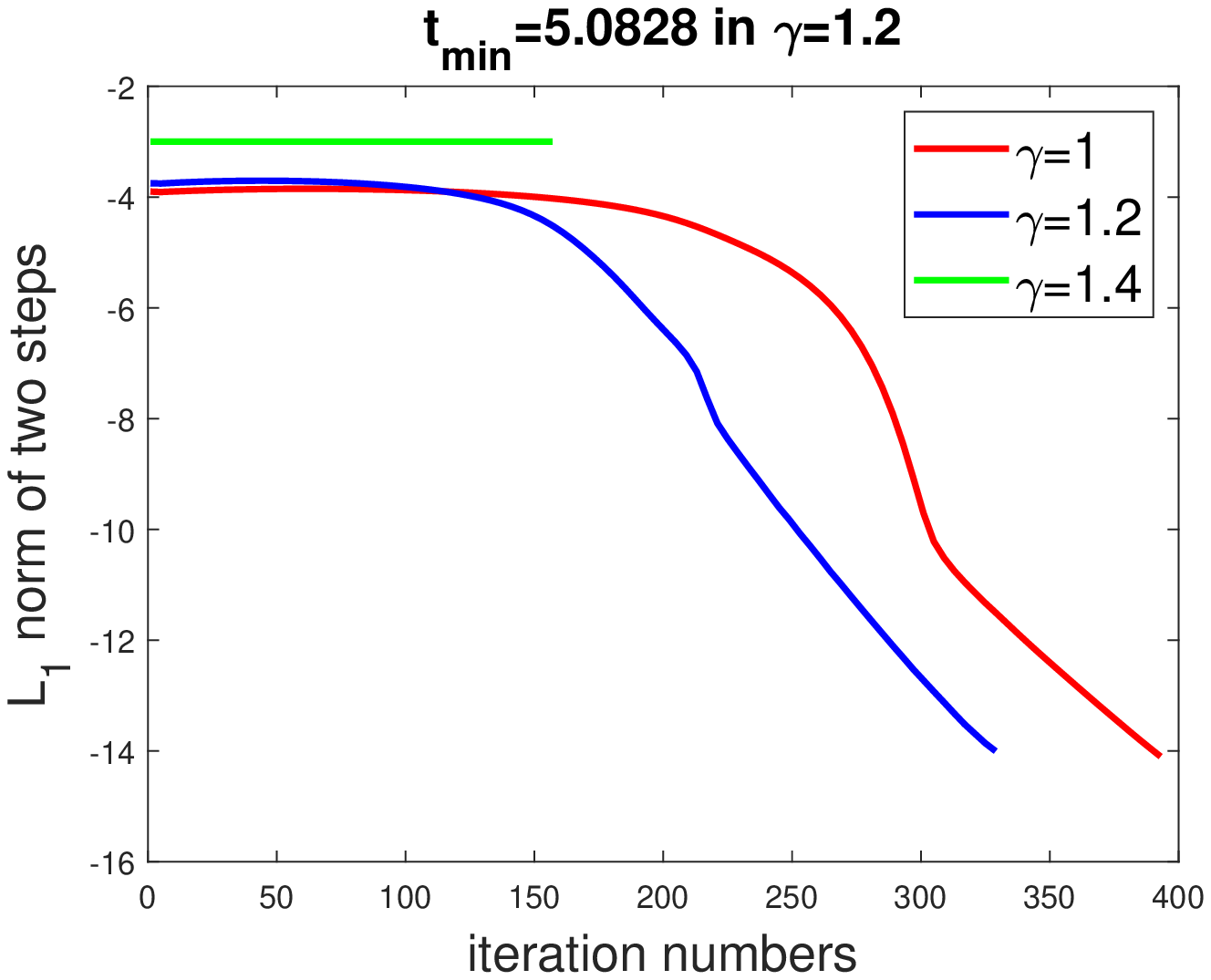}
}
\subfigure[FE-FSM-h-Residual.]{
\includegraphics[width=5.1cm]{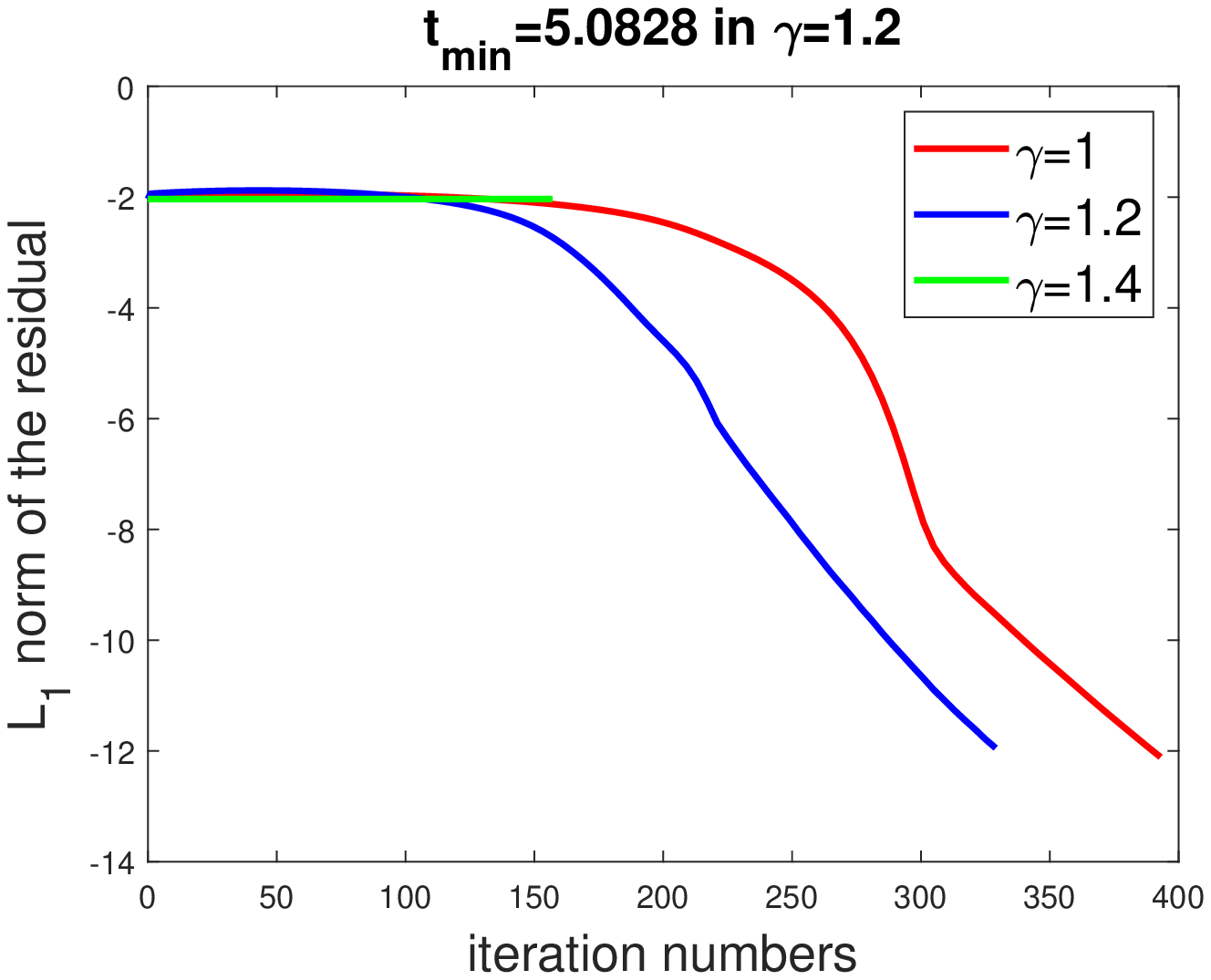}
}
\subfigure[FE-FSM-h-Numerical error.]{
\includegraphics[width=5.1cm]{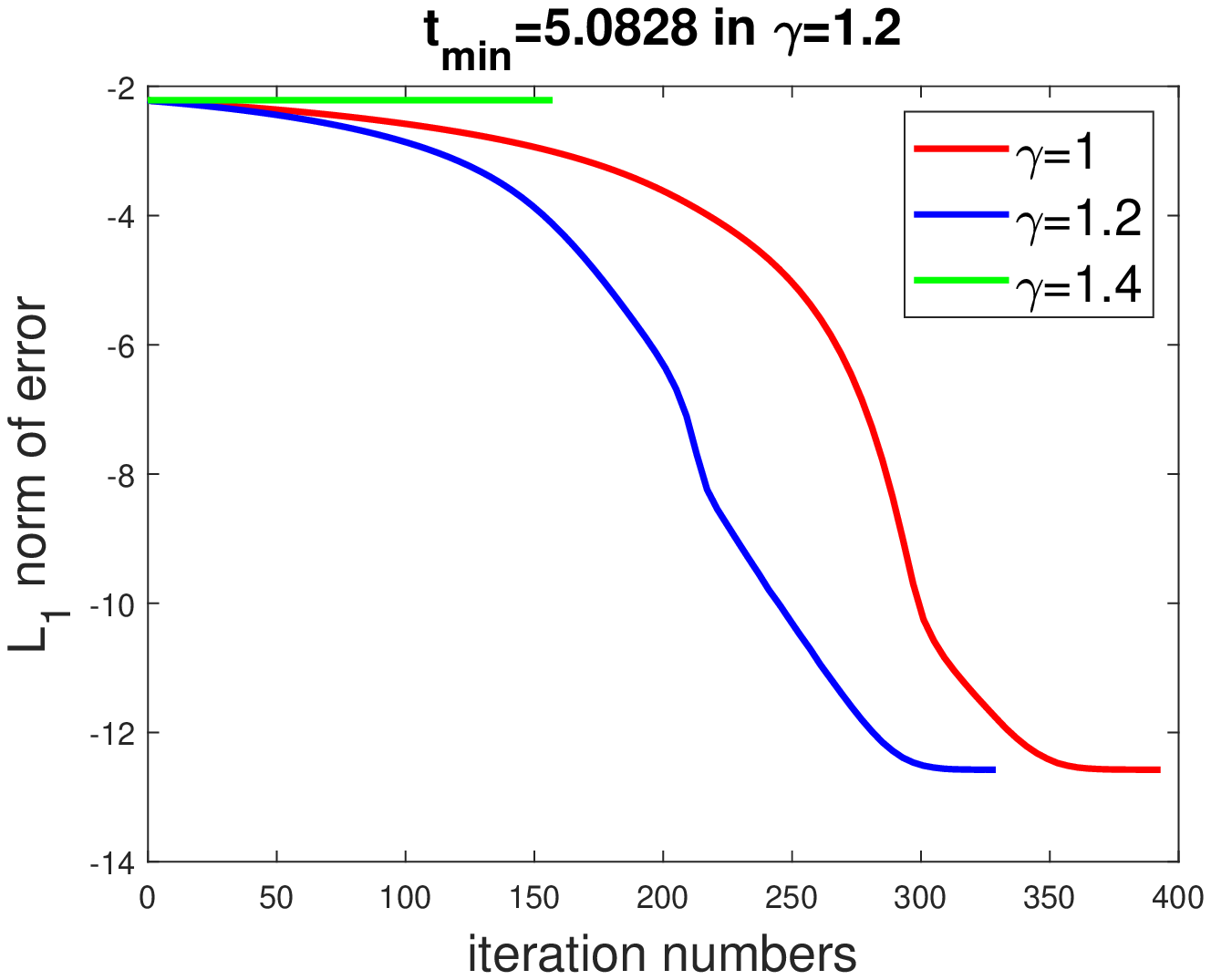}
}
\subfigure[RK-FSM-h-Convergence.]{
\includegraphics[width=5.1cm]{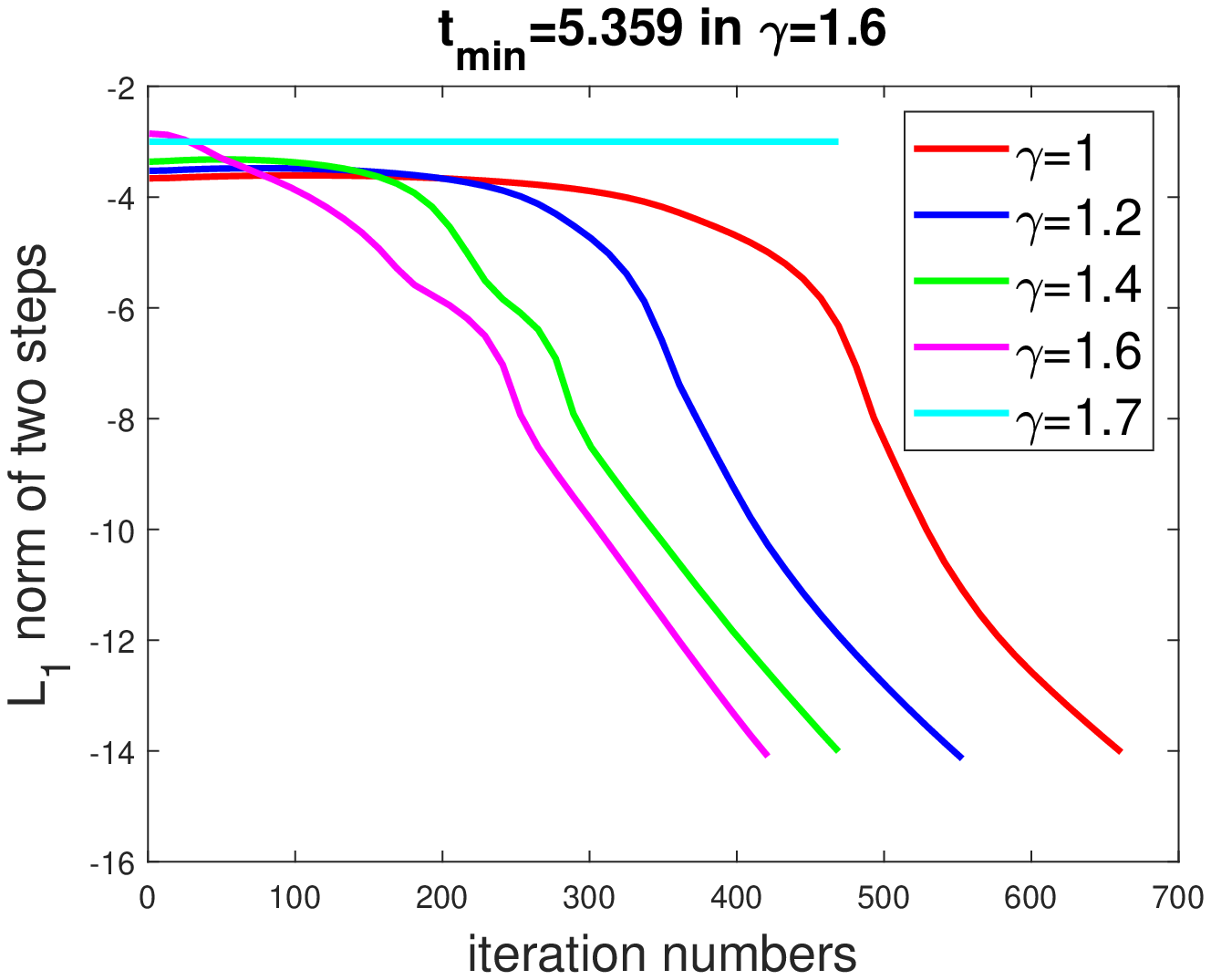}
}
\subfigure[RK-FSM-h-Residual.]{
\includegraphics[width=5.1cm]{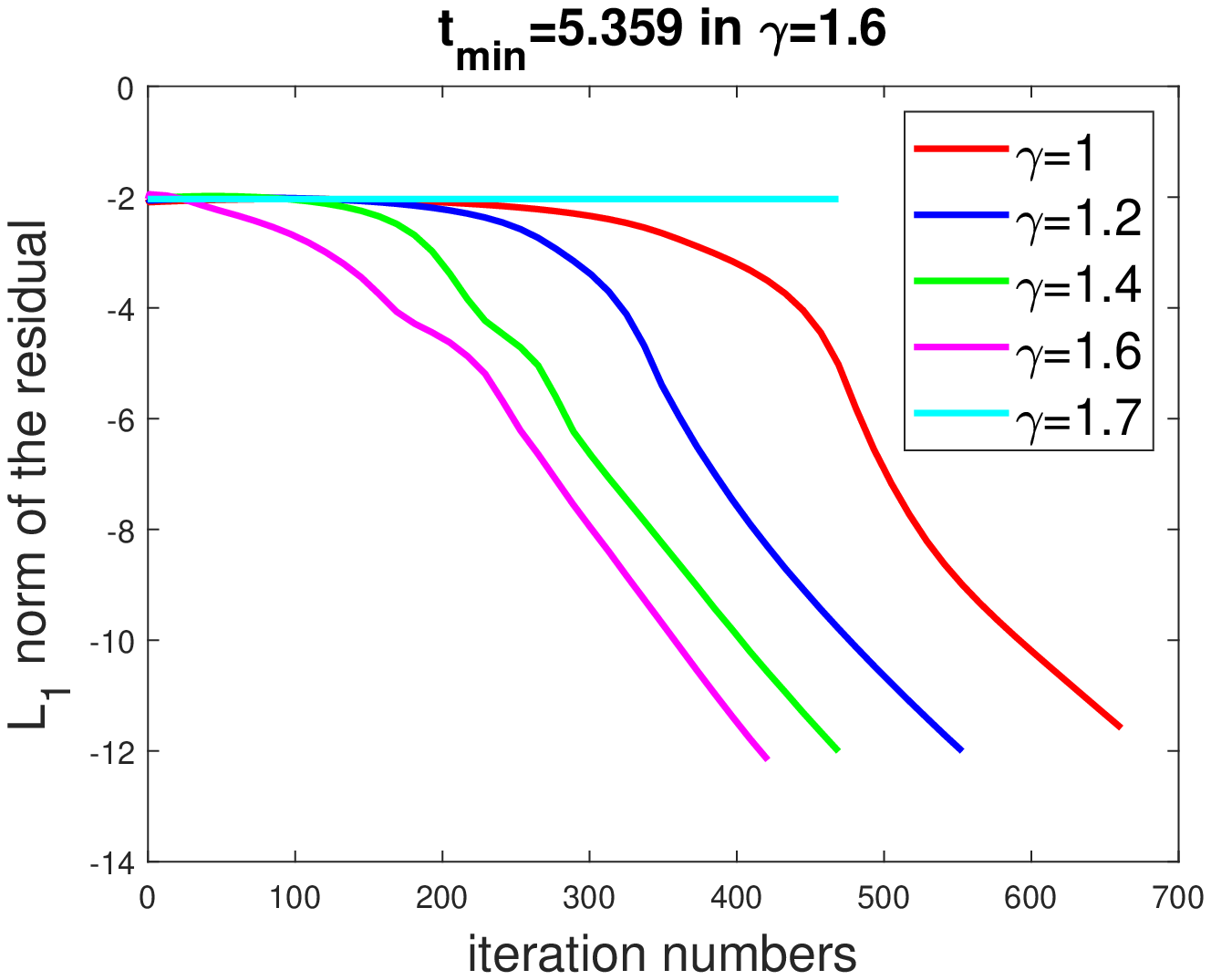}
}
\subfigure[RK-FSM-h-Numerical error.]{
\includegraphics[width=5.1cm]{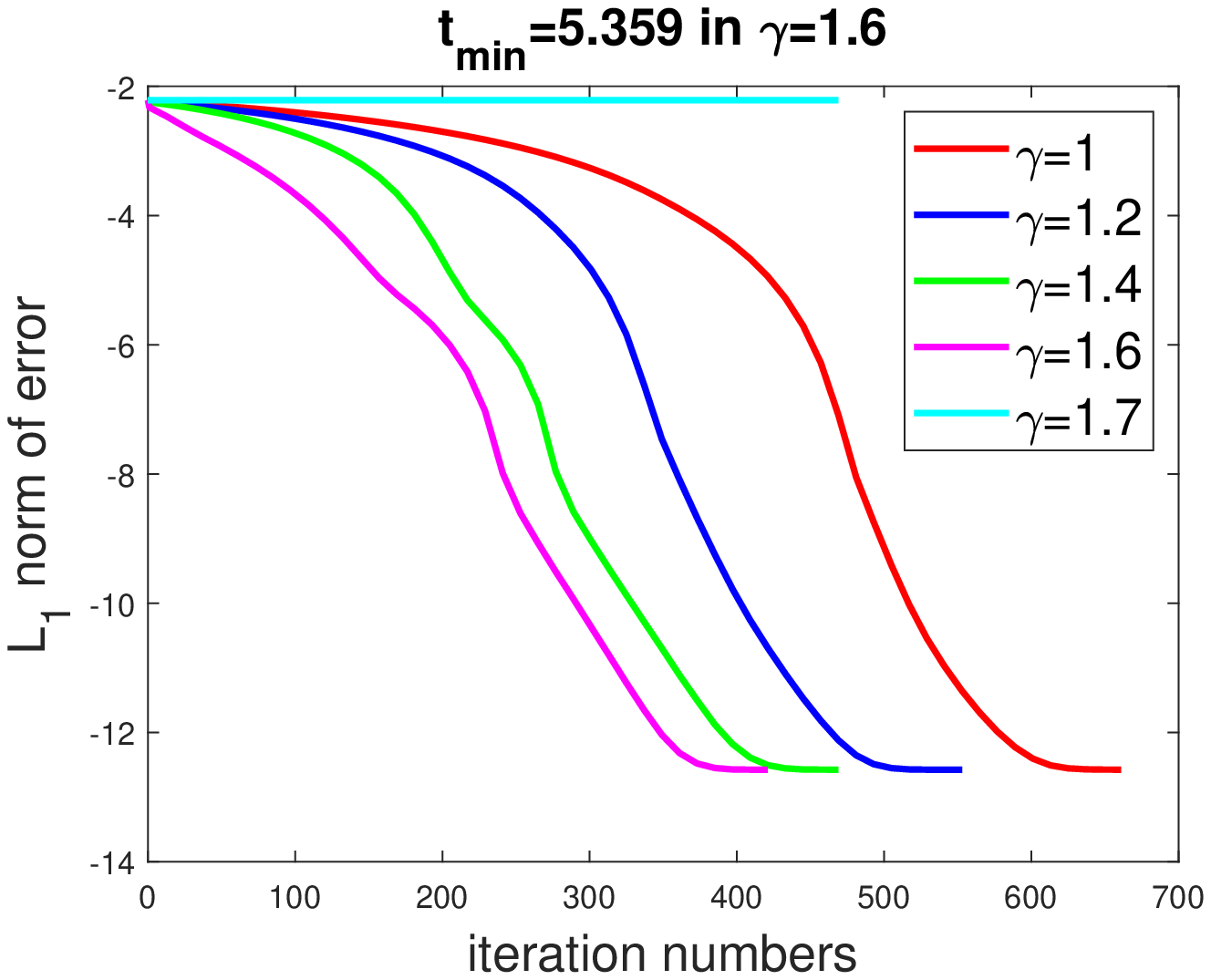}
}
\caption{Example 1. Study on different $\gamma's$}\label{fige1}
\end{figure}

Finally, we want to remark that, for the three schemes other than the FE-Jacobi scheme, the CFL number can be taken to be greater than $1$.
Figure \ref{fige1} shows that the convergence history of the FE-FSM and RK-FSM with larger CFL numbers on mesh $N=320$, reporting
the errors between two consecutive iteration steps, the residual and numerical errors. Here, ``FE-FSM-h'' represents FE-FSM with hybrid strategy, and similarly for ``RK-FSM-h''. We have listed the optimal CFL number and CPU time of each scheme in the title of each sub-figure. For example, ``$t_{min}=24.7683$ (unit: second) in $\gamma=1.2$'' on first line of Figure \ref{fige1} represents the optimal CFL number for FE-FSM is $1.2$, because it requires the least CPU time $t_{min}=24.7683$ among all the choices of $\gamma$ that lead to convergent results.
We also observe that the RK-FSM tends to admit a larger CFL number than the FE-FSM, however the RK-FSM still costs more CPU time than the FE-FSM.

\begin{table}
\caption{Example 1. Comparison of the four methods: The errors of the numerical solution, the accuracy obtained and the number of iterations for convergence}\label{tab1}
	\begin{center}
\resizebox{\textwidth}{55mm}{
		\begin{tabular}{|c|cccccc|cccccc|}
			\hline
\multicolumn{7}{|c|}{FE-Jacobi $\gamma=0.1$}&\multicolumn{6}{|c|}{FE-Jacobi $\gamma=0.1$ with hybrid strategy}\\ \hline
N & $L_{1}$&order&  $L_{\infty}$&order&iter &time& $L_{1}$&order&  $L_{\infty}$&order&iter&time\\ \hline
40 &3.65e-06 &- &2.97e-05 &- &1690 &1.8799&3.62e-06 &- &2.93e-05 &- &1690 &0.73858\\
80 &5.75e-08 &5.99 &1.07e-06 &4.78 &2113 &9.4906&5.74e-08 &5.97 &1.07e-06 &4.76 &2113 &2.5132\\
160 &2.91e-10 &7.62 &1.16e-08 &6.52 &3281 &58.4108&- &- &- &- &- &-\\
320 &- &- &- &- &- &-&- &- &- &- &- &-\\ \hline
\multicolumn{7}{|c|}{FE-FSM $\gamma=1$}&\multicolumn{6}{|c|}{FE-FSM $\gamma=1$ with hybrid strategy}\\ \hline
$N$ & $L_{1}$~error&order& $L_{\infty}$~error&order&iter&time& $L_{1}$&order&  $L_{\infty}$&order&iter&time\\ \hline
40 &3.65e-06 &- &2.97e-05 &- &240 &0.29083&3.62e-06 &- &2.93e-05 &- &240 &0.12409\\
80 &5.75e-08 &5.99 &1.07e-06 &4.78 &272 &1.0191&5.74e-08 &5.97 &1.07e-06 &4.76 &272 &0.35377\\
160 &2.91e-10 &7.62 &1.16e-08 &6.53 &348 &5.2703&2.91e-10 &7.62 &1.16e-08 &6.53 &332 &1.1835\\
320 &2.66e-13 &10.09 &9.62e-12 &10.23 &524 &32.0309&2.65e-13 &10.09 &1.03e-11 &10.12 &396 &6.8875\\ \hline
\multicolumn{7}{|c|}{RK-Jacobi $\gamma=1$}&\multicolumn{6}{|c|}{RK-Jacobi $\gamma=1$ with hybrid strategy}\\ \hline
$N$ & $L_{1}$~error&order& $L_{\infty}$~error&order&iter&time& $L_{1}$&order&  $L_{\infty}$&order&iter&time\\ \hline
40 &3.65e-06 &- &2.97e-05 &- &399 &0.68148&3.62e-06 &- &2.93e-05 &- &399 &0.19815\\
80 &5.75e-08 &5.99 &1.07e-06 &4.78 &471 &1.9727&5.74e-08 &5.97 &1.07e-06 &4.76 &471 &0.57089\\
160 &2.91e-10 &7.62 &1.15e-08 &6.54 &780 &12.9281&2.91e-10 &7.62 &1.15e-08 &6.54 &780 &3.5232\\
320 &3.04e-13 &9.90 &8.34e-11 &7.11 &1413 &102.4285&3.04e-13 &9.90 &8.34e-11 &7.11 &1413 &27.8141\\ \hline
\multicolumn{7}{|c|}{RK-FSM $\gamma=1$}&\multicolumn{6}{|c|}{RK-FSM $\gamma=1$ with hybrid strategy}\\ \hline
$N$ & $L_{1}$~error&order& $L_{\infty}$~error&order&iter&time& $L_{1}$&order&  $L_{\infty}$&order&iter&time\\ \hline
40 &3.65e-06 &- &2.97e-05 &- &360 &0.61658&3.62e-06 &- &2.93e-05 &- &360 &0.35017\\
80 &5.75e-08 &5.99 &1.07e-06 &4.78 &420 &1.7162&5.74e-08 &5.97 &1.07e-06 &4.76 &420 &0.73419\\
160 &2.91e-10 &7.62 &1.16e-08 &6.53 &516 &8.4034&2.91e-10 &7.62 &1.16e-08 &6.53 &516 &2.1677\\
320 &2.65e-13 &10.09 &9.33e-12 &10.28 &744 &46.2991&2.65e-13 &10.09 &1.10e-11 &10.03 &672 &11.431\\ \hline
\end{tabular}}
\end{center}
\end{table}
\newpage

\begin{figure}[!h]
\begin{center}
	\includegraphics[width=6.6cm]{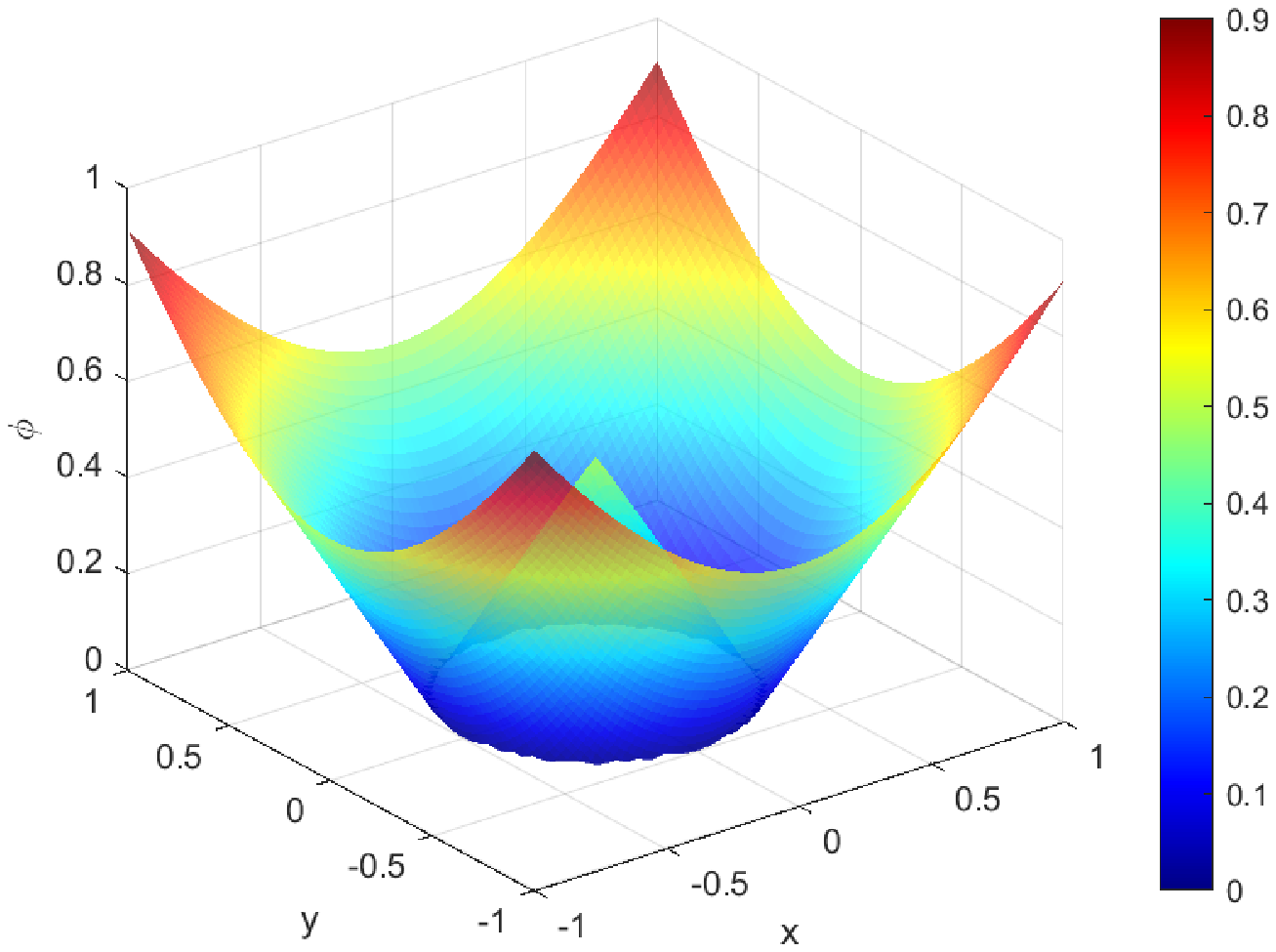}
	\includegraphics[width=6.6cm]{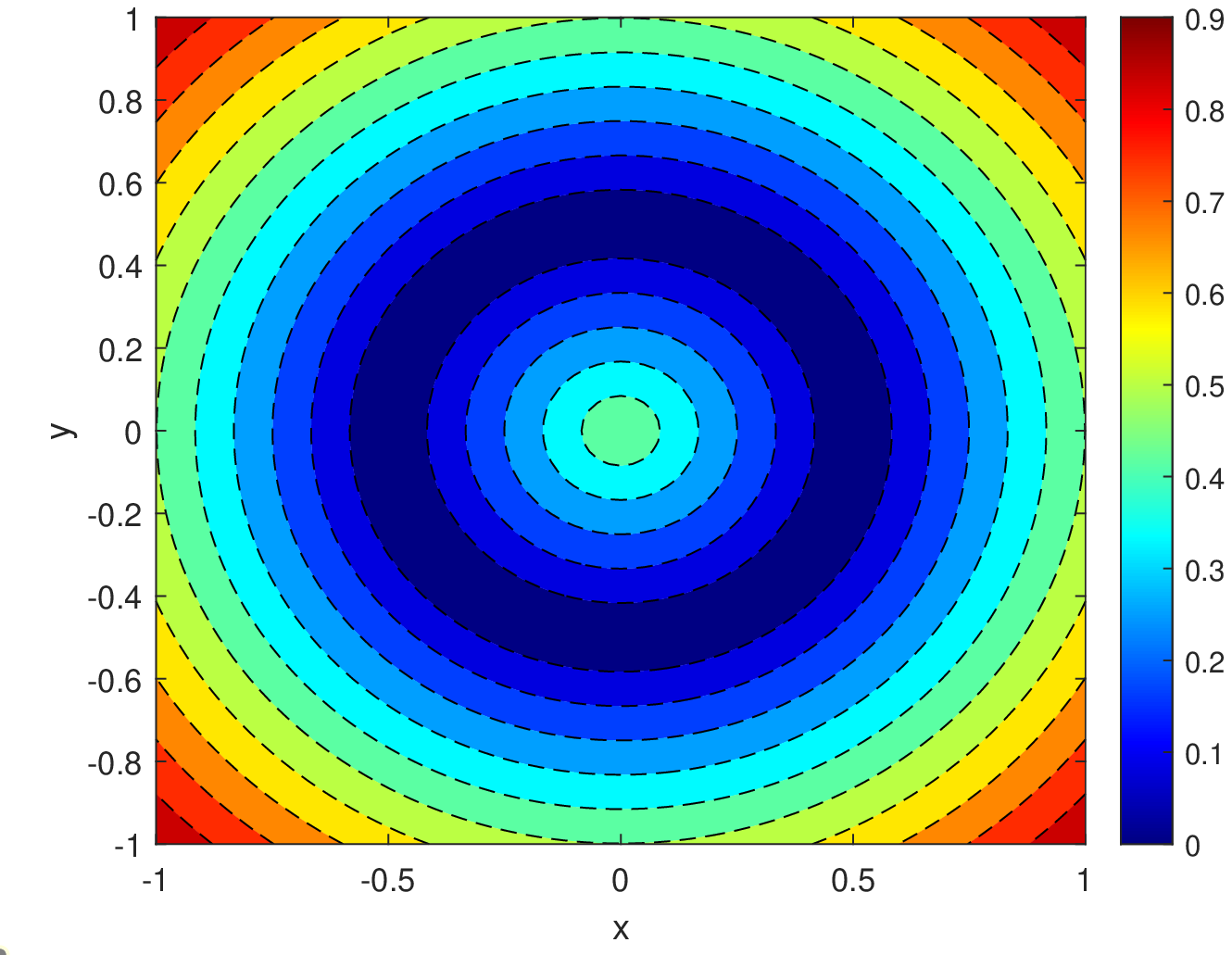}
\end{center}
\caption{Example 2. The numerical solution by FE-FSM on mesh $N=80$. Left: the 3D plot of numerical solution $\phi$; Right: the contour plot for $\phi$ .}\label{fige2sc}
\end{figure}

\noindent\textbf{Example 2}. We solve the Eikonal equation with $f(x,y)=1$. The computational domain is set as $[-1,1]^2$, and the inflow boundary $\Gamma$ is the circle with center at $(0,0)$ and radius $0.5$, that is
$$\Gamma=\left\{(x,y)|x^2+y^2=\frac{1}{4}\right\}.$$
The boundary condition is given as $\phi(x,y)=0$ on $\Gamma$. The exact solution is a distance function to the circle $\Gamma$, and it has a singularity at the center of the circle (due to the intersection of characteristic lines).
The Godunov numerical Hamiltonian is used, and the numerical errors are measured in the box $[-0.9,0.9]^2$ and outside the box $[-0.15,0.15]^2$, which aims to remove the influence of singularity and outflow boundary treatment.
The surface and contour of numerical solution by FE-FSM are shown in Figure \ref{fige2sc}.

The numerical results without the hybrid strategy are reported on the left side of Table \ref{tab2}, and those with the hybrid strategy are reported on the right side. Again, the FE-Jacobi iteration requires a smaller CFL number of $0.1$ due to the linear instability, and fails to converge on the refined mesh. For the other three methods, a larger CFL number can be taken. Again, we can observe that the fast sweeping technique can improve the convergence of the Jacobi scheme. On the same refined mesh, it can be seen that the RK-FSM only takes about $50\%$ CPU time of the RK-Jacobi scheme. Furthermore, the FE-FSM costs even less CPU time than RK-FSM. The numerical results obtained with the hybrid strategy can be seen on the right side of Table \ref{tab2}, which indicates that the hybrid strategy can save about $50\%-75\%$ CPU time on the refined mesh.

\begin{figure}[H]
\centering
\subfigure[FE-FSM-Convergence.]{
\includegraphics[width=5.1cm]{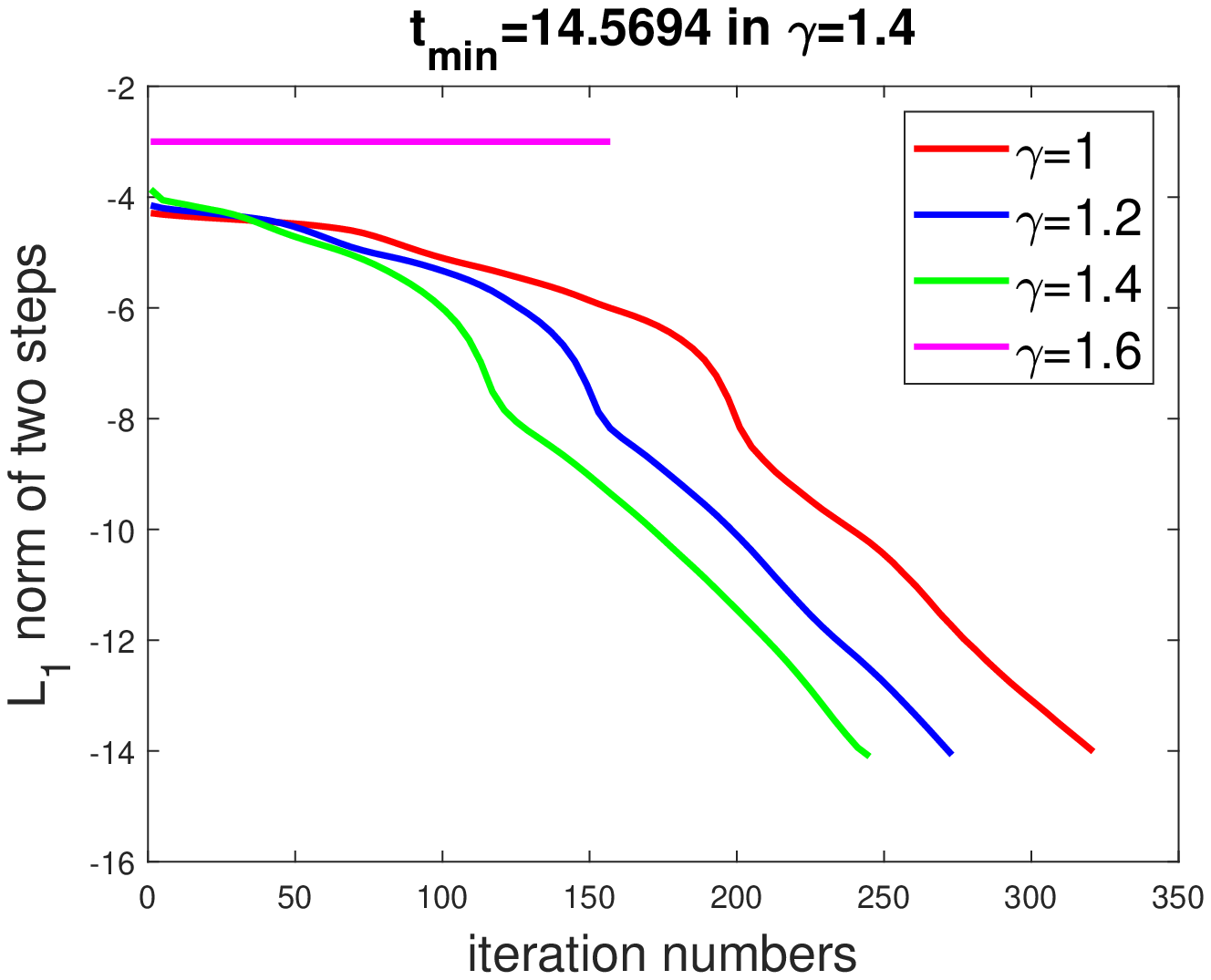}
}
\subfigure[FE-FSM-Residual.]{
\includegraphics[width=5.1cm]{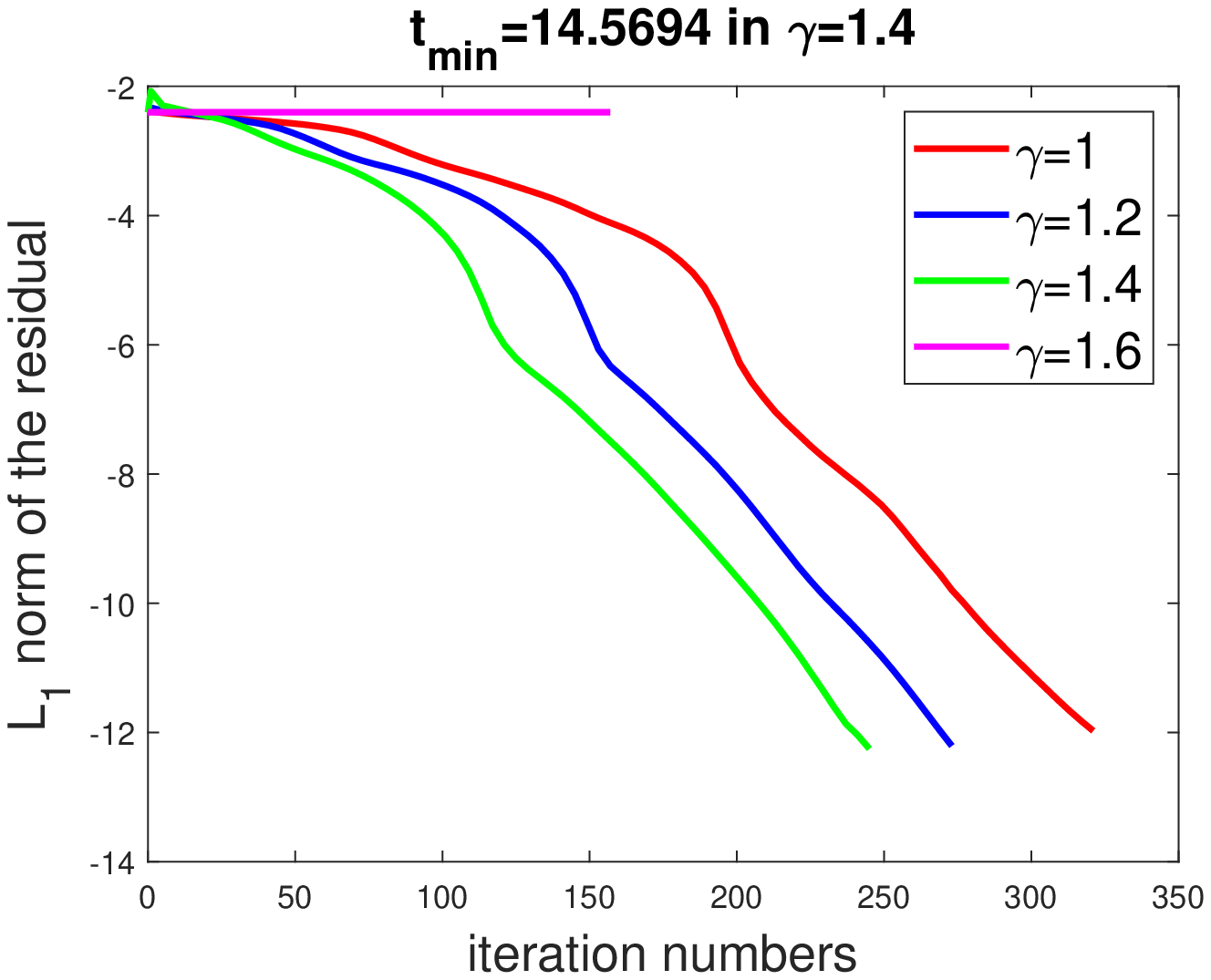}
}
\subfigure[FE-FSM-Numerical error.]{
\includegraphics[width=5.1cm]{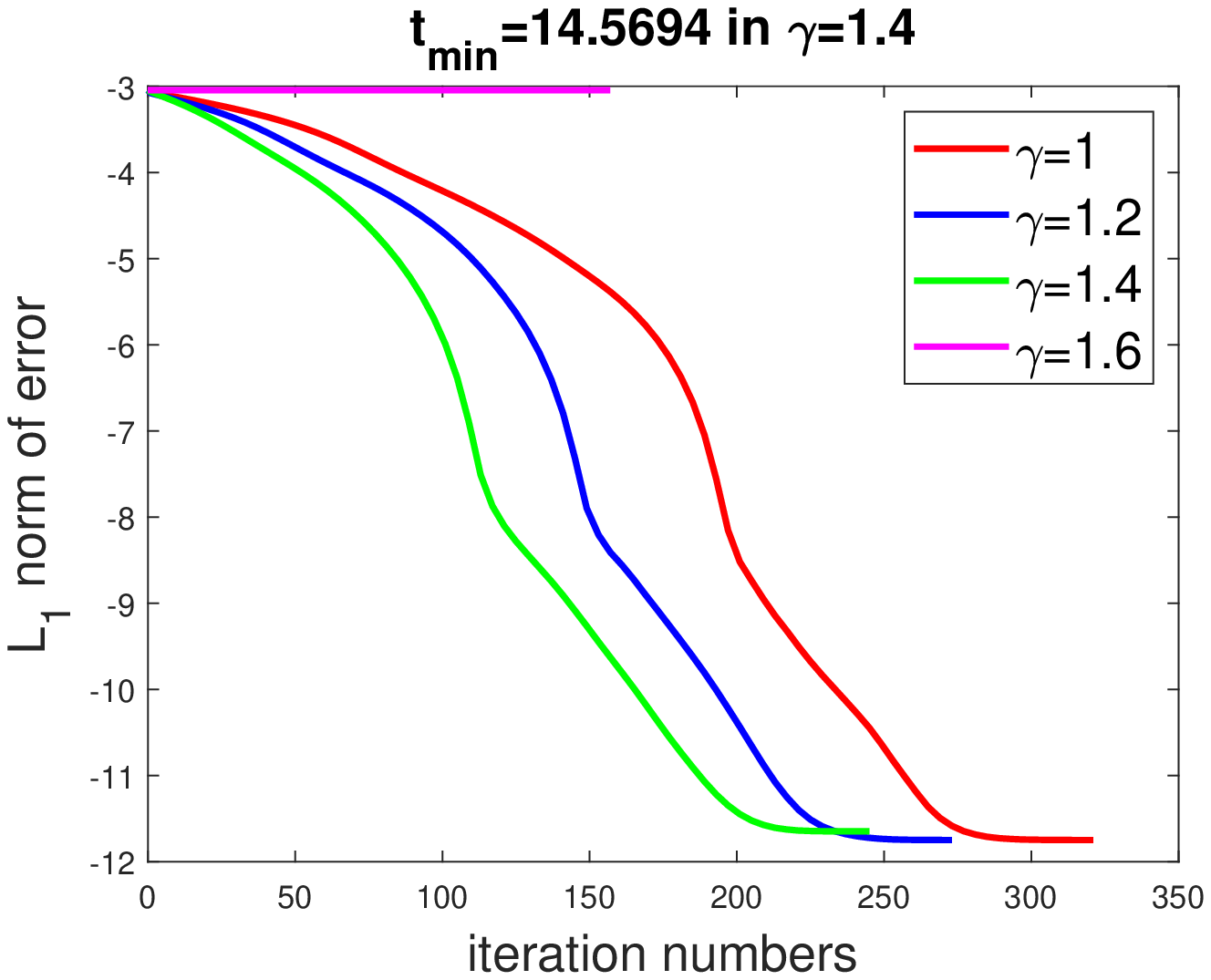}
}
\subfigure[RK-FSM-Convergence.]{
\includegraphics[width=5.1cm]{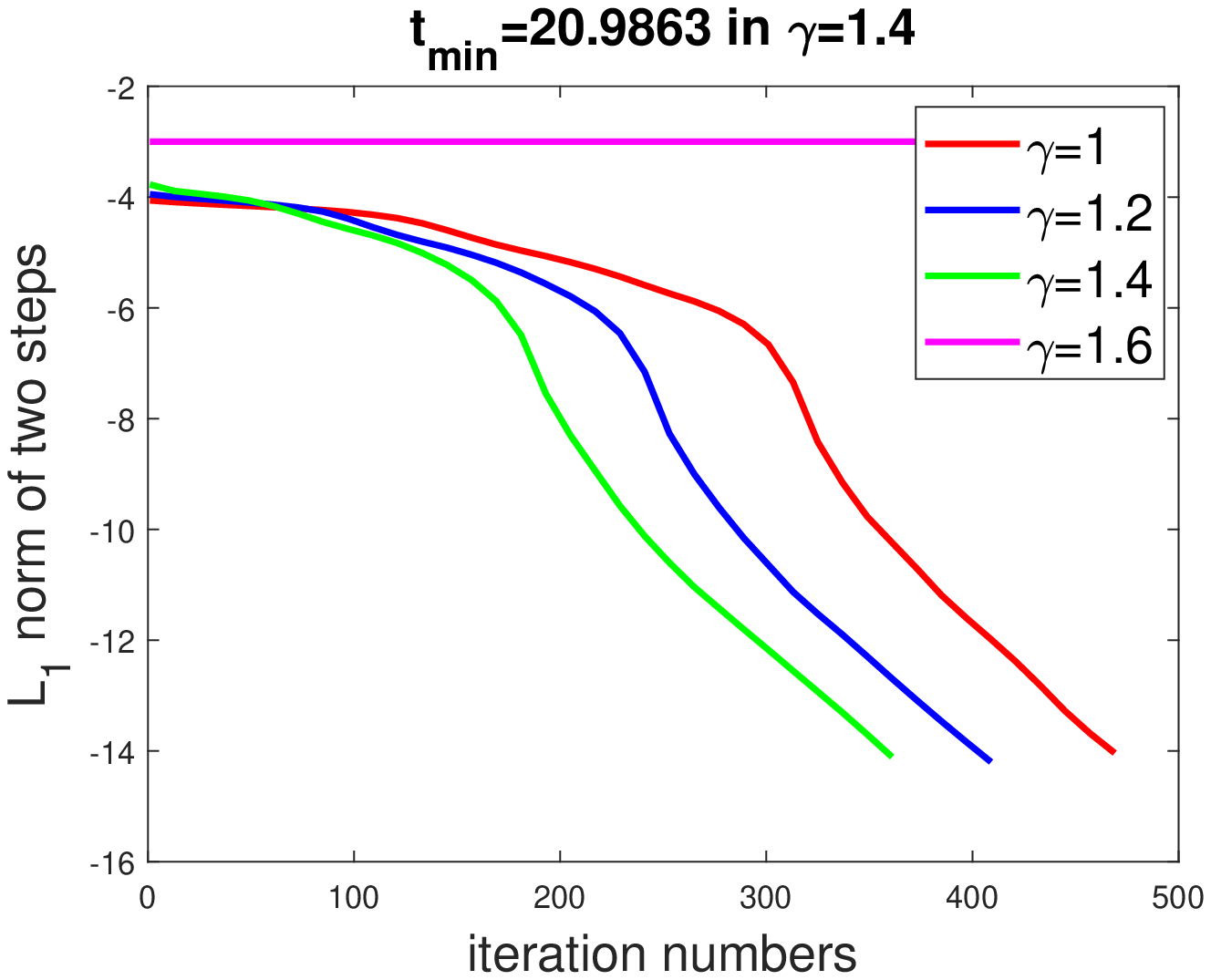}
}
\subfigure[RK-FSM-Residual.]{
\includegraphics[width=5.1cm]{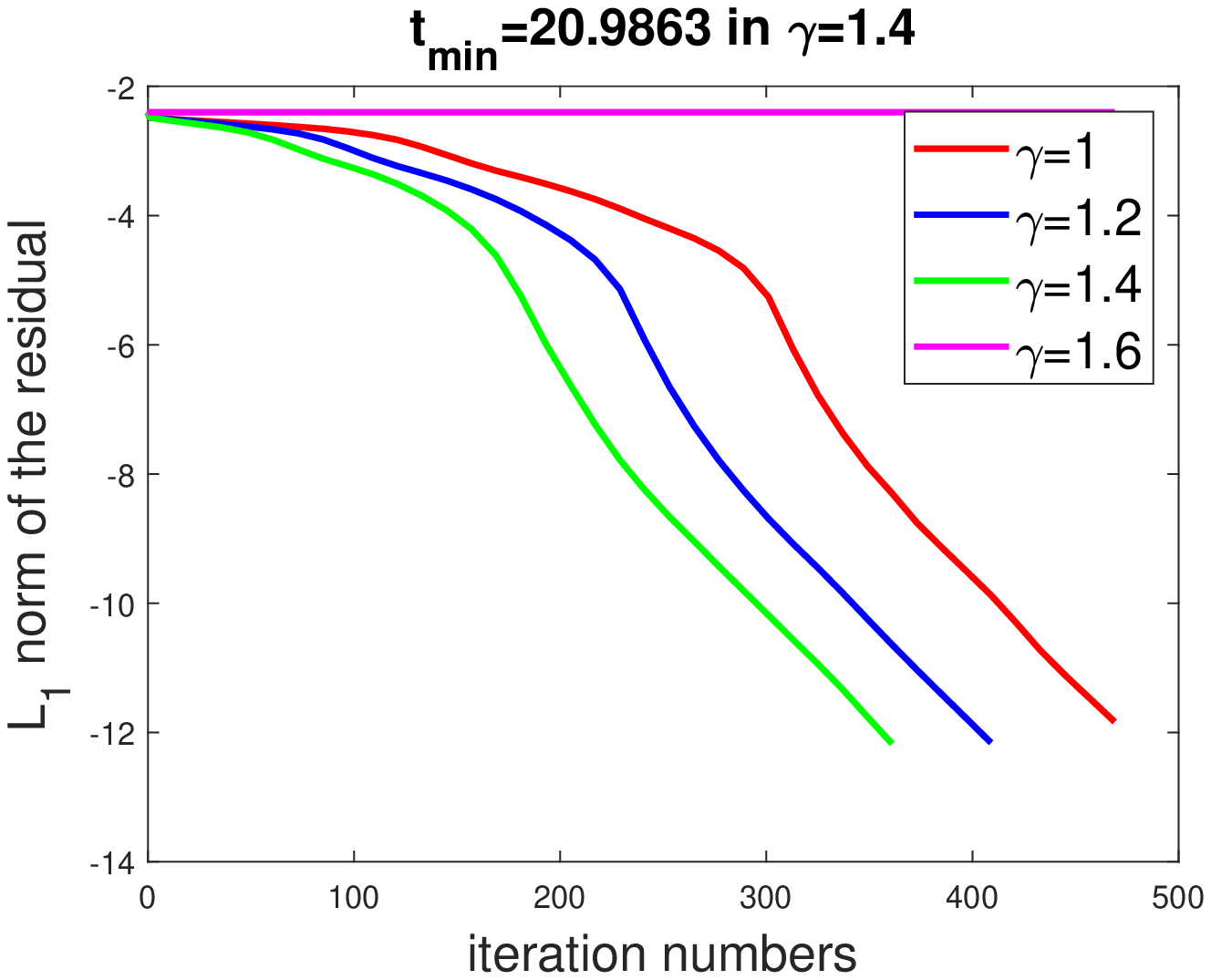}
}
\subfigure[RK-FSM-Numerical error.]{
\includegraphics[width=5.1cm]{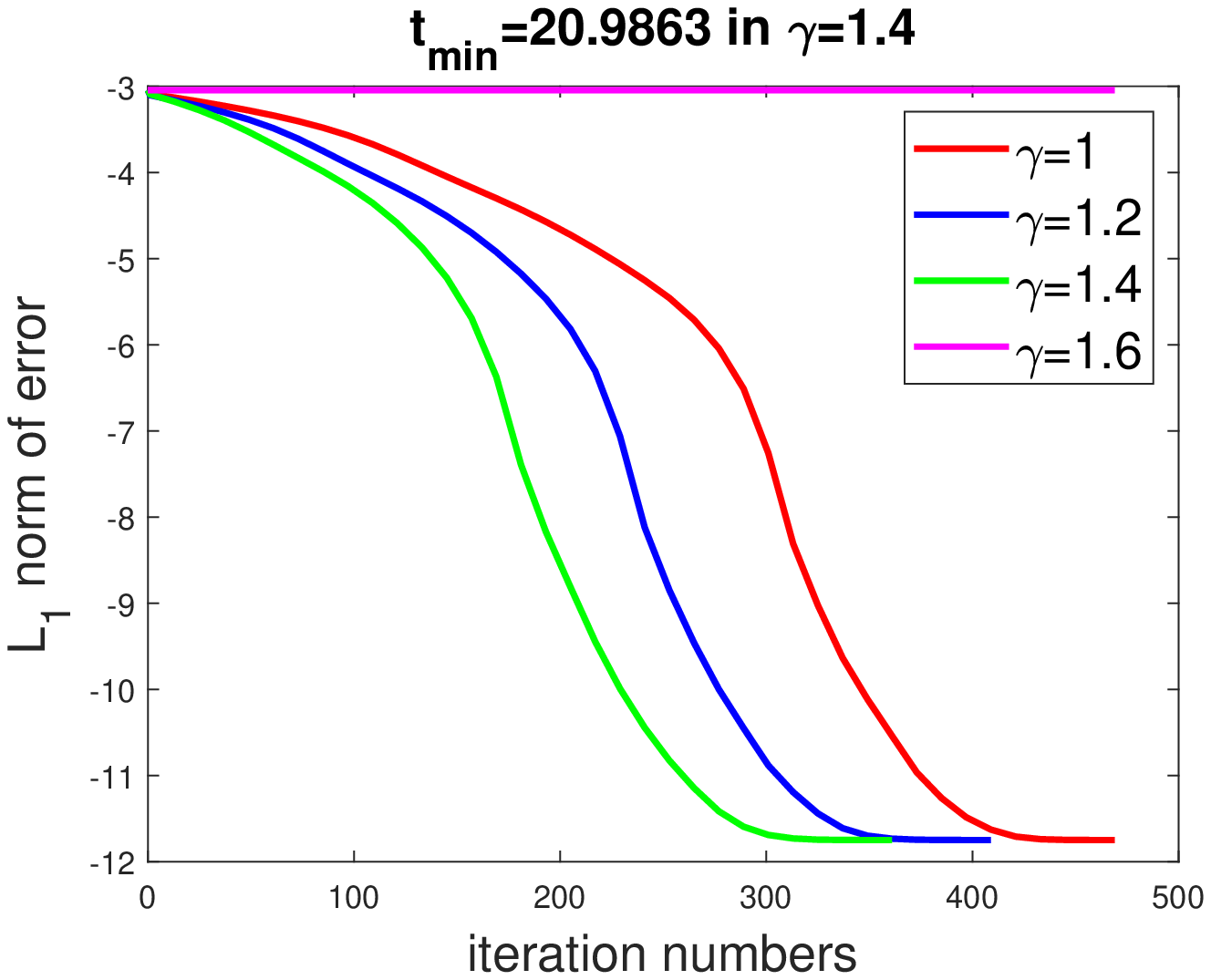}
}
\subfigure[FE-FSM-h-Convergence.]{
\includegraphics[width=5.1cm]{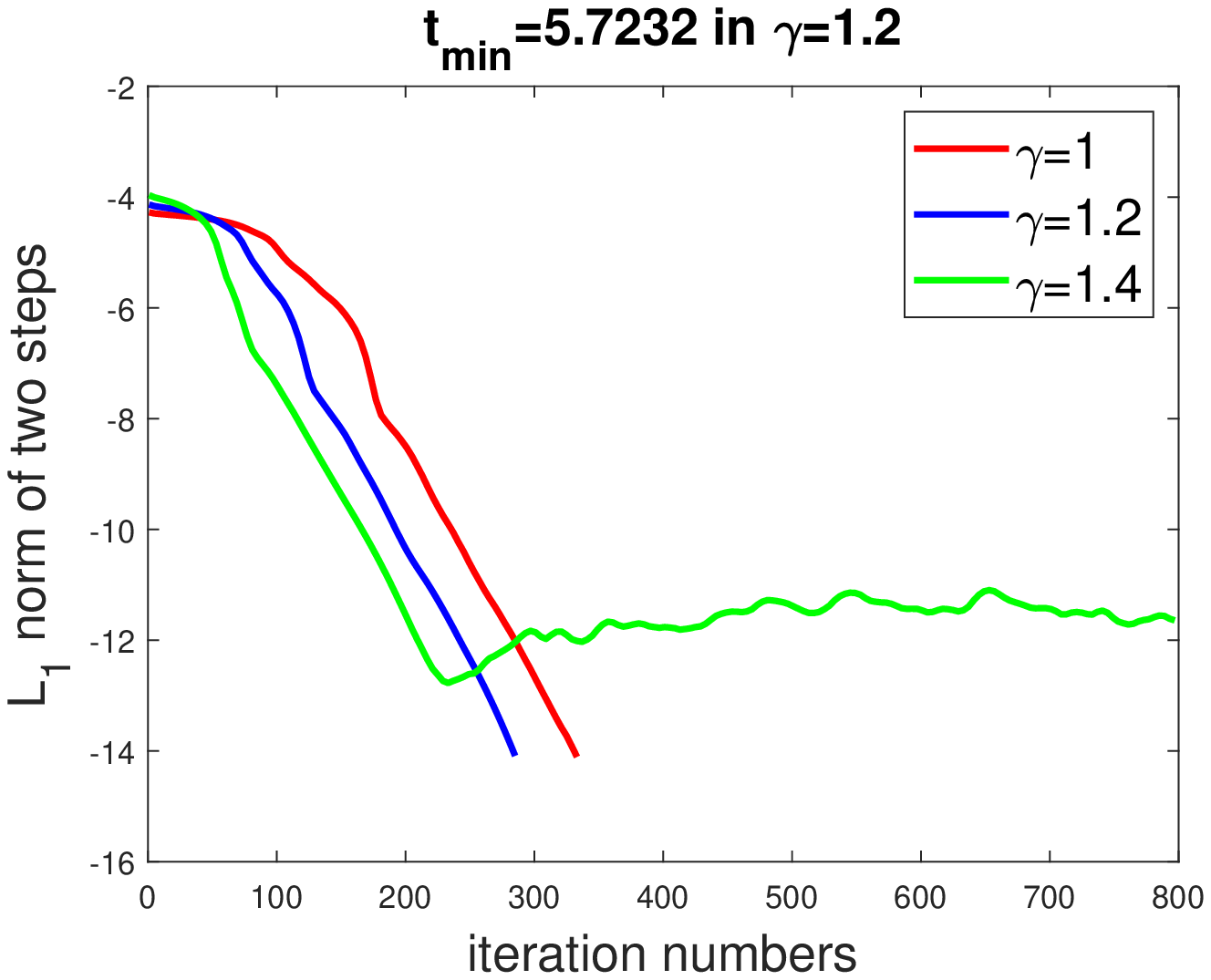}
}
\subfigure[FE-FSM-h-Residual.]{
\includegraphics[width=5.1cm]{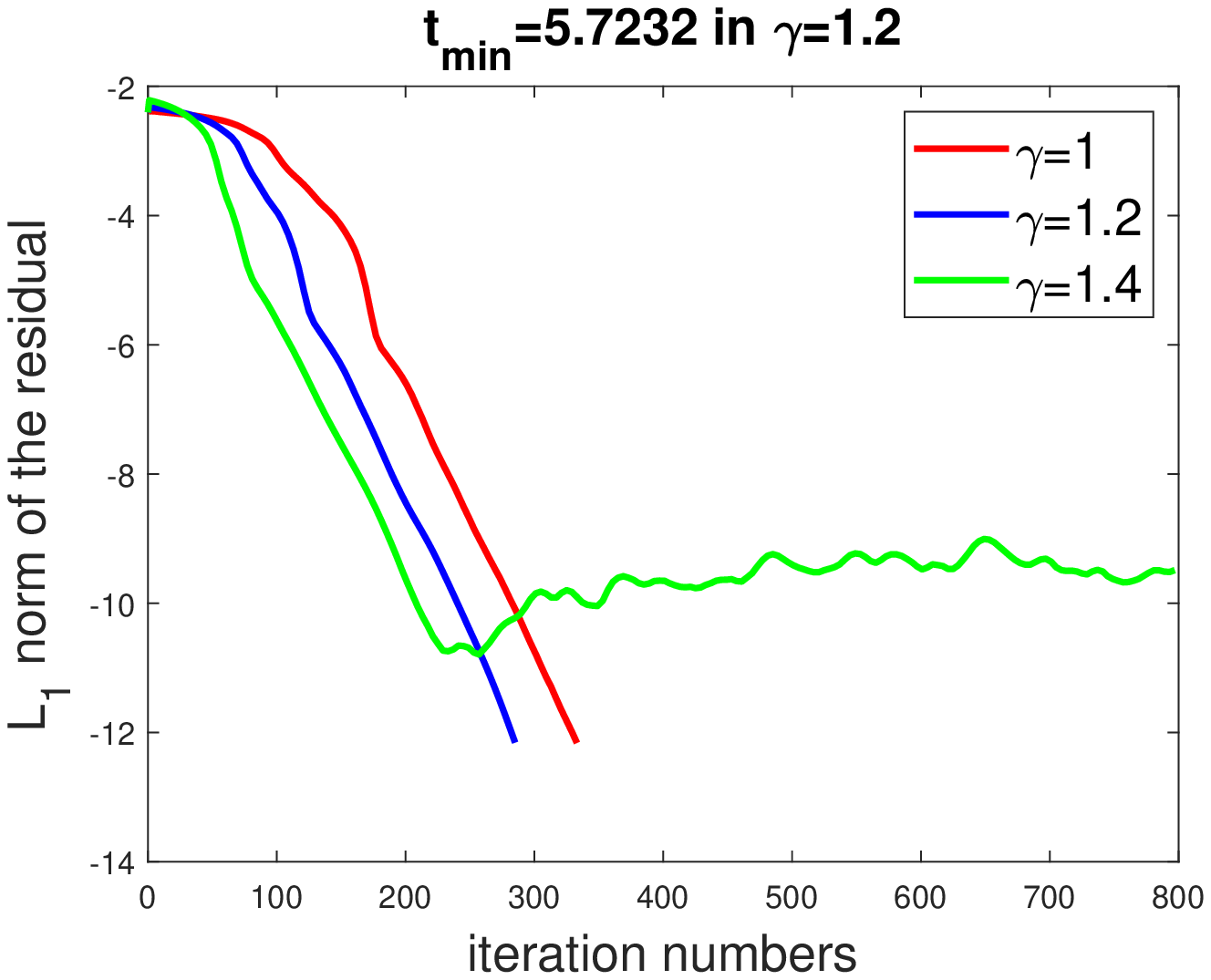}
}
\subfigure[FE-FSM-h-Numerical error.]{
\includegraphics[width=5.1cm]{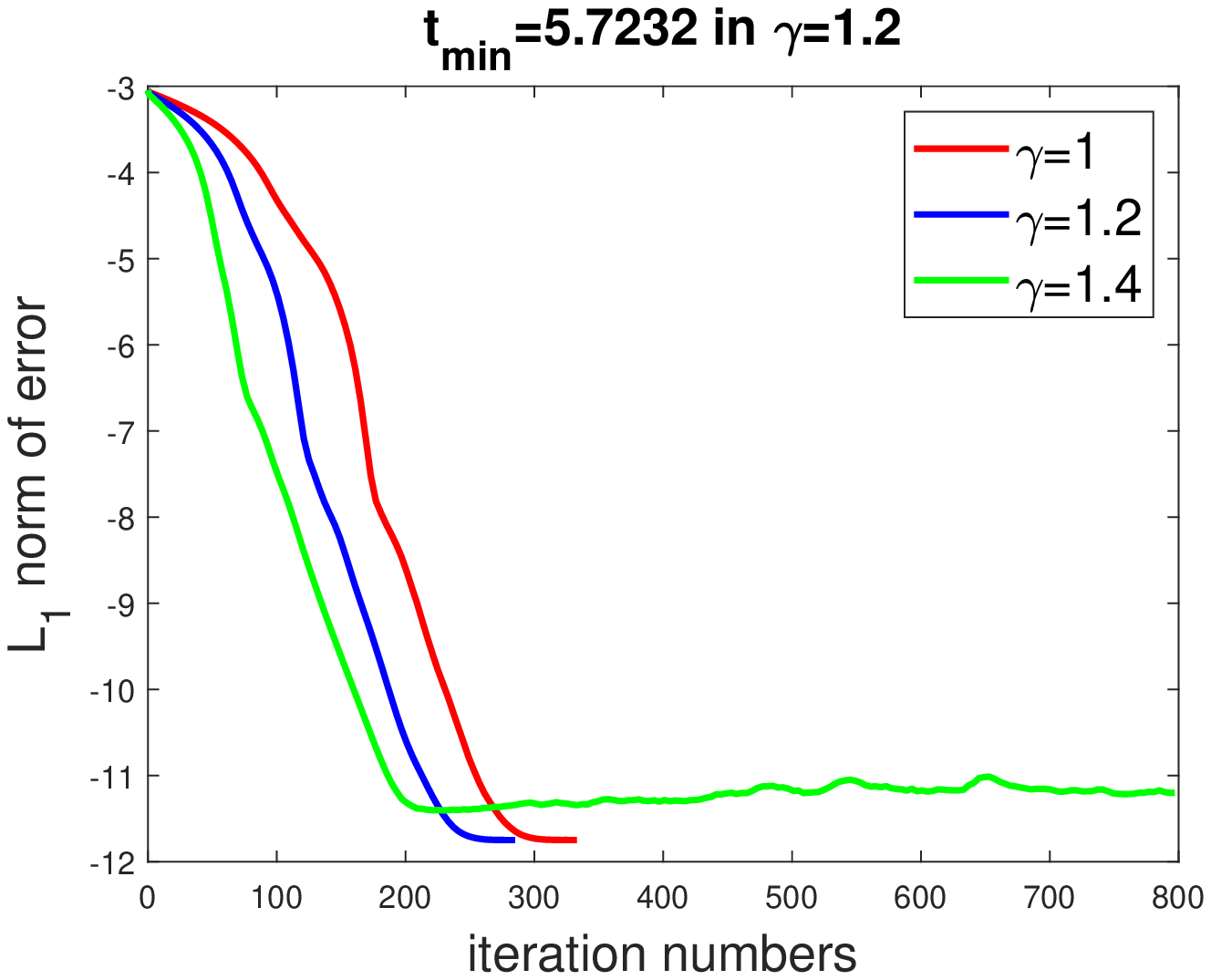}
}
\subfigure[RK-FSM-h-Convergence.]{
\includegraphics[width=5.1cm]{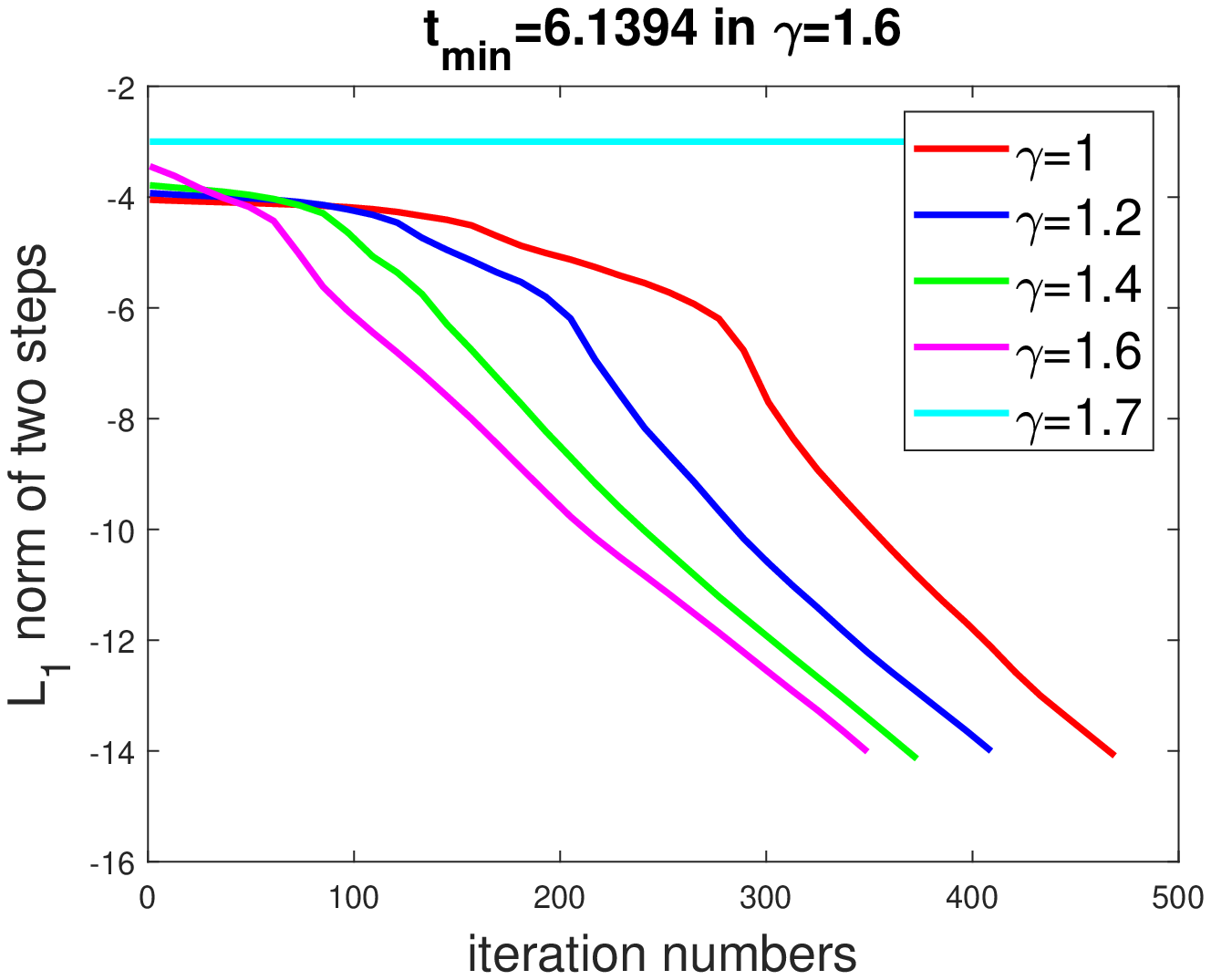}
}
\subfigure[RK-FSM-h-Residual.]{
\includegraphics[width=5.1cm]{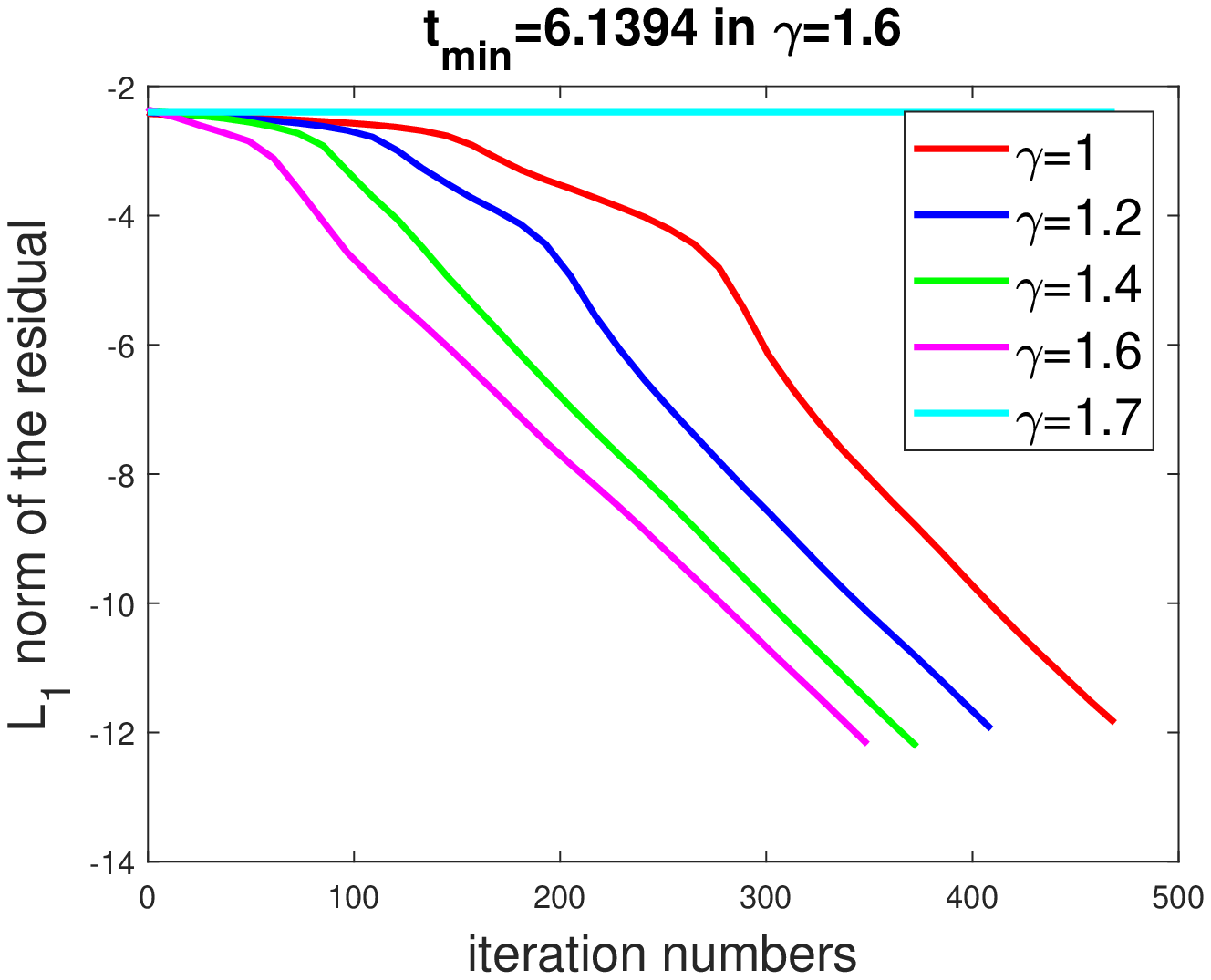}
}
\subfigure[RK-FSM-h-Numerical error.]{
\includegraphics[width=5.1cm]{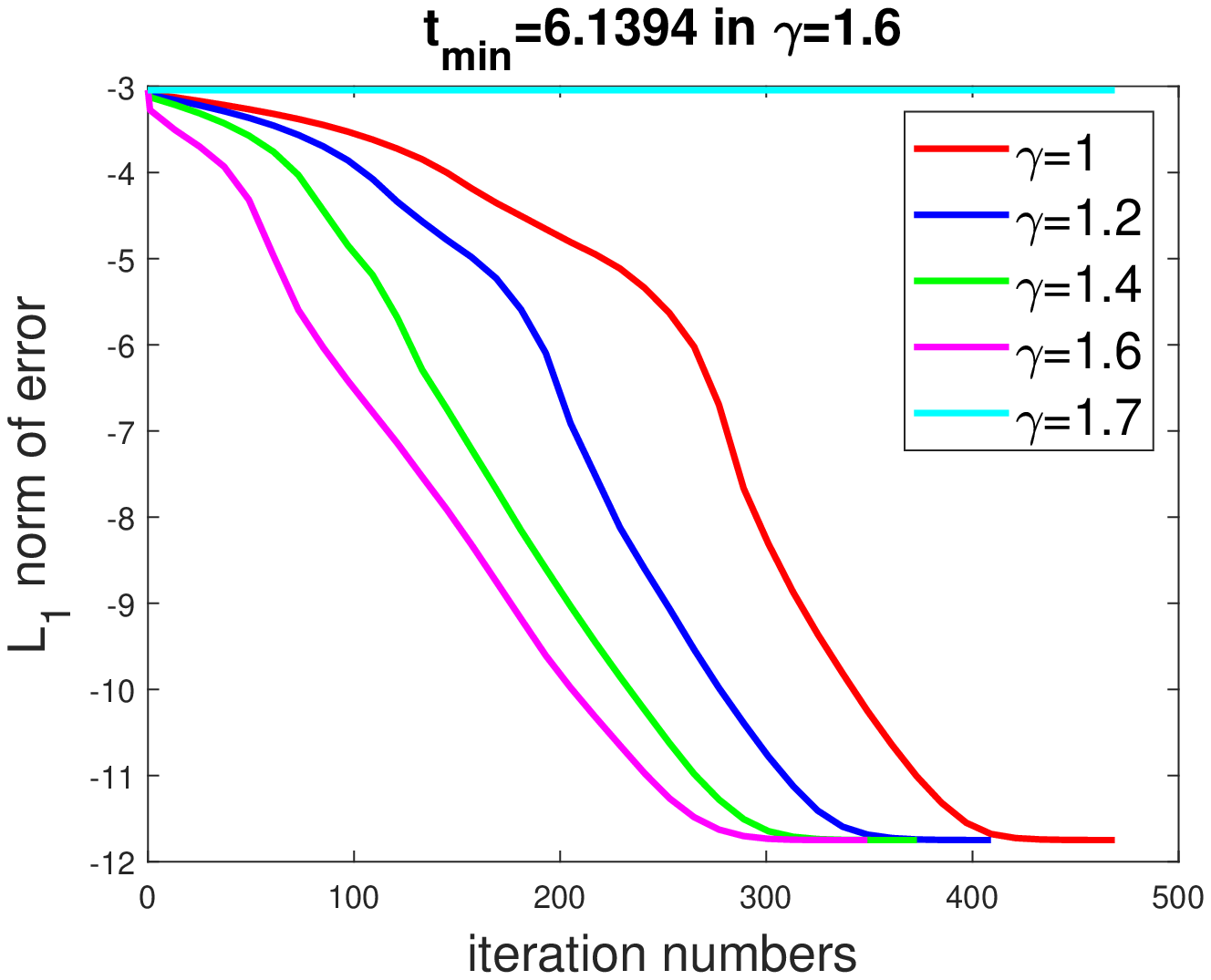}
}
\caption{Example 2. Study on different $\gamma's$}\label{fige2}
\end{figure}

Again, we want to remark that, for the three schemes other than the FE-Jacobi scheme, the CFL number can be taken to be greater than 1. Figure \ref{fige2} shows that the convergence history of the FE-FSM and RK-FSM with larger CFL numbers on mesh $N=320$, reporting the errors between two consecutive iteration steps, the residual and numerical errors. Similar to Example 1, we have presented the optimal CFL number and CPU time of each scheme in the title of each sub-figure. With larger CFL number, fewer iterations are needed for convergence. We observe that the RK-FSM tends to admit a larger CFL number than the FE-FSM, but it still underperforms in terms of the computational time.

In summary, the FE-FSM demonstrates to yield the best performance out of these four methods, which is consistent with the observation from Example 1. Also, the hybrid strategy can further reduce the computational cost.

\begin{table}
\caption{Example 2. Comparison of the four methods: The errors of the numerical solution, the accuracy obtained and the number of iterations for convergence}\label{tab2}
	\begin{center}
\resizebox{\textwidth}{55mm}{
		\begin{tabular}{|c|cccccc|cccccc|}
			\hline
\multicolumn{7}{|c|}{FE-Jacobi $\gamma=0.1$}&\multicolumn{6}{|c|}{FE-Jacobi $\gamma=0.1$ with hybrid strategy}\\ \hline
N & $L_{1}$&order&  $L_{\infty}$&order&iter &time& $L_{1}$&order&  $L_{\infty}$&order&iter&time\\ \hline
40 &6.09e-07 &- &2.33e-05 &- &1121 &0.9706&6.09e-07 &- &2.33e-05 &- &1106 &0.8406\\
80 &1.16e-08 &5.70 &1.74e-06 &3.74 &1553 &5.3702&1.16e-08 &5.70 &1.74e-06 &3.74 &1541 &3.1265\\
160 &8.70e-11 &7.06 &2.20e-08 &6.30 &2260 &32.6669&1.04e-10 &6.80 &3.21e-08 &5.76 &2207 &15.6215\\
320 &- &- &- &- &- &-&- &- &- &- &- &-\\ \hline
\multicolumn{7}{|c|}{FE-FSM $\gamma=1$}&\multicolumn{6}{|c|}{FE-FSM $\gamma=1$ with hybrid strategy}\\ \hline
$N$ & $L_{1}$~error&order& $L_{\infty}$~error&order&iter&time& $L_{1}$&order&  $L_{\infty}$&order&iter&time\\ \hline
 40 &6.75e-07 &- &5.41e-05 &- &160 &0.1337 &6.04e-07 &- &2.33e-05 &- &280 &0.2301\\
80 &1.16e-08 &5.85 &1.74e-06 &4.95 &208 &0.6374 &1.15e-08 &5.70 &1.74e-06 &3.74 &280 &0.46745\\
160 &8.76e-11 &7.05 &2.21e-08 &6.30 &244 &3.1495 &9.98e-11 &6.85 &4.55e-08 &5.26 &248 &1.2264\\
320 &1.78e-12 &5.61 &1.30e-10 &7.40 &324 &19.7232 &1.78e-12 &5.80 &1.30e-10 &8.44 &336 &5.6339\\ \hline
\multicolumn{7}{|c|}{RK-Jacobi $\gamma=1$}&\multicolumn{6}{|c|}{RK-Jacobi $\gamma=1$ with hybrid strategy}\\ \hline
$N$ & $L_{1}$~error&order& $L_{\infty}$~error&order&iter&time& $L_{1}$&order&  $L_{\infty}$&order&iter&time\\ \hline
 40 &9.00e-07 &- &6.26e-05 &- &282 &0.3819&8.95e-07 &- &6.25e-05 &- &282 &0.4601\\
80 &1.16e-08 &6.26 &1.74e-06 &5.16 &363 &1.2983&1.15e-08 &6.27 &1.74e-06 &5.16 &369 &1.0393\\
160 &8.76e-11 &7.05 &2.21e-08 &6.30 &513 &7.0839&8.77e-11 &7.03 &2.21e-08 &6.30 &516 &3.4253\\
320 &1.80e-12 &5.60 &1.30e-10 &7.40 &903 &54.8592&1.80e-12 &5.60 &1.30e-10 &7.40 &903 &26.9378\\ \hline
\multicolumn{7}{|c|}{RK-FSM $\gamma=1$}&\multicolumn{6}{|c|}{RK-FSM $\gamma=1$ with hybrid strategy}\\ \hline
$N$ & $L_{1}$~error&order& $L_{\infty}$~error&order&iter&time& $L_{1}$&order&  $L_{\infty}$&order&iter&time\\ \hline
40 &5.97e-07 &- &2.33e-05 &- &252 &0.2615 &5.96e-07 &- &2.33e-05 &- &252 &0.1610\\
80 &1.41e-08 &5.40 &4.14e-06 &2.49 &312 &0.9659 &1.17e-08 &5.66 &1.74e-06 &3.74 &312 &0.5905\\
160 &8.72e-11 &7.33 &2.20e-08 &7.55 &372 &4.8890 &9.28e-11 &6.98 &3.42e-08 &5.67 &372 &1.8319\\
320 &1.78e-12 &5.61 &1.30e-10 &7.39 &480 &27.7620 &1.78e-12 &5.70 &1.30e-10 &8.03 &480 &7.6022\\ \hline
\end{tabular}}
\end{center}
\end{table}
\newpage

\noindent \textbf{Example 3}. In this test case, we solve the Eikonal equation with $f(x,y)=1$. The computational domain is set as $[-3,3]^2$, and the inflow boundary $\Gamma$ consists of two circles of equal radius $0.5$ with the centers located at $(-1,0)$ and $(\sqrt{1.5},0)$, respectively, which leads to
$$\Gamma=\left\{(x,y)|(x+1)^2+y^2=\frac{1}{4}\quad or\quad (x-\sqrt{1.5})^2+y^2=\frac{1}{4}\right\}.$$
The exact solution is a distance function to the inflow boundary $\Gamma$, containing the singularities at the center of each circle and the line $x=0.5(\sqrt{1.5}-1)$ that is of equal distance to two circle centers.

Again, the Godunov numerical Hamiltonian is used. We measure the numerical errors within the
box of $[-2.85,2.85]^2$, which also excludes the boxes $[-1.15,-0.85]\times[-0.15,0.15]$, $[\sqrt{1.5}-0.15,\sqrt{1.5}+0.15]\times[-0.15,0.15]$
and $[\sqrt{0.375}-0.65,\sqrt{0.375}-0.35]\times[-2.85,2.85]$. These excluded boxes contain
two centers of $\Gamma$ and the singular line.

\begin{figure}[!h]
\begin{center}
	\includegraphics[width=6cm]{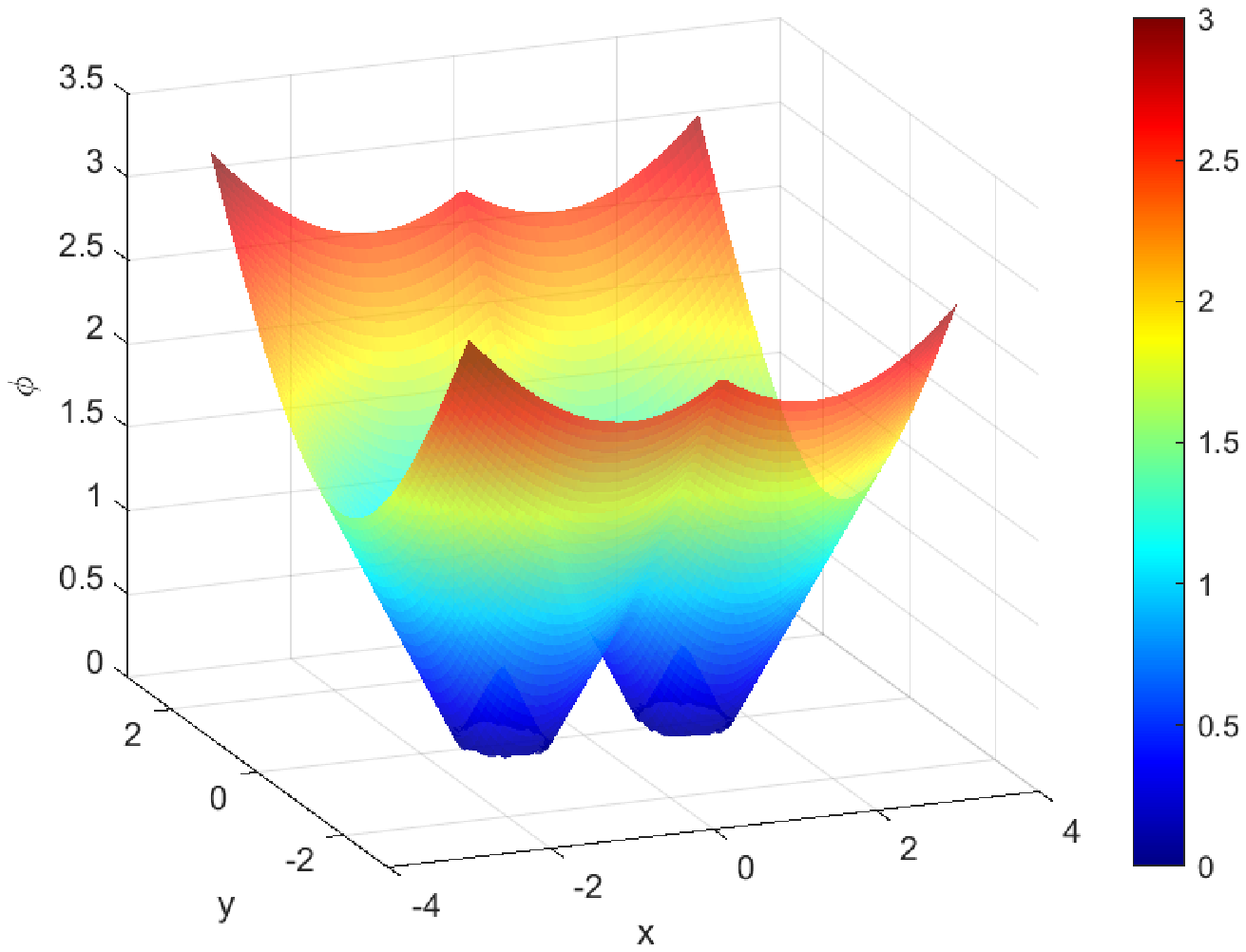}
	\includegraphics[width=6cm]{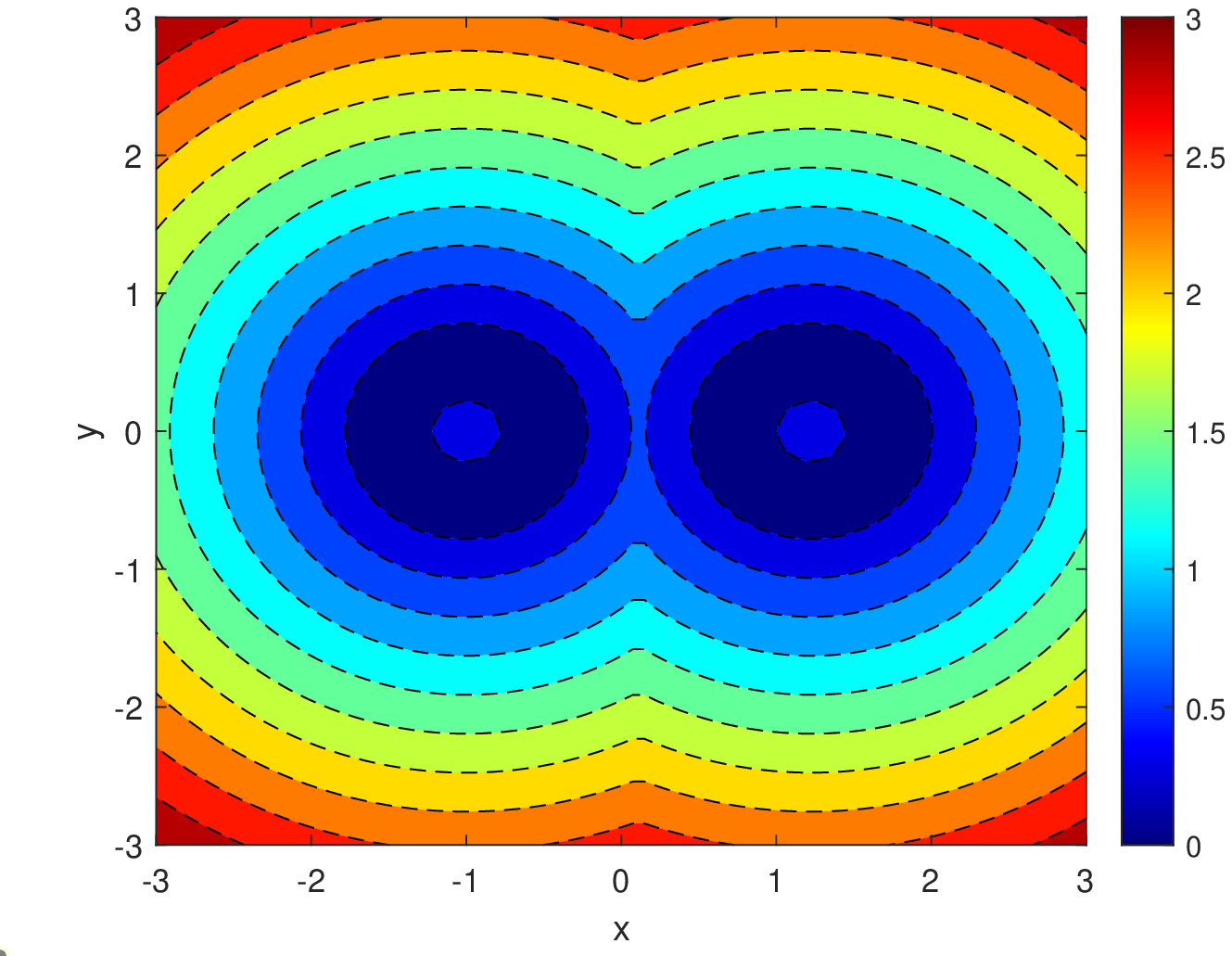}
\end{center}
\caption{Example 3. The numerical solution by FE-FSM on mesh $N=80$. Left: the 3D plot of numerical solution $\phi$; Right: the contour plot for $\phi$ .}\label{fige3sc}
\end{figure}

Figure \ref{fige3sc} shows that the numerical solution by FE-FSM.
The numerical results with hybrid
strategy and without hybrid strategy are reported on the Table
\ref{tab3}. Again, We can observe that the FE-Jacobi method requires the smaller CFL number of value 0.1, but still fails to converge on the refined mesh. The other three schemes can take a larger CFL number and the fast sweeping technique can improve the convergence of the Jacobi scheme. On the same refined mesh, we can see that
the RK-FSM only takes about $50\%$ CPU time of the RK-Jacobi scheme. Furthermore, the FE-FSM costs even less CPU time than RK-FSM. When the hybrid strategy is used, all four schemes can save about $75\%$ CPU time on the refined mesh.

In addition, we want to remark that, for the three schemes other than the FE-Jacobi scheme, the CFL number can be taken to be greater than 1.
With larger CFL number, fewer iterations are needed for convergence. This observation is exactly the same as the previous example, hence we will no longer report the convergence history with different CFL numbers to save space.

\begin{table}
\caption{Example 3. Comparison of the four methods: The errors of the numerical solution, the accuracy obtained and the number of iterations for convergence}\label{tab3}
	\begin{center}
\resizebox{\textwidth}{47mm}{
		\begin{tabular}{|c|cccccc|cccccc|}
			\hline
\multicolumn{7}{|c|}{FE-Jacobi $\gamma=0.1$}&\multicolumn{6}{|c|}{FE-Jacobi $\gamma=0.1$ with hybrid strategy}\\ \hline
N & $L_{1}$&order&  $L_{\infty}$&order&iter &time& $L_{1}$&order&  $L_{\infty}$&order&iter&time\\ \hline
80 &7.09e-07 &- &1.17e-04 &- &1884 &7.6586&7.09e-07 &- &1.17e-04 &- &1867 &3.3288\\
160 &1.10e-07 &2.68 &9.25e-06 &3.67 &2859 &49.5274&1.94e-07 &1.86 &1.81e-05 &2.70 &2852 &14.9261\\
320 &- &- &- &- &- &-&- &- &- &- &- &-\\ \hline
\multicolumn{7}{|c|}{FE-FSM $\gamma=1$}&\multicolumn{6}{|c|}{FE-FSM $\gamma=1$ with hybrid strategy}\\ \hline
$N$ & $L_{1}$~error&order& $L_{\infty}$~error&order&iter&time& $L_{1}$&order&  $L_{\infty}$&order&iter&time\\ \hline
80 &7.02e-07 &- &1.17e-04 &- &236 &0.88271&7.03e-07 &- &1.17e-04 &- &236 &0.34559\\
160 &1.10e-07 &2.66 &9.71e-06 &3.59 &292 &4.599&2.01e-07 &1.80 &1.81e-05 &2.69 &288 &1.3335\\
320 &2.98e-10 &8.53 &3.86e-07 &4.65 &408 &26.4953 &3.18e-10 &9.30 &4.03e-07 &5.48 &368 &6.1812\\ \hline
\multicolumn{7}{|c|}{RK-Jacobi $\gamma=1$}&\multicolumn{6}{|c|}{RK-Jacobi $\gamma=1$ with hybrid strategy}\\ \hline
$N$ & $L_{1}$~error&order& $L_{\infty}$~error&order&iter&time& $L_{1}$&order&  $L_{\infty}$&order&iter&time\\ \hline
 80 &7.09e-07 &- &1.18e-04 &- &426 &1.6574&7.09e-07 &- &1.18e-04 &- &423 &0.65345\\
160 &1.19e-07 &2.57 &1.75e-05 &2.74 &681 &11.6017&2.05e-07 &1.78 &1.81e-05 &2.70 &678 &3.5579\\
320 &3.54e-10 &8.39 &3.74e-07 &5.55 &1197 &83.9635&3.86e-10 &9.05 &4.64e-07 &5.28 &1194 &24.711\\ \hline
\multicolumn{7}{|c|}{RK-FSM $\gamma=1$}&\multicolumn{6}{|c|}{RK-FSM $\gamma=1$ with hybrid strategy}\\ \hline
$N$ & $L_{1}$~error&order& $L_{\infty}$~error&order&iter&time& $L_{1}$&order&  $L_{\infty}$&order&iter&time\\ \hline
80 &7.10e-07 &- &1.01e-04 &- &360 &1.4243&7.11e-07 &- &1.01e-04 &- &360 &0.53576\\
160 &1.09e-07 &2.70 &9.53e-06 &3.41 &480 &9.2416&1.94e-07 &1.87 &1.81e-05 &2.48 &468 &2.4039\\
320 &3.01e-10 &8.50 &3.85e-07 &4.63 &648 &41.7686&3.16e-10 &9.26 &3.84e-07 &5.55 &624 &11.6644\\ \hline
\end{tabular}}
\end{center}
\end{table}

\noindent \textbf{Example 4}. Here we again consider the Eikonal equation with $f(x,y)=1$. The computational
domain is set as $[-1,1]^2$, and the inflow boundary is given by $\Gamma=(0,0)$. The
exact solution for this problem is a distance function to $\Gamma$,
and it contains a singularity at $\Gamma$.

The Godunov numerical Hamiltonian is used. Due to the singularity, we follow the setup in \cite{xiong}, and pre-assign the exact solution in a small box with
length $0.3$ around the source point. Numerical errors and orders are listed in Table \ref{tab4}. The same behavior as in the previous examples can be observed,
namely, a smaller CFL number is needed for FE-Jacobi method, and the fast sweeping technique can improve the convergence of the Jacobi scheme. Again, the FE-FSM performs the best out of these four methods, and the hybrid strategy can further reduce the computational cost.

\begin{table}
\caption{Example 4. Comparison of the four methods: The errors of the numerical solution, the accuracy obtained and the number of iterations for convergence}\label{tab4}
	\begin{center}
\resizebox{\textwidth}{55mm}{
		\begin{tabular}{|c|cccccc|cccccc|}
			\hline
\multicolumn{7}{|c|}{FE-Jacobi $\gamma=0.1$}&\multicolumn{6}{|c|}{FE-Jacobi $\gamma=0.1$ with hybrid strategy}\\ \hline
N & $L_{1}$&order&  $L_{\infty}$&order&iter &time& $L_{1}$&order&  $L_{\infty}$&order&iter&time\\ \hline
40 &1.29e-07 &- &4.17e-07 &- &1354 &1.4606&3.10e-07 &- &4.60e-06 &- &1281 &0.85274\\
80 &3.58e-09 &5.17 &1.07e-08 &5.27 &1910 &7.7459&6.95e-09 &5.48 &1.62e-07 &4.82 &1765 &2.8821\\
160 &1.06e-10 &5.06 &2.85e-10 &5.24 &2707 &44.3751&- &- &- &- &- &-\\
320 &- &- &- &- &- &-&- &- &- &- &- &-\\ \hline
\multicolumn{7}{|c|}{FE-FSM $\gamma=1$}&\multicolumn{6}{|c|}{FE-FSM $\gamma=1$ with hybrid strategy}\\ \hline
$N$ & $L_{1}$~error&order& $L_{\infty}$~error&order&iter&time& $L_{1}$&order&  $L_{\infty}$&order&iter&time\\ \hline
40 &3.10e-07 &- &4.60e-06 &- &192 &0.19561&3.10e-07 &- &4.60e-06 &- &292 &0.16633\\
80 &6.95e-09 &5.48 &1.62e-07 &4.82 &240 &0.91323&6.95e-09 &5.48 &1.62e-07 &4.82 &336 &0.49803\\
160 &1.19e-10 &5.85 &1.78e-09 &6.51 &296 &4.466&1.19e-10 &5.85 &1.79e-09 &6.50 &340 &1.5533\\
320 &3.26e-12 &5.20 &8.67e-12 &7.68 &412 &22.457&3.26e-12 &5.20 &8.67e-12 &7.69 &380 &5.9912\\ \hline
\multicolumn{7}{|c|}{RK-Jacobi $\gamma=1$}&\multicolumn{6}{|c|}{RK-Jacobi $\gamma=1$ with hybrid strategy}\\ \hline
$N$ & $L_{1}$~error&order& $L_{\infty}$~error&order&iter&time& $L_{1}$&order&  $L_{\infty}$&order&iter&time\\ \hline
 40 &1.29e-07 &- &4.27e-07 &- &339 &0.35881&3.10e-07 &- &4.60e-06 &- &345 &0.23766\\
80 &3.58e-09 &5.17 &1.07e-08 &5.31 &441 &1.8265&6.95e-09 &5.48 &1.62e-07 &4.82 &444 &0.75318\\
160 &1.06e-10 &5.06 &2.85e-10 &5.23 &615 &9.8505&1.19e-10 &5.85 &1.78e-09 &6.51 &621 &3.3396\\
320 &3.26e-12 &5.03 &8.54e-12 &5.06 &1047 &69.828&3.26e-12 &5.20 &8.66e-12 &7.68 &1053 &21.902\\ \hline
\multicolumn{7}{|c|}{RK-FSM $\gamma=1$}&\multicolumn{6}{|c|}{RK-FSM $\gamma=1$ with hybrid strategy}\\ \hline
$N$ & $L_{1}$~error&order& $L_{\infty}$~error&order&iter&time& $L_{1}$&order&  $L_{\infty}$&order&iter&time\\ \hline
40 &3.10e-07 &- &4.60e-06 &- &300 &0.30303&3.10e-07 &- &4.60e-06 &- &300 &0.19219\\
80 &6.95e-09 &5.48 &1.62e-07 &4.82 &372 &1.4035 &6.95e-09 &5.48 &1.62e-07 &4.82 &372 &0.56225\\
160 &1.19e-10 &5.85 &1.78e-09 &6.51 &468 &7.0361&1.19e-10 &5.85 &1.78e-09 &6.51 &468 &2.2295\\
320 &3.26e-12 &5.20 &8.67e-12 &7.68 &612 &32.526&3.26e-12 &5.20 &8.67e-12 &7.68 &624 &10.4386\\ \hline
\end{tabular}}
\end{center}
\end{table}

\noindent\textbf{Example 5}. Consider the Eikonal equation with $f(x,y)=1$ on the computational domain $[-1,1]^2$. The inflow boundary
$\Gamma$ is a sector of three quarters of the circle centered at $(0,0)$ with radius $0.5$, closed with the x-axis and y-axis in the first quadrant, which can be described as
$$\Gamma=\left\{(x,y):\sqrt{x^2+y^2=0.5},~\mathrm{if}~x<0 \text{ or } y<0\}\cup\{(x,0):0\leq x\leq0.5\}\cup\{(0,y):0\leq y\leq0.5\right\}.$$
The exact solution is still the distance function to $\Gamma$. Singularities appear at the two corners in $\Gamma$, which give rise to both shock and rarefaction wave in the solution.

\begin{figure}[!h]
\begin{center}
	\includegraphics[width=6.cm]{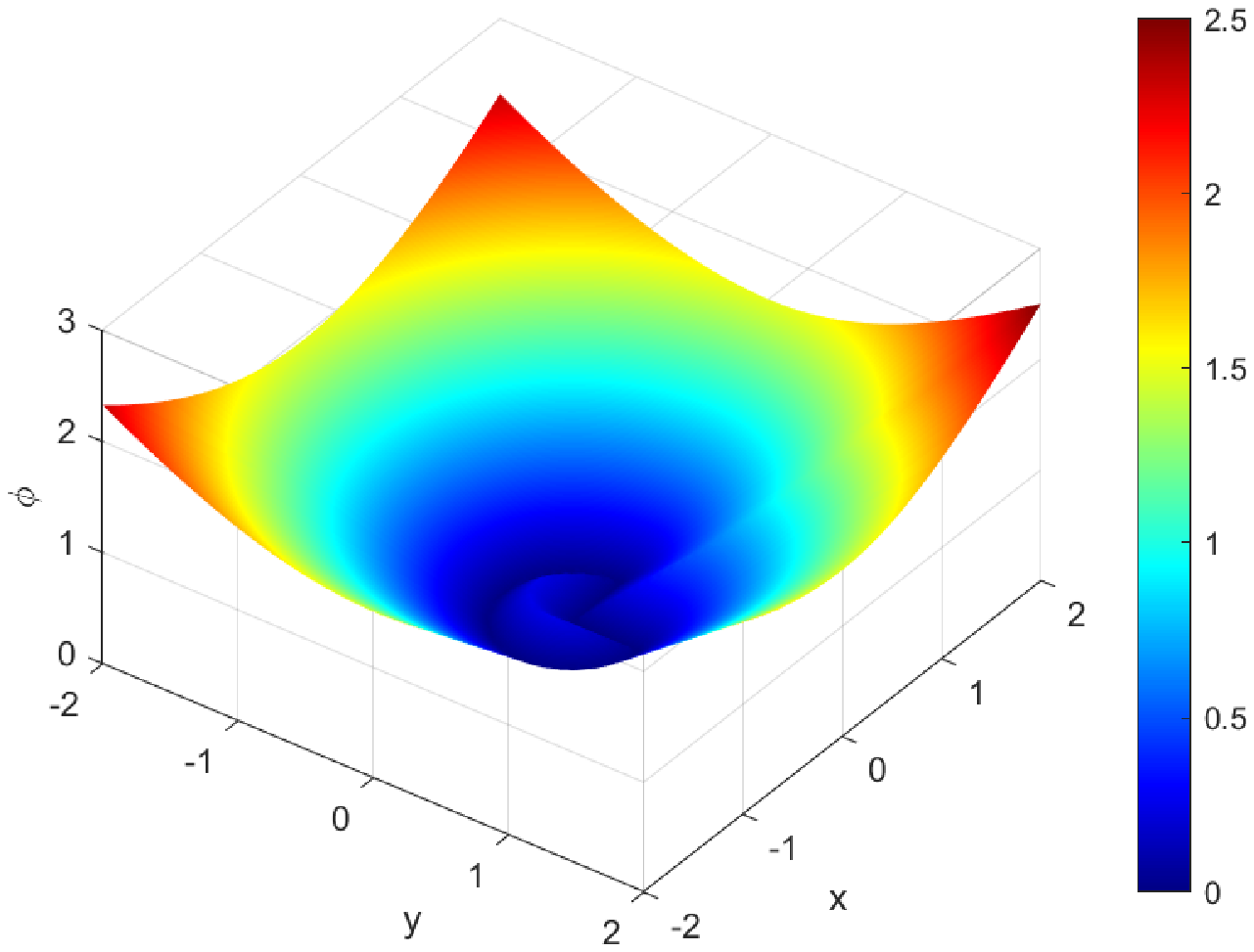}
	\includegraphics[width=6.cm]{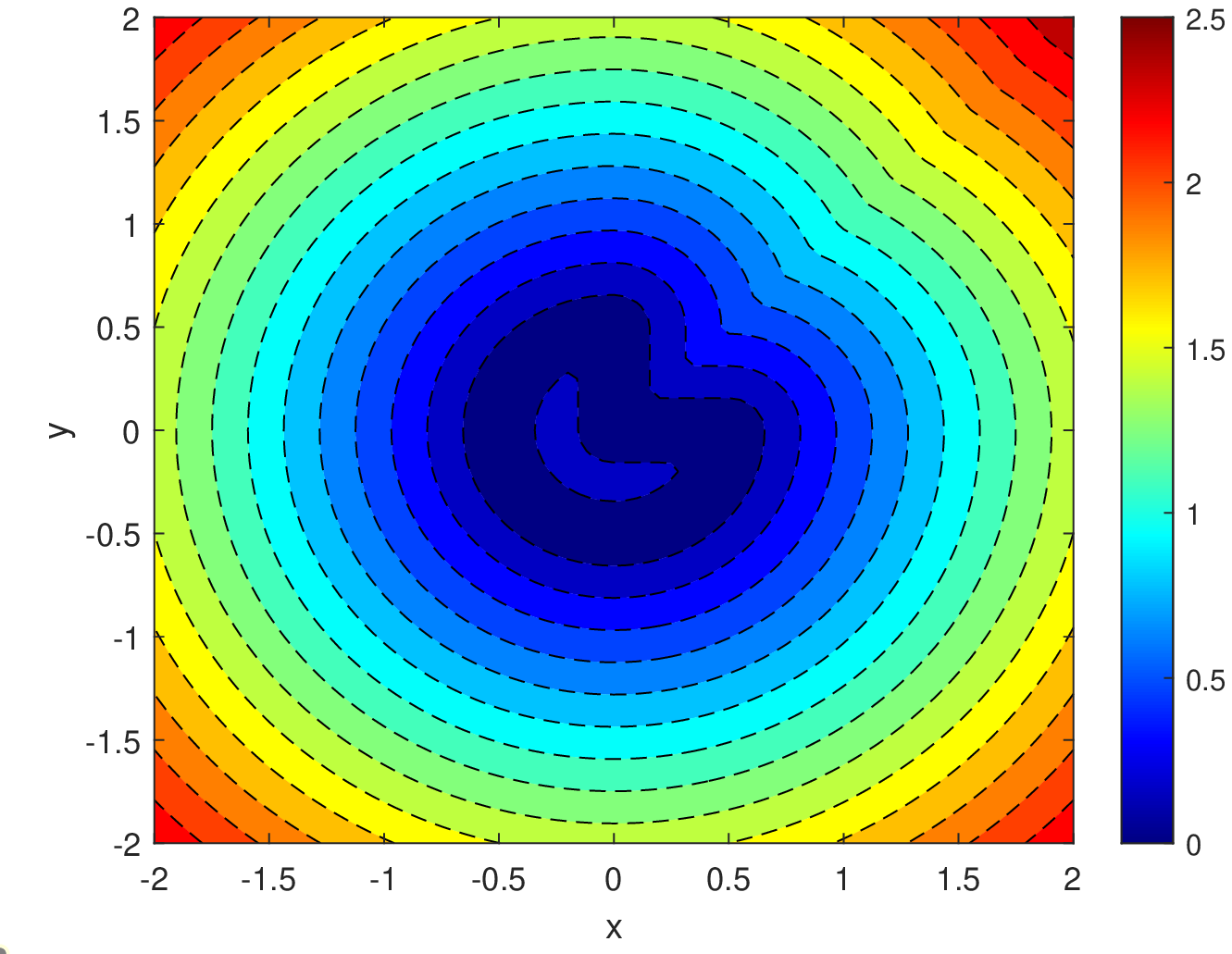}
\end{center}
\caption{Example 5. The numerical solution by FE-FSM on mesh $N=80$. Left: the 3D plot of numerical solution $\phi$; Right: the contour plot for $\phi$ .}\label{fige5sc}
\end{figure}

The Godunov numerical Hamiltonian is used. We measure the errors in smooth regions inside the box of $[-1.9,1.9]^2$ with $x\leq 0$ or $y\leq 0$, and outside the box $[-0.5,0.5]^2$. The surface and contour of the numerical solution by FE-FSM are shown in Figure \ref{fige5sc}.
The numerical errors and orders of convergence are shown in Table \ref{tab5}. Again, we observe that
the FE-Jacobi method requires the smaller CFL number of value 0.1, the other three schemes can take
a lager CFL number, and the fast sweeping technique can improve the convergence of the Jacobi scheme.
On the same refined mesh, we can see that the RK-FSM only takes about the $50\%$ CPU time of the RK-Jacobi scheme. Furthermore, the FE-FSM costs
even less CPU time than RK-FSM. The numerical results with the hybrid strategy can be seen on the right side
of Table \ref{tab5}, which suggests that the hybrid strategy can save 60\%-75\% CPU time on refined mesh. As the previous examples, the FE-FSM is more efficient scheme for the example, and the hybrid strategy can further reduce the computational cost.

\begin{table}
\caption{Example 5. Comparison of the four methods: The errors of the numerical solution, the accuracy obtained and the number of iterations for convergence}\label{tab5}
	\begin{center}
\resizebox{\textwidth}{55mm}{
		\begin{tabular}{|c|cccccc|cccccc|}
			\hline
\multicolumn{7}{|c|}{FE-Jacobi $\gamma=0.1$}&\multicolumn{6}{|c|}{FE-Jacobi $\gamma=0.1$ with hybrid strategy}\\ \hline
N & $L_{1}$&order&  $L_{\infty}$&order&iter &time& $L_{1}$&order&  $L_{\infty}$&order&iter&time\\ \hline
40 &1.12e-06 &- &1.61e-05 &- &1338 &0.81886&1.12e-06 &- &1.61e-05 &- &1338 &0.72088\\
80 &4.28e-08 &4.71 &1.08e-06 &3.89 &1978 &5.7302&4.28e-08 &4.71 &1.08e-06 &3.89 &1978 &2.8225\\
160 &1.05e-09 &5.33 &3.93e-08 &4.78 &3182 &41.18&1.05e-09 &5.33 &3.93e-08 &4.78 &3183 &14.5589\\
320 &2.11e-11 &5.64 &4.44e-10 &6.46 &5530 &323.1253 &- &- &- &- &- &-\\ \hline
\multicolumn{7}{|c|}{FE-FSM $\gamma=1$}&\multicolumn{6}{|c|}{FE-FSM $\gamma=1$ with hybrid strategy}\\ \hline
$N$ & $L_{1}$~error&order& $L_{\infty}$~error&order&iter&time& $L_{1}$&order&  $L_{\infty}$&order&iter&time\\ \hline
40 &1.12e-06 &- &1.61e-05 &- &236 &0.13161&1.12e-06 &- &1.61e-05 &- &236 &0.11094\\
80 &4.28e-08 &4.71 &1.08e-06 &3.89 &264 &0.72234&4.28e-08 &4.71 &1.08e-06 &3.89 &244 &0.30641\\
160 &1.05e-09 &5.33 &3.93e-08 &4.78 &296 &3.5751 &1.05e-09 &5.33 &3.93e-08 &4.78 &292 &1.1575\\
320 &2.11e-11 &5.64 &4.44e-10 &6.46 &444 &25.457&2.11e-11 &5.64 &4.44e-10 &6.46 &396 &6.0395\\ \hline
\multicolumn{7}{|c|}{RK-Jacobi $\gamma=1$}&\multicolumn{6}{|c|}{RK-Jacobi $\gamma=1$ with hybrid strategy}\\ \hline
$N$ & $L_{1}$~error&order& $L_{\infty}$~error&order&iter&time& $L_{1}$&order&  $L_{\infty}$&order&iter&time\\ \hline
40 &1.12e-06 &- &1.61e-05 &- &327 &0.23457&1.12e-06 &- &1.61e-05 &- &327 &0.16019\\
80 &4.28e-08 &4.71 &1.08e-06 &3.89 &438 &1.3605&4.28e-08 &4.710 &1.08e-06 &3.89 &438 &0.68153\\
160 &1.05e-09 &5.33 &3.93e-08 &4.78 &732 &10.2797&1.05e-09 &5.33 &3.93e-08 &4.78 &732 &3.3886\\
320 &2.11e-11 &5.64 &4.44e-10 &6.46 &1311 &75.9122&2.11e-11 &5.64 &4.44e-10 &6.46 &1311 &22.6901\\ \hline
\multicolumn{7}{|c|}{RK-FSM $\gamma=1$}&\multicolumn{6}{|c|}{RK-FSM $\gamma=1$ with hybrid strategy}\\ \hline
$N$ & $L_{1}$~error&order& $L_{\infty}$~error&order&iter&time& $L_{1}$&order&  $L_{\infty}$&order&iter&time\\ \hline
40 &1.12e-06 &- &1.61e-05 &- &288 &0.16363&1.12e-06 &- &1.61e-05 &- &300 &0.13274\\
80 &4.28e-08 &4.71 &1.08e-06 &3.89 &360 &1.0679&4.28e-08 &4.71 &1.08e-06 &3.89 &372 &0.44366\\
160 &1.05e-09 &5.33 &3.93e-08 &4.78 &468 &5.6948&1.05e-09 &5.33 &3.93e-08 &4.78 &480 &2.0111\\
320 &2.11e-11 &5.64 &4.44e-10 &6.46 &648 &36.3658&2.11e-11 &5.64 &4.44e-10 &6.46 &624 &8.9956\\ \hline
\end{tabular}}
\end{center}
\end{table}

 \noindent\textbf{Example 6}. In this test, we solve the Eikonal equation with
 $$f(x,y)=2\pi\sqrt{[\cos(2\pi x)\sin(2\pi y)]^2+[\sin(2\pi x)\cos(2\pi y)]^2}.$$
The computational domain is set as $\Omega=[0,1]^2$, and the inflow boundary condition is given by $\Gamma=\{(\frac{1}{4},\frac{1}{4}),(\frac{3}{4},\frac{3}{4}),(\frac{1}{4},\frac{3}{4}),(\frac{3}{4},\frac{1}{4}),(\frac{1}{2},\frac{1}{2})\},$ consisting of five isolated points.
$\phi(x,y)=0$ is prescribed at the boundary of the unit square. The exact solution of this problem is the shape function \cite{zhang2006}. Two cases are considered here, based on different boundary conditions.
\indent {Case (a):}
$$g\left(\frac{1}{4},\frac{1}{4}\right)=g\left(\frac{3}{4},\frac{3}{4}\right)=1,\quad g\left(\frac{1}{4},\frac{3}{4}\right)=g\left(\frac{3}{4},\frac{1}{4}\right)=-1,\quad
g\left(\frac{1}{2},\frac{1}{2}\right)=0,$$
with the exact solution being
$$\phi(x,y)=\sin(2\pi x)\sin(2\pi y);$$
\indent{Case (b):}
 $$g\left(\frac{1}{4},\frac{1}{4}\right)=g\left(\frac{3}{4},\frac{3}{4}\right)=g\left(\frac{1}{4},\frac{3}{4}\right)=g\left(\frac{3}{4},\frac{1}{4}\right)=1,\quad
g\left(\frac{1}{2},\frac{1}{2}\right)=2,$$
with the exact solution being
\begin{equation*}
  \phi(x,y)=\begin{cases}
\max(|\sin(2\pi x)\sin(2\pi y)|,1+\cos(2\pi x)\cos(2\pi y)),
&\mathrm{if}~|x+y-1|<\frac{1}{2}~\mathrm{and}~|x-y|<\frac{1}{2},\\
|\sin(2\pi x)\sin(2\pi y)|,&\mathrm{otherwise},
\end{cases}
\end{equation*}
which is not smooth.

Due to the singularity of these point sources, the exact solutions are placed in a small box with a length $2h$ around these isolated points in both test cases. The Godunov numerical Hamiltonian is used in this test.

\begin{figure}[!h]
\begin{center}
	\includegraphics[width=6.cm]{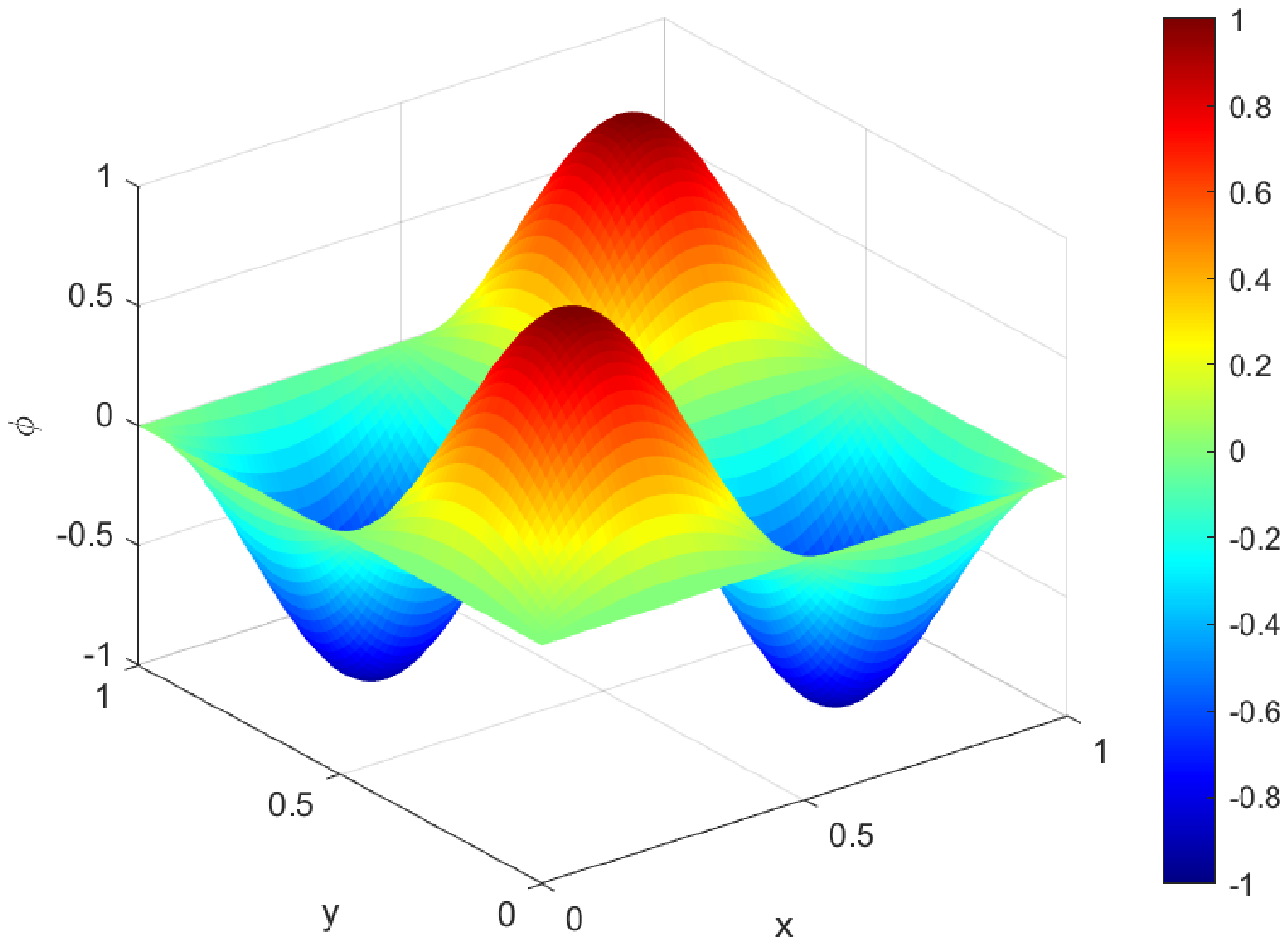}
	\includegraphics[width=6.cm]{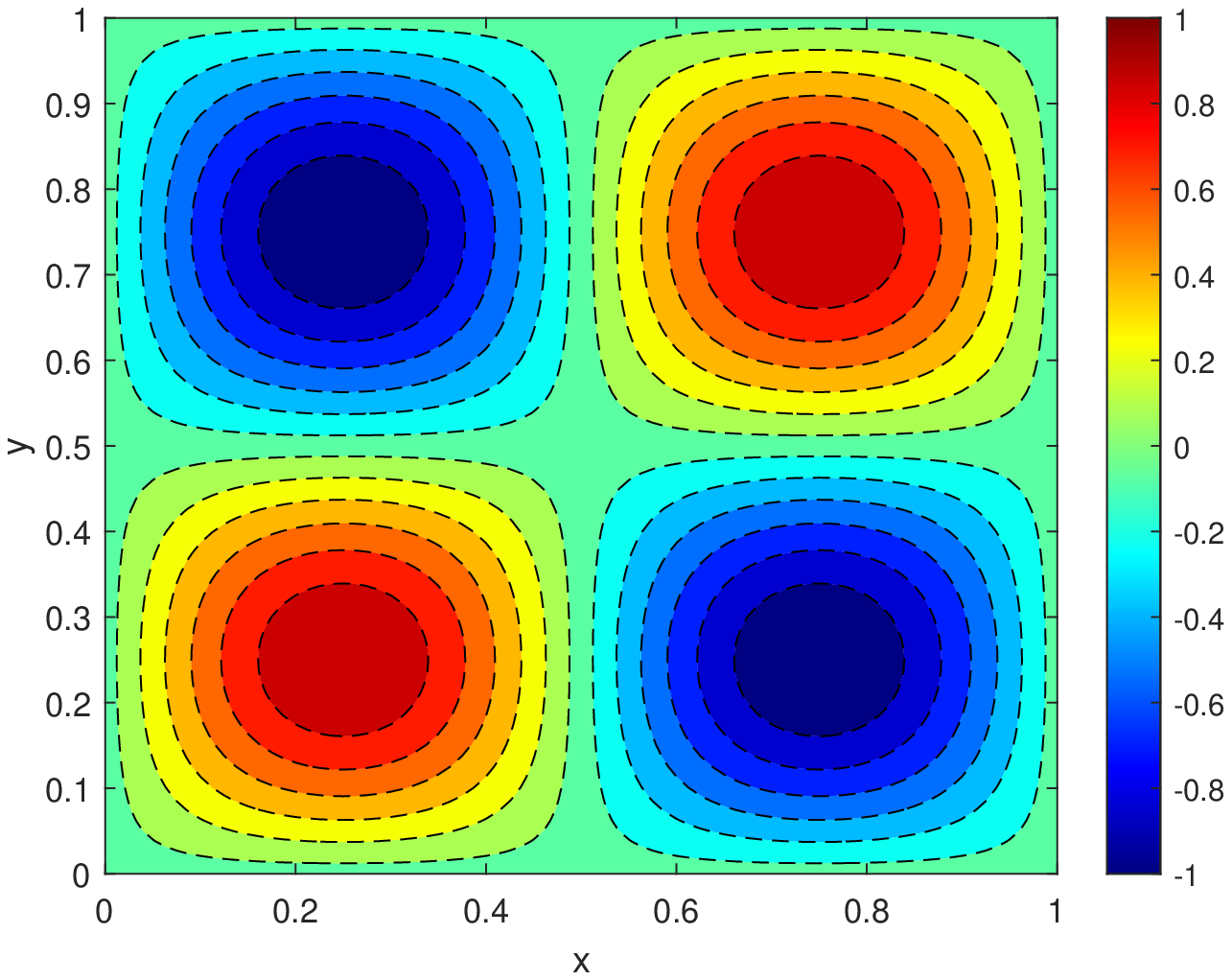}
\end{center}
\caption{Example 6 case (a). The numerical solution by FE-FSM on mesh $N=80$. Left: the 3D plot of numerical solution $\phi$; Right: the contour plot for $\phi$ .}\label{fige61sc}
\end{figure}

For the case (a), Figure \ref{fige61sc} shows the surface and contour of numerical solution by FE-FSM. The numerical
errors and orders of convergence of four methods are listed in Table \ref{tab6a}. We
can observe that the fifth order accuracy can be obtained, and the same behavior as the previous
examples can be observe for these four schemes.
On the same refined mesh, we can see that the RK-FSM only takes about $50\%$ CPU time of the RK-Jacobi scheme, furthermore, the FE-FSM costs even less CPU time than RK-FSM. In addition, when the fast sweeping technique is used, the FE-FSM can now use a larger CFL number than FE-Jacobi scheme. The numerical results after the hybrid strategy is used can be seen on the right side of Table \ref{tab6a}, which suggests that the hybrid strategy can save $50\%-75\%$ of the CPU time on the refined mesh.

\begin{figure}[!h]
\begin{center}
	\includegraphics[width=6.cm]{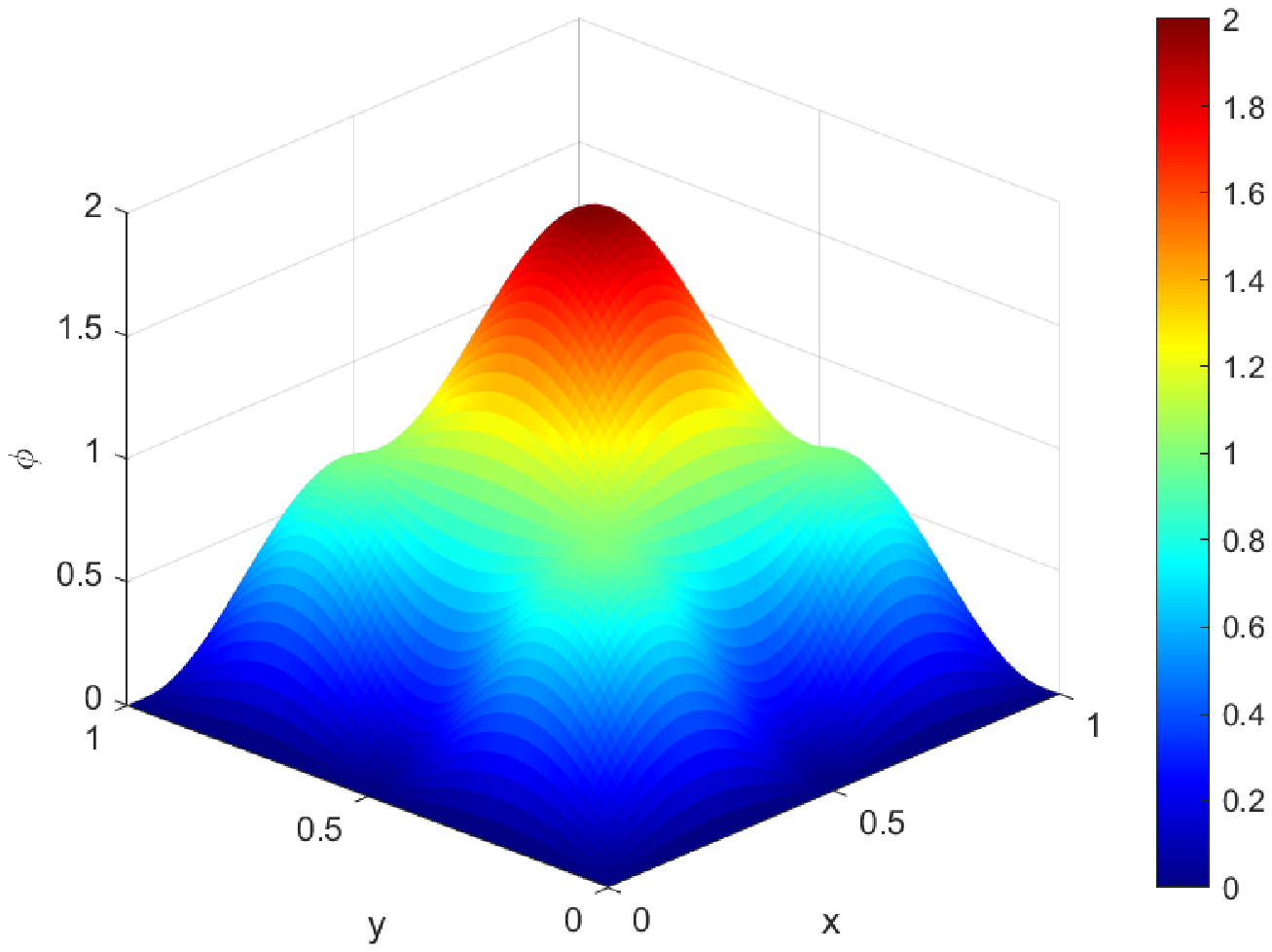}
	\includegraphics[width=6.cm]{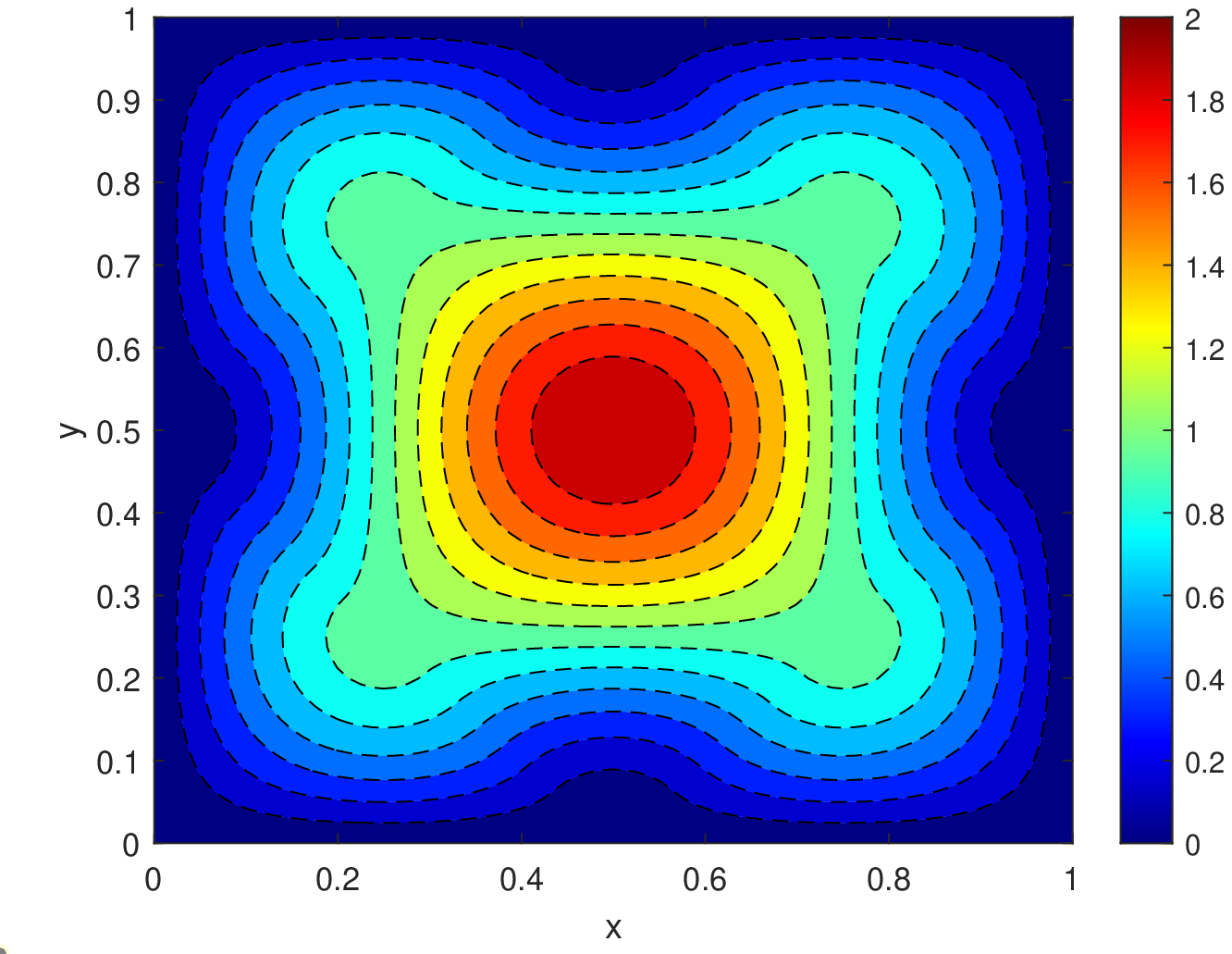}
\end{center}
\caption{Example 6 case (b). The numerical solution by FE-FSM on mesh $N=80$. Left: the 3D plot of numerical solution $\phi$; Right: the contour plot for $\phi$ .}\label{fige62sc}
\end{figure}

For case (b), Figure \ref{fige62sc} shows the surface and contour of numerical solution by FE-FSM.
The numerical errors and orders of convergence
are listed in Table \ref{tab6b}. Due to the non-smoothness of the exact solution,
we can only achieve second order accuracy. Again, the same behavior as the previous
examples can be observe for these four methods. 

\begin{table}
\caption{Example 6 case (a). Comparison of the four methods: The errors of the numerical solution, the accuracy obtained and the number of iterations for convergence}\label{tab6a}
	\begin{center}
\resizebox{\textwidth}{55mm}{
		\begin{tabular}{|c|cccccc|cccccc|}
			\hline
\multicolumn{7}{|c|}{FE-Jacobi $\gamma=0.1$}&\multicolumn{6}{|c|}{FE-Jacobi $\gamma=0.1$ with hybrid strategy}\\ \hline
N & $L_{1}$&order&  $L_{\infty}$&order&iter &time& $L_{1}$&order&  $L_{\infty}$&order&iter&time\\ \hline
40 &7.46e-08 &- &3.55e-07 &- &1519 &1.3384&7.46e-08 &- &3.54e-07 &- &1526 &1.7981\\
80 &3.41e-09 &4.45 &1.49e-08 &4.56 &1974 &7.7447&3.41e-09 &4.45 &1.49e-08 &4.56 &1974 &8.0321\\
160 &- &- &- &- &- &-&- &- &- &- &- &-\\
320 &- &- &- &- &- &- &- &- &- &- &- &-\\ \hline
\multicolumn{7}{|c|}{FE-FSM $\gamma=1$}&\multicolumn{6}{|c|}{FE-FSM $\gamma=1$ with hybrid strategy}\\ \hline
$N$ & $L_{1}$~error&order& $L_{\infty}$~error&order&iter&time& $L_{1}$&order&  $L_{\infty}$&order&iter&time\\ \hline
40 &7.46e-08 &- &3.54e-07 &- &216 &0.29623&7.46e-08 &- &3.54e-07 &- &216 &0.15314\\
80 &3.41e-09 &4.45 &1.49e-08 &4.56 &244 &0.8919&3.41e-09 &4.45 &1.49e-08 &4.56 &244 &0.53176\\
160 &1.09e-10 &4.96 &4.67e-10 &5.00 &288 &4.4203&1.09e-10 &4.96 &4.67e-10 &5.00 &288 &2.1206\\
320 &3.46e-12 &4.98 &1.45e-11 &5.00 &372 &24.2263&3.46e-12 &4.98 &1.45e-11 &5.00 &372 &9.6624\\ \hline
\multicolumn{7}{|c|}{RK-Jacobi $\gamma=1$}&\multicolumn{6}{|c|}{RK-Jacobi $\gamma=1$ with hybrid strategy}\\ \hline
$N$ & $L_{1}$~error&order& $L_{\infty}$~error&order&iter&time& $L_{1}$&order&  $L_{\infty}$&order&iter&time\\ \hline
 40 &7.46e-08 &- &3.56e-07 &- &417 &0.34677&7.46e-08 &- &3.56e-07 &- &417 &0.33492\\
80 &3.41e-09 &4.45 &1.49e-08 &4.57 &501 &1.9144&3.41e-09 &4.45 &1.49e-08 &4.57 &501 &1.2738\\
160 &1.09e-10 &4.96 &4.67e-10 &5.00 &771 &12.7353&1.09e-10 &4.96 &4.67e-10 &5.00 &771 &6.5294\\
320 &3.49e-12 &4.97 &6.73e-11 &2.79 &1392 &97.0302&3.49e-12 &4.97 &6.74e-11 &2.79 &1392 &45.1382\\ \hline
\multicolumn{7}{|c|}{RK-FSM $\gamma=1$}&\multicolumn{6}{|c|}{RK-FSM $\gamma=1$ with hybrid strategy}\\ \hline
$N$ & $L_{1}$~error&order& $L_{\infty}$~error&order&iter&time& $L_{1}$&order&  $L_{\infty}$&order&iter&time\\ \hline
40 &7.46e-08 &- &3.55e-07 &- &324 &0.37401&7.46e-08 &- &3.56e-07 &- &324 &0.20401\\
80 &3.41e-09 &4.45 &1.49e-08 &4.56 &372 &1.4448&3.41e-09 &4.45 &1.50e-08 &4.57 &372 &0.66739\\
160 &1.09e-10 &4.96 &4.67e-10 &5.00 &468 &7.4771&1.09e-10 &4.96 &4.67e-10 &5.00 &468 &2.5455\\
320 &3.46e-12 &4.98 &1.45e-11 &5.00 &636 &41.7842&3.47e-12 &4.98 &1.45e-11 &5.00 &636 &11.4085\\ \hline
\end{tabular}}
\end{center}
\end{table}
\begin{table}
\caption{Example 6 case (b). Comparison of the four methods: The errors of the numerical solution, the accuracy obtained and the number of iterations for convergence}\label{tab6b}
	\begin{center}
\resizebox{\textwidth}{55mm}{
		\begin{tabular}{|c|cccccc|cccccc|}
			\hline
\multicolumn{7}{|c|}{FE-Jacobi $\gamma=0.1$}&\multicolumn{6}{|c|}{FE-Jacobi $\gamma=0.1$ with hybrid strategy}\\ \hline
N & $L_{1}$&order&  $L_{\infty}$&order&iter &time& $L_{1}$&order&  $L_{\infty}$&order&iter&time\\ \hline
40 &2.08e-04 &- &1.39e-03 &- &1528 &1.325&2.13e-04 &- &1.37e-03 &- &1527 &1.6133\\
80 &6.69e-05 &1.63 &5.75e-04 &1.27 &1907 &7.5108&6.70e-05 &1.67 &5.74e-04 &1.26 &1905 &7.2127\\
160 &1.82e-05 &1.87 &2.03e-04 &1.50 &2992 &49.5382&1.82e-05 &1.87 &2.03e-04 &1.49 &2992 &43.7957\\
320 &- &- &- &- &- &- &- &- &- &- &- &-\\ \hline
\multicolumn{7}{|c|}{FE-FSM $\gamma=1$}&\multicolumn{6}{|c|}{FE-FSM $\gamma=1$ with hybrid strategy}\\ \hline
$N$ & $L_{1}$~error&order& $L_{\infty}$~error&order&iter&time& $L_{1}$&order&  $L_{\infty}$&order&iter&time\\ \hline
40 &2.08e-04 &- &1.41e-03 &- &212 &0.23206&2.14e-04 &- &1.40e-03 &- &212 &0.19333\\
80 &6.72e-05 &1.63 &5.79e-04 &1.29 &240 &0.88731&6.73e-05 &1.66 &5.79e-04 &1.27 &240 &0.56582\\
160 &1.83e-05 &1.87 &2.04e-04 &1.50 &284 &4.3167&1.83e-05 &1.87 &2.04e-04 &1.50 &284 &1.8842\\
320 &4.78e-06 &1.93 &6.51e-05 &1.64 &360 &22.9094&4.78e-06 &1.93 &6.51e-05 &1.64 &360 &8.9642\\ \hline
\multicolumn{7}{|c|}{RK-Jacobi $\gamma=1$}&\multicolumn{6}{|c|}{RK-Jacobi $\gamma=1$ with hybrid strategy}\\ \hline
$N$ & $L_{1}$~error&order& $L_{\infty}$~error&order&iter&time& $L_{1}$&order&  $L_{\infty}$&order&iter&time\\ \hline
 40 &2.08e-04 &- &1.38e-03 &- &402 &0.38608&2.13e-04 &- &1.37e-03 &- &396 &0.37579\\
80 &6.69e-05 &1.63 &5.74e-04 &1.27 &456 &1.7332&6.70e-05 &1.67 &5.75e-04 &1.25 &456 &1.0565\\
160 &1.82e-05 &1.87 &2.03e-04 &1.49 &756 &12.0887&1.82e-05 &1.87 &2.03e-04 &1.50 &756 &6.0166\\
320 &4.76e-06 &1.93 &6.49e-05 &1.64 &1380 &94.8398&4.76e-06 &1.93 &6.49e-05 &1.64 &1380 &43.1701\\ \hline
\multicolumn{7}{|c|}{RK-FSM $\gamma=1$}&\multicolumn{6}{|c|}{RK-FSM $\gamma=1$ with hybrid strategy}\\ \hline
$N$ & $L_{1}$~error&order& $L_{\infty}$~error&order&iter&time& $L_{1}$&order&  $L_{\infty}$&order&iter&time\\ \hline
40 &2.08e-04 &- &1.41e-03 &- &324 &0.36728&2.14e-04 &- &1.40e-03 &- &324 &0.22781\\
80 &6.72e-05 &1.63 &5.79e-04 &1.29 &372 &1.4231&6.73e-05 &1.66 &5.79e-04 &1.27 &372 &0.59763\\
160 &1.83e-05 &1.87 &2.04e-04 &1.50 &456 &7.1213&1.83e-05 &1.87 &2.04e-04 &1.50 &456 &2.1953\\
320 &4.78e-06 &1.93 &6.51e-05 &1.64 &612 &39.3192&4.78e-06 &1.93 &6.51e-05 &1.64 &612 &10.2167\\ \hline
\end{tabular}}
\end{center}
\end{table}

\noindent \textbf{Example 7.} We solve the Eikonal equation with
\begin{equation*}
 f(x,y)=2\sqrt{y^2(1-x^2)^2+x^2(1-y^2)^2}.
\end{equation*}
 The computational domain is $\Omega=[-1,1]^2$, and the inflow boundary is the whole outside boundary of the box $[-1,1]^2$, namely $\Gamma=\{(x,y)\mid|x|=1~or~|y|=1\}.$ The boundary condition $\phi(x,y)=0$ is prescribed on $\Gamma$, and an additional boundary condition $\phi(0,0)=1$ is also prescribed at the center of domain. The exact solutions is given by
 \begin{equation*}
\phi(x,y)=(1-x^2)(1-y^2).
\end{equation*}
The exact values are arranged in a small box with length $3h$ around the point $(0,0)$. The errors are measured on the whole domain.

Since the exact solution is a bi-quadratic polynomial, the proposed fifth order HWENO method can recover the exact solution, and the round-off errors can be observed, as shown in Table \ref{tab7}. We can also observe that the fast sweeping technique can improve the convergence of the Jacobi scheme. On the same refined mesh, the RK-FSM only takes about $40\%$ CPU time of the RK-Jacobi scheme, and the FE-FSM costs even less CPU time than RK-FSM. The numerical results after the hybrid strategy is used can be seen on the right side of Table \ref{tab7}, which suggests that the hybrid strategy can save about $80\%$ of the CPU time on the refined mesh.

\begin{table}
\caption{Example 7. Comparison of the four methods: The errors of the numerical solution, the accuracy obtained and the number of iterations for convergence}\label{tab7}
	\begin{center}
\resizebox{\textwidth}{55mm}{
		\begin{tabular}{|c|cccccc|cccccc|}
			\hline
\multicolumn{7}{|c|}{FE-Jacobi $\gamma=0.1$}&\multicolumn{6}{|c|}{FE-Jacobi $\gamma=0.1$ with hybrid strategy}\\ \hline
N & $L_{1}$&order&  $L_{\infty}$&order&iter &time& $L_{1}$&order&  $L_{\infty}$&order&iter&time\\ \hline
40 &1.05e-14 &- &3.15e-13 &- &1547 &1.4301 &1.05e-14 &-  &3.15e-13 &- &1547 &0.5632\\
80 &9.90e-15 &9.41 &8.99e-13 &-1.51 &1850 &6.9849 &9.73e-15 &1.19  &9.00e-13 &-1.51 &1850 &1.8203\\
160 &9.85e-15 &7.11 &1.58e-12 &-8.18 &3016 &51.576 &8.71e-15 &1.59  &1.42e-12 &-6.64 &3019 &10.1809\\
320 &1.04e-14 &-8.33 &4.06e-12 &-1.35 &5559 & 428.72 &9.51e-15 &-1.27  &4.08e-12 &-1.51 &5559 &89.2658\\ \hline
\multicolumn{7}{|c|}{FE-FSM $\gamma=1$}&\multicolumn{6}{|c|}{FE-FSM $\gamma=1$ with hybrid strategy}\\ \hline
$N$ & $L_{1}$~error&order& $L_{\infty}$~error&order&iter&time& $L_{1}$&order&  $L_{\infty}$&order&iter&time\\ \hline
40 &4.69e-16 &- &8.32e-15 &- &212 &0.2448&3.49e-16 &-  &7.54e-15 &- &212 &0.1366\\
80 &8.31e-16 &-8.24 &3.48e-14 &-2.06 &236 & 0.8822&3.73e-16 &-9.41  &3.68e-14 &-2.28 &236 &0.3952\\
160 &1.51e-15 &-8.63 &1.21e-13 &-1.80 &284 &4.5492 &4.01e-16 &-1.03  &1.17e-13 &-1.67 &284 &0.8534\\
320 &2.64e-15 &-8.06 &1.10e-13 &1.38 &376 &25.933&3.62e-16 &1.48 &1.16e-13 &8.20 &376 &4.2301\\ \hline
\multicolumn{7}{|c|}{RK-Jacobi $\gamma=1$}&\multicolumn{6}{|c|}{RK-Jacobi $\gamma=1$ with hybrid strategy}\\ \hline
$N$ & $L_{1}$~error&order& $L_{\infty}$~error&order&iter&time& $L_{1}$&order&  $L_{\infty}$&order&iter&time\\ \hline
 40 &1.12e-15 &- &4.96e-14 &- &393 &0.4107&1.00e-15 &-  &3.49e-14 &- &402 &0.2235\\
80 &6.93e-15 &-2.62 &6.47e-13 &-3.70 &477 &1.8100&7.96e-15 &-2.99  &6.73e-13 &-4.26 &474 &0.64002\\
160 &2.62e-15 &1.40 &1.21e-12 &-9.07 &741 &12.404&3.03e-15 &1.39 &1.22e-12 &-8.65 &741 &2.6417\\
320 &3.37e-15 &-3.66 &5.18e-12 &-2.09 &1356 &101.06 &4.50e-15 &-5.68  &5.20e-12 &-2.08 &1356 &20.7209\\ \hline
\multicolumn{7}{|c|}{RK-FSM $\gamma=1$}&\multicolumn{6}{|c|}{RK-FSM $\gamma=1$ with hybrid strategy}\\ \hline
$N$ & $L_{1}$~error&order& $L_{\infty}$~error&order&iter&time& $L_{1}$&order&  $L_{\infty}$&order&iter&time\\ \hline
40 &4.32e-16 &- &8.10e-15 &- &324 &0.3867&3.22e-16 &-  &8.54e-15 &- &324 &0.1814\\
80 &8.51e-16 &-9.77 &3.87e-14 &-2.25 &360 &1.3585&3.91e-16 &-2.78 &4.09e-14 &-2.26 &360 &0.3997\\
160 &1.72e-15 &-1.01 &2.24e-13 &-2.53 &444 &7.1011&6.35e-16 &-6.98  &2.18e-13 &-2.41 &444 &1.2264\\
320 &2.76e-15 &-6.81 &2.14e-13 &7.00 &612 &41.072&5.05e-16 &3.29  &2.27e-13 &-5.53 &612 &6.576\\ \hline
\end{tabular}}
\end{center}
\end{table}

\noindent \textbf{Example 8} The travel-time problem in elastic wave propagation is considered in this example. The
quasi-P and the quasi-SV slowness surfaces are defined as follows \cite{qianchengosher}
\begin{equation*}
  c_{1}\phi_{x}^{4}+c_{2}\phi_{x}^{2}\phi_{y}^{2}+c_{3}\phi_{y}^{4}+c_{4}\phi_{x}^{2}+c_{5}\phi_{y}^{2}+1=0,
\end{equation*}
where
\begin{equation*}
\begin{split}
  c_{1}=a_{11}a_{44}, \quad c_{2}=a_{11}a_{33}+a_{44}^{2}-(a_{13}+a_{44})^{2},&  \\
    c_{3}=a_{33}a_{44},\quad c_{4}=-(a_{11}+a_{44}),\quad c_{5}=-(a_{33}+a_{44}),&
\end{split}
\end{equation*}
in which $a_{i,j}$ are given elastic parameters. The quasi-P wave Eikonal equation is
\begin{equation*}
  \sqrt{-\frac{1}{2}(c_{4}\phi_{x}^{2}+c_{5}\phi_{y}^{2})+\sqrt{\frac{1}{4}(c_{4}\phi_{x}^{2}+c_{5}\phi_{y}^{2})^2-(c_{1}\phi_{x}^{4}+c_{2}\phi_{x}^{2}\phi_{y}^{2}+c_{3}\phi_{y}^{4})}}=1,
\end{equation*}
which is a convex HJ equation, and the elastic parameters are taken to be
\begin{equation*}
  a_{11}=15.0638,~ a_{33}=10.8373,~a_{13}=1.6381,~a_{44}=3.1258.
\end{equation*}
The corresponding quasi-SV wave Eikonal equation is given by
 \begin{equation*}
  \sqrt{-\frac{1}{2}(c_{4}\phi_{x}^{2}+c_{5}\phi_{y}^{2})-\sqrt{\frac{1}{4}(c_{4}\phi_{x}^{2}+c_{5}\phi_{y}^{2})^2-(c_{1}\phi_{x}^{4}+c_{2}\phi_{x}^{2}\phi_{y}^{2}+c_{3}\phi_{y}^{4})}}=1,
\end{equation*}
which is a nonconvex HJ equation, and the elastic parameters are taken to be
\begin{equation*}
  a_{11}=15.90,~ a_{33}=6.21,~a_{13}=4.82,~a_{44}=4.00.
\end{equation*}
The computational domain is set as $\Omega=[-1,1]^2$, and the inflow boundary is given by $\Gamma={(0,0)}$. Exact values are assigned in a small box with length $0.3$ around the source point. Because these Hamiltonians are in complicated forms, we use the LF numerical Hamiltonian for both equations. In addition, since we only know the numerical solution of $\phi$, the ``exact solution'' of $u$ and $v$ on \emph{Category I} will be obtained by fifth order WENO reconstruction.

\begin{figure}[!h]
\begin{center}
	\includegraphics[width=6.cm]{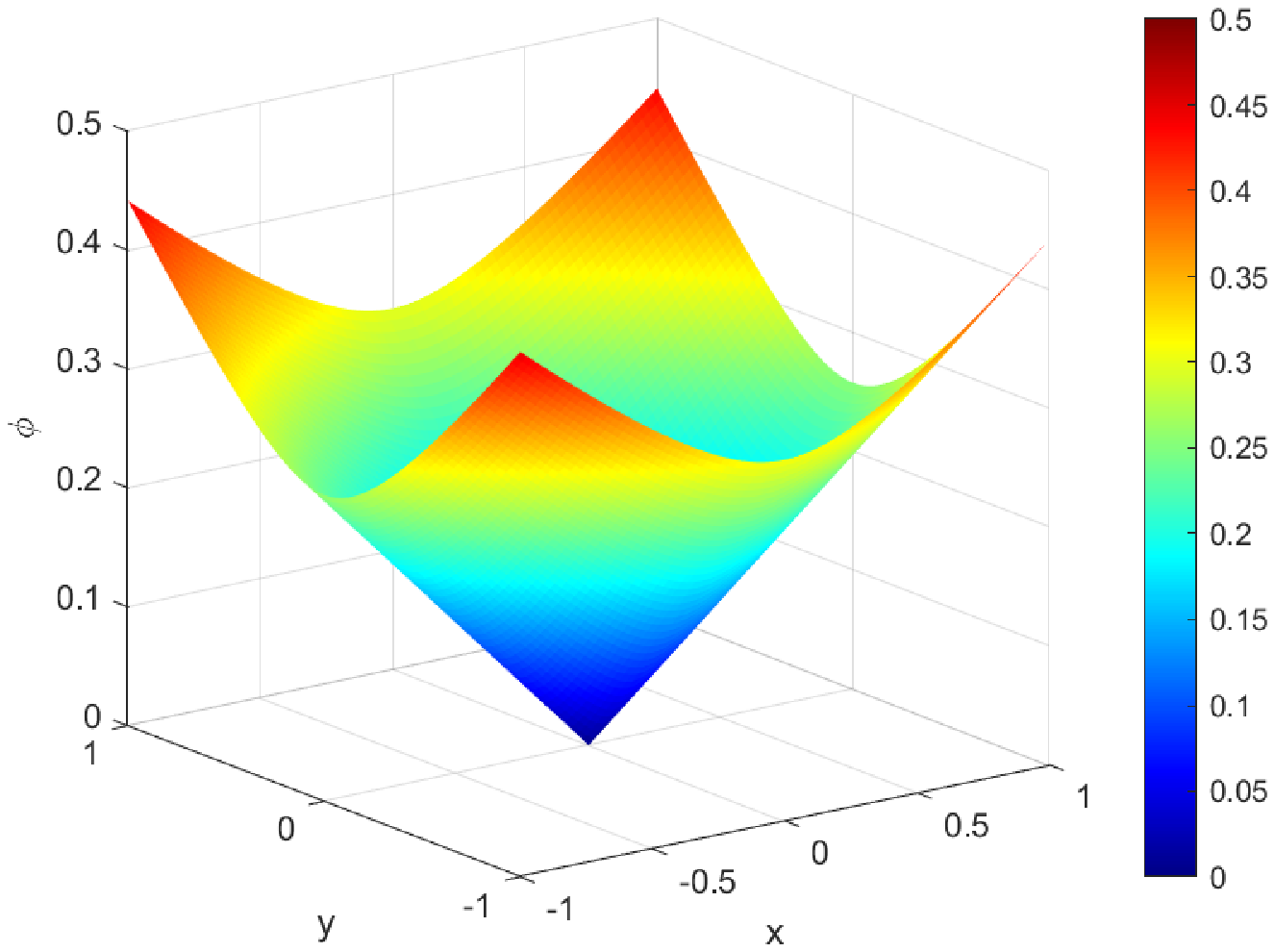}
	\includegraphics[width=6.cm]{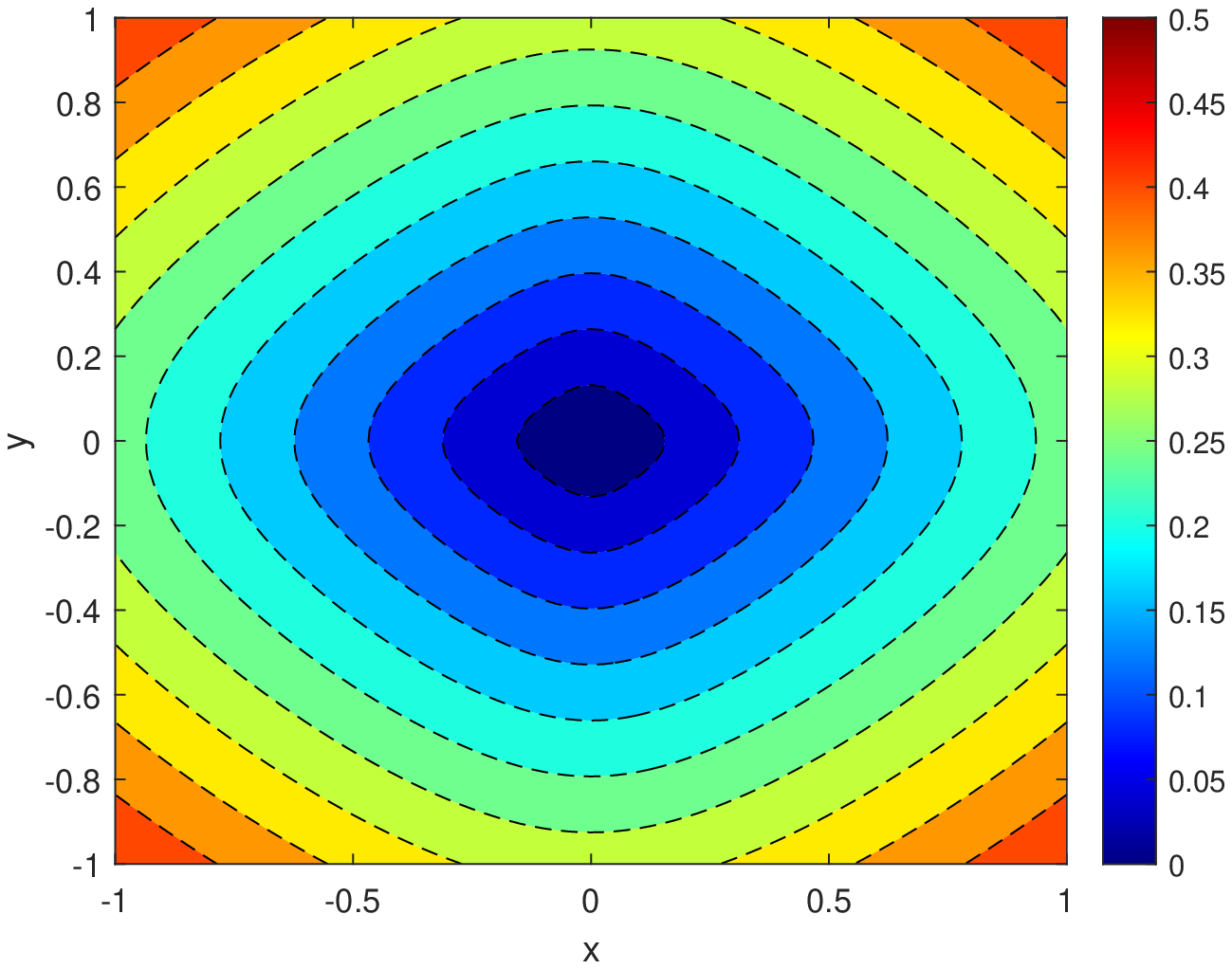}
\end{center}
\caption{Example 8 P-wave. The numerical solution by FE-FSM on mesh $N=80$. Left: the 3D plot of numerical solution $\phi$; Right: the contour plot for $\phi$ .}\label{fige8psc}
\end{figure}

For the P-wave equation, the surface and contour of numerical solution by FE-FSM are shown in Figure \ref{fige8psc}. The numerical errors and orders of convergence are presented in Table
\ref{tab8p} for four methods. We can observe that the FE-Jacobi method requires the smaller CFL number of value $0.1$. When the third RK time discretization is used, the RK-Jacobi scheme can take a larger CFL number of value $1$. Also, the number of iterations of RK-Jacobi method is smaller than FE-Jacobi method on the same mesh. The fast sweeping technique can improve the convergence of the Jacobi scheme. On the same refined mesh, we can see that the RK-FSM only takes about $50\%$ CPU time of the RK-Jacobi scheme. Furthermore, the FE-FSM costs even less CPU time than RK-FSM. The numerical results after the hybrid strategy is used can be seen on the right side of Table \ref{tab8p}, which suggests that the hybrid strategy can save $50\%$ of the CPU time on the refined mesh.

\begin{figure}[!h]
\begin{center}
	\includegraphics[width=6.cm]{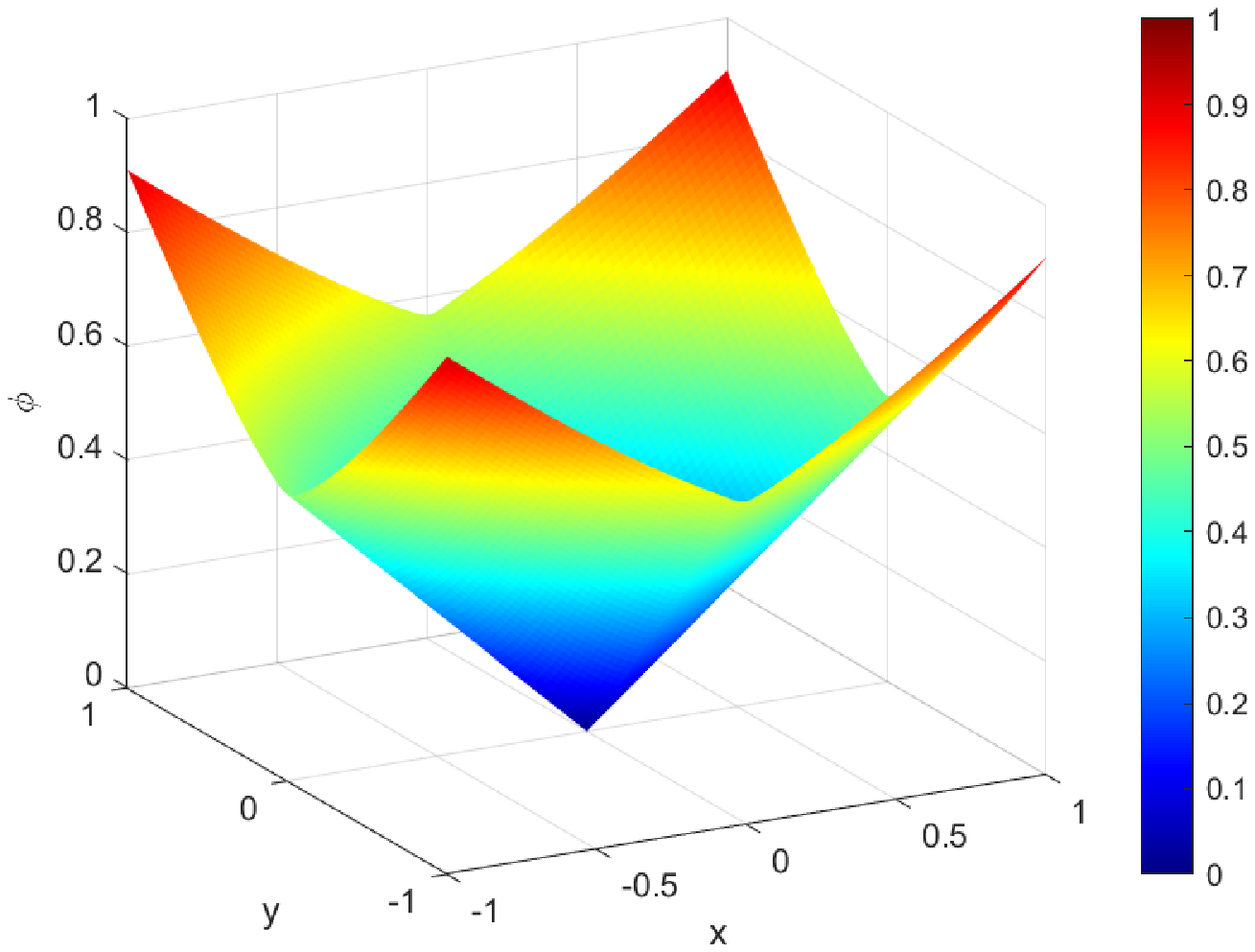}
	\includegraphics[width=6.cm]{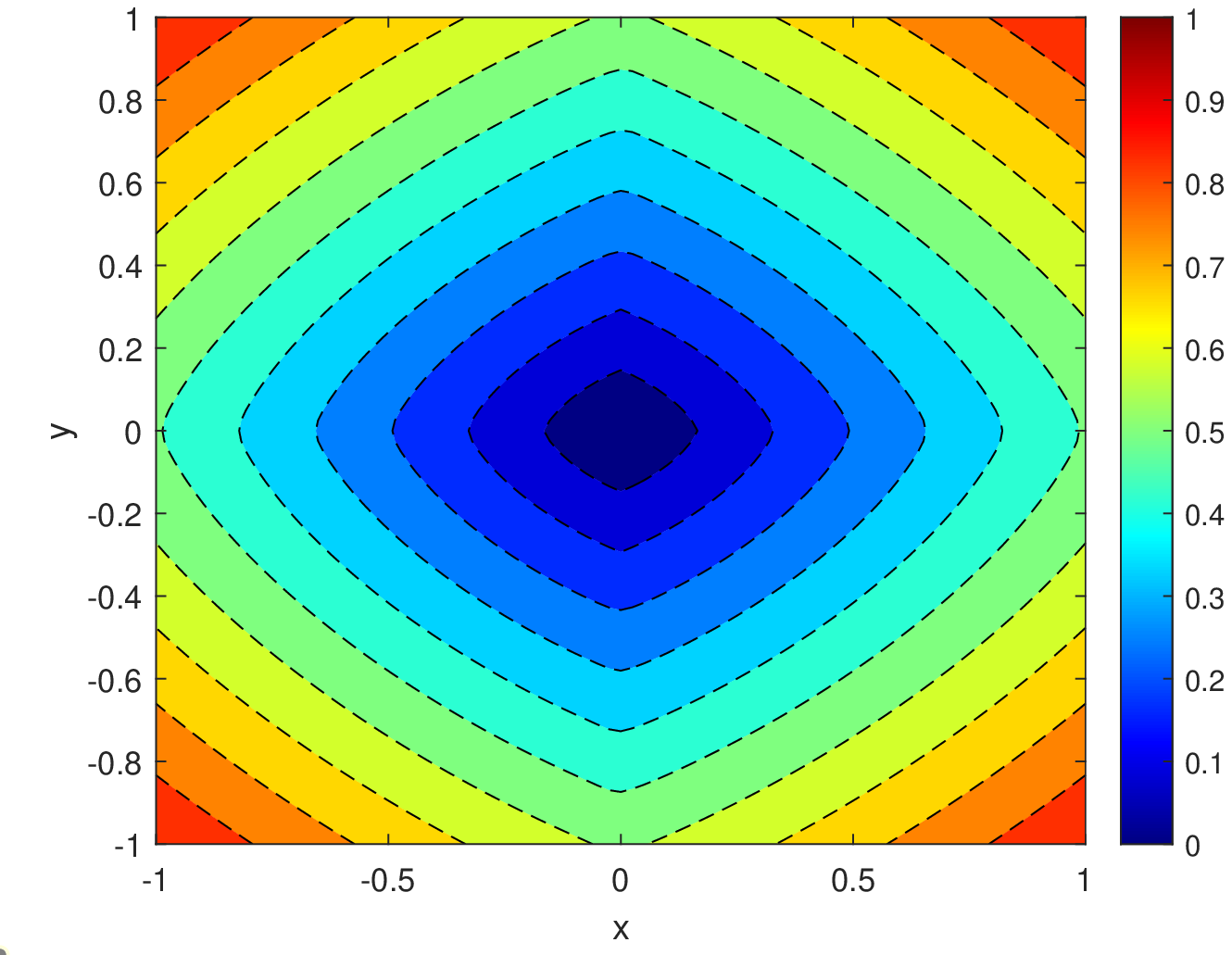}
\end{center}
\caption{Example 8 SV-wave. The numerical solution by FE-FSM on mesh $N=80$. Left: the 3D plot of numerical solution $\phi$; Right: the contour plot for $\phi$ .}\label{fige8svsc}
\end{figure}

For the SV-wave equation, Figure \ref{fige8svsc} shows the surface and contour of numerical solution by FE-FSM. The numerical errors and orders of convergence are listed in Table
\ref{tab8sv} for four methods. Again, the FE-FSM performs the best out of these four methods, and the hybrid strategy can further reduce the computational cost.

\begin{table}
\caption{Example 8 P-wave. Comparison of the four methods: The errors of the numerical solution, the accuracy obtained and the number of iterations for convergence}\label{tab8p}
	\begin{center}
\resizebox{\textwidth}{55mm}{
		\begin{tabular}{|c|cccccc|cccccc|}
			\hline
\multicolumn{7}{|c|}{FE-Jacobi $\gamma=0.1$}&\multicolumn{6}{|c|}{FE-Jacobi $\gamma=0.1$ with hybrid strategy}\\ \hline
N & $L_{1}$&order&  $L_{\infty}$&order&iter &time& $L_{1}$&order&  $L_{\infty}$&order&iter&time\\ \hline
40 &4.77e-06 &- &3.95e-05 &- &1469 &2.1388&4.77e-06 &- &3.95e-05 &- &1463 &2.0531\\
80 &2.07e-07 &4.52 &2.25e-06 &4.13 &1781 &9.9557&2.07e-07 &4.52 &2.25e-06 &4.13 &1938 &9.4128\\
160 &7.05e-09 &4.87 &8.09e-08 &4.79 &2778 &61.2665 &7.05e-09 &4.87 &8.09e-08 &4.79 &2941 &57.1844\\
320 &2.29e-10 &4.93 &2.61e-09 &4.95 &4785 &435.3249 &2.29e-10 &4.93 &2.61e-09 &4.95 &4856 &391.72\\ \hline
\multicolumn{7}{|c|}{FE-FSM $\gamma=1$}&\multicolumn{6}{|c|}{FE-FSM $\gamma=1$ with hybrid strategy}\\ \hline
$N$ & $L_{1}$~error&order& $L_{\infty}$~error&order&iter&time& $L_{1}$&order&  $L_{\infty}$&order&iter&time\\ \hline
40 &4.77e-06 &- &3.95e-05 &- &184 &0.24145 &4.77e-06 &- &3.95e-05 &- &184 &0.1427\\
80 &2.07e-07 &4.52 &2.25e-06 &4.13 &216 &1.1029&2.07e-07 &4.52 &2.25e-06 &4.13 &212 &0.54758\\
160 &7.05e-09 &4.87 &8.09e-08 &4.79 &272 &5.6603&7.05e-09 &4.87 &8.09e-08 &4.79 &272 &2.8627\\
320 &2.29e-10 &4.93 &2.61e-09 &4.95 &376 &30.9783&2.29e-10 &4.93 &2.61e-09 &4.95 &400 &15.4856\\ \hline
\multicolumn{7}{|c|}{RK-Jacobi $\gamma=1$}&\multicolumn{6}{|c|}{RK-Jacobi $\gamma=1$ with hybrid strategy}\\ \hline
$N$ & $L_{1}$~error&order& $L_{\infty}$~error&order&iter&time& $L_{1}$&order&  $L_{\infty}$&order&iter&time\\ \hline
 40 &4.77e-06 &- &3.95e-05 &- &429 &0.5894&4.78e-06 &- &3.95e-05 &- &429 &0.42699\\
80 &2.07e-07 &4.52 &2.25e-06 &4.13 &495 &2.6246&2.07e-07 &4.52 &2.25e-06 &4.13 &498 &1.7419\\
160 &7.05e-09 &4.87 &8.09e-08 &4.79 &738 &15.6503&7.05e-09 &4.87 &8.09e-08 &4.79 &723 &7.9175\\
320 &2.29e-10 &4.93 &2.61e-09 &4.95 &1296 &112.4767&2.29e-10 &4.93 &2.61e-09 &4.95 &1296 &56.0106\\ \hline
\multicolumn{7}{|c|}{RK-FSM $\gamma=1$}&\multicolumn{6}{|c|}{RK-FSM $\gamma=1$ with hybrid strategy}\\ \hline
$N$ & $L_{1}$~error&order& $L_{\infty}$~error&order&iter&time& $L_{1}$&order&  $L_{\infty}$&order&iter&time\\ \hline
40 &4.77e-06 &- &3.95e-05 &- &300 &0.3865&4.77e-06 &- &3.95e-05 &- &300 &0.2366\\
80 &2.07e-07 &4.52 &2.25e-06 &4.13 &336 &1.6854&2.07e-07 &4.52 &2.25e-06 &4.13 &336 &0.9660\\
160 &7.05e-09 &4.87 &8.09e-08 &4.79 &420 &8.638&7.05e-09 &4.87 &8.09e-08 &4.79 &420 &4.2972\\
320 &2.29e-10 &4.93 &2.61e-09 &4.95 &612 &48.9773&2.29e-10 &4.93 &2.61e-09 &4.95 &612 &24.5722\\ \hline
\end{tabular}}
\end{center}
\end{table}
\begin{table}
\caption{Example 8 SV-wave. Comparison of the four methods: The errors of the numerical solution, the accuracy obtained and the number of iterations for convergence}\label{tab8sv}
	\begin{center}
\resizebox{\textwidth}{47mm}{
		\begin{tabular}{|c|cccccc|cccccc|}
			\hline
\multicolumn{7}{|c|}{FE-Jacobi $\gamma=0.1$}&\multicolumn{6}{|c|}{FE-Jacobi $\gamma=0.1$ with hybrid strategy}\\ \hline
N & $L_{1}$&order&  $L_{\infty}$&order&iter &time& $L_{1}$&order&  $L_{\infty}$&order&iter&time\\ \hline
80 &1.00e-06 &- &1.65e-05 &- &1891 &9.7947&1.01e-06 &- &1.65e-05 &- &1967 &9.3438\\
160 &2.58e-08 &5.28 &1.03e-06 &4.00 &2865 &59.192 &2.59e-08 &5.28 &1.03e-06 &3.99 &2974 &55.317\\
320 &1.08e-10 &7.89 &1.15e-08 &6.48 &4995 &415.46&1.09e-10 &7.89 &1.15e-08 &6.49 &4999 &384.21\\ \hline
\multicolumn{7}{|c|}{FE-FSM $\gamma=1$}&\multicolumn{6}{|c|}{FE-FSM $\gamma=1$ with hybrid strategy}\\ \hline
$N$ & $L_{1}$~error&order& $L_{\infty}$~error&order&iter&time& $L_{1}$&order&  $L_{\infty}$&order&iter&time\\ \hline
80&9.63e-07 &- &1.60e-05 &-&284&1.5094 &8.61e-07 &- &2.19e-05 &- &272 &0.75179\\
160 &1.95e-08 &5.61 &8.97e-07 &4.15 &368&7.6795 &2.25e-08 &5.25 &1.44e-06 &3.92 &312 &2.959\\
320 &6.37e-11 &8.26 &9.12e-09 &6.61 &448 &35.524 &7.85e-11 &8.16 &7.90e-09 &7.51 &452 &17.488\\ \hline
\multicolumn{7}{|c|}{RK-Jacobi $\gamma=1$}&\multicolumn{6}{|c|}{RK-Jacobi $\gamma=1$ with hybrid strategy}\\ \hline
$N$ & $L_{1}$~error&order& $L_{\infty}$~error&order&iter&time& $L_{1}$&order&  $L_{\infty}$&order&iter&time\\ \hline
 80 &1.00e-06 &- &1.65e-05 &- &567 &3.1826&9.27e-07 &-Inf &1.60e-05 &- &570 &2.2798\\
160 &2.58e-08 &5.28 &1.03e-06 &4.00 &855 &18.276 &2.191-08 &5.40 &9.80e-07 &4.03 &846 &8.9459\\
320 &1.08e-10 &7.89 &1.15e-08 &6.48 &1476 &125.21&9.06e-11 &7.91 &1.12e-08 &6.44 &1467 &60.973\\ \hline
\multicolumn{7}{|c|}{RK-FSM $\gamma=1$}&\multicolumn{6}{|c|}{RK-FSM $\gamma=1$ with hybrid strategy}\\ \hline
$N$ & $L_{1}$~error&order& $L_{\infty}$~error&order&iter&time& $L_{1}$&order&  $L_{\infty}$&order&iter&time\\ \hline
80 &9.51e-07 &-&1.50e-05 &- &360 &1.8518&9.17e-07 &- &1.65e-05 &- &552 &1.7596\\
160 &2.10e-08 &5.49 &9.78e-07 &3.9473 &480 &9.5202&2.66e-08 &5.10 &1.21e-06 &3.76 &480 &4.6605\\
320 &6.60e-11 &8.31 &1.00e-08 &6.604 &720 &57.929&8.52e-11 &8.28 &1.08e-08 &6.79 &720 &27.38\\ \hline
\end{tabular}}
\end{center}
\end{table}

\section{Conclusion Remark}\label{sec:conclusion}

In this paper, we design a fifth order HWENO fixed-point fast sweeping method for solving the static HJ equations, by combining the fixed-point iteration with the fast sweeping strategy in a novel HWENO framework. The numerical study suggests that the fast sweeping technique can greatly improve the stability of the high-order spatial scheme. We presented a large number of numerical experiments to test four different methods, including FE method and RK time marching method, and the methods combined with fast sweeping technique. Numerical results show that the FE time-marching method with fast sweeping technique is the most effective method to solve the static HJ equations. In addition, a hybrid strategy which combines both linear and HWENO reconstruction is also proposed and tested, which yields additional savings in computational time.


\end{document}